\documentclass[12pt]{amsart}
\usepackage{latexsym, amssymb, amsfonts, amscd}

\theoremstyle{plain}
\newtheorem{Thm}{Theorem}[subsection]
\newtheorem{Cor}[Thm]{Corollary}

\newtheorem{Prop}[Thm]{Proposition}
\newtheorem{Lem}[Thm]{Lemma}

\newtheorem{Cl}[Thm]{Claim}

\newtheorem{Thm'}{Theorem}[section]
\newtheorem{Cor'}[Thm']{Corollary}
\newtheorem{Prop'}[Thm']{Proposition}
\newtheorem{Lem'}[Thm']{Lemma}
\newtheorem{Cl'}[Thm']{Claim}

\theoremstyle{definition}
\newtheorem{Def}[Thm]{Definition}
\newtheorem{Rem}[Thm]{Remark}

\newtheorem{Emp}[Thm]{}
\newtheorem{Ex}[Thm]{Example}

\newtheorem{Not}[Thm]{Notation}

\newtheorem{Def'}[Thm']{Definition}
\newtheorem{Rem'}[Thm']{Remark}
\newtheorem{Rem1'}[Thm']{Remarks}
\newtheorem{Emp'}[Thm']{}
\newtheorem{Ex'}[Thm']{Example}
\newtheorem{Exs'}[Thm']{Examples}
\newtheorem{Con'}[Thm']{Construction}
\newtheorem{Not'}[Thm']{Notation}
\newtheorem{Q'}[Thm']{Question}

\numberwithin{equation}{section}

\newcommand{\qlbar}{\overline{\B{Q}_l}}

\newcommand{\om}{\omega}
\newcommand{\diam}{\diamondsuit}
\newcommand{\La}{\Lambda}
\newcommand{\lan}{\left\langle}
\newcommand{\ran}{\right\rangle}
\newcommand{\pp}{\boxtimes}\newcommand{\ov}{\overline}
\newcommand{\un}{\underline}

\newcommand{\fq}{\B{F}_q}

\newcommand{\B}[1]{\mathbb#1}
\newcommand{\cal}[1]{\mathcal{#1}}
\newcommand{\form}[1]{(\ref{Eq:#1})}
\newcommand{\C}[1]{\cal#1}

\newcommand{\isom}{\overset {\thicksim}{\to}}

\newcommand{\lra}{\longrightarrow}
\newcommand{\lla}{\longleftarrow}
\newcommand{\hra}{\hookrightarrow}
\newcommand{\wt}{\widetilde}
\newcommand{\diez}{\natural}

\newcommand{\Gm}{\Gamma}

\newcommand{\Dt}{\Delta}

\newcommand{\rs}[1]{Section \ref{S:#1}}
\newcommand{\rl}[1]{Lemma \ref{L:#1}}
\newcommand{\rn}[1]{Notation \ref{N:#1}}
\newcommand{\rcl}[1]{Claim \ref{C:#1}}
\newcommand{\rp}[1]{Proposition \ref{P:#1}}
\newcommand{\rr}[1]{Remark \ref{R:#1}}

\newcommand{\re}[1]{\ref{E:#1}}
\newcommand{\rco}[1]{Corollary \ref{C:#1}}

\newcommand{\rt}[1] {Theorem \ref{T:#1}}
\newcommand{\rd}[1]{Definition \ref{D:#1}}
\newcommand{\sm}{\smallsetminus}

\newcommand{\pt}{\operatorname{pt}}

\newcommand{\Ob}{\operatorname{Ob}}
\newcommand{\Mor}{\operatorname{Mor}}

\newcommand{\Spec}{\operatorname{Spec}}

\newcommand{\Tr}{\operatorname{Tr}}
\newcommand{\Fr}{\operatorname{Fr}}

\newcommand{\Id}{\operatorname{Id}}

\newcommand{\Hom}{\operatorname{Hom}}

\newcommand{\RHom}{\operatorname{RHom}}

\begin{document}

\title[Lefschetz--Verdier trace formula]%
{Lefschetz--Verdier trace formula 
and a generalization of a theorem of Fujiwara}
\author{Yakov Varshavsky}
\address{Institute of Mathematics\\
Hebrew University\\
Givat-Ram, Jerusalem,  91904\\
Israel}
\email{vyakov@math.huji.ac.il }

\thanks{The work was supported by 
THE ISRAEL SCIENCE FOUNDATION (Grant No. 555/04)}
\date{November 2005}

\begin{abstract}

The goal of this paper is to generalize a theorem of Fujiwara (formerly 
Deligne's conjecture) to the situation appearing in a joint work \cite{KV} 
with David Kazhdan on the global Langlands correspondence over function fields.  
Moreover, our proof is more elementary than the original one and stays in the 
realm of ordinary algebraic geometry, that is, does not use rigid geometry.
We also give a proof of the Lefschetz--Verdier 
trace formula and of the additivity of filtered trace maps, 
thus making the paper essentially self-contained. 
\end{abstract}
\maketitle

\tableofcontents

\section*{Introduction}
Suppose we are given a diagram
$X\overset{c_1}{\lla}C\overset{c_2}{\lra}X$
of separates schemes of finite type over a separably closed field $k$, which we call a 
{\em correspondence},   
an ``$\ell$-adic sheaf'' $\C{F}\in D^b_{c}(X,\qlbar)$ and a morphism 
$u:c_{2!}c_1^*\C{F}\to\C{F}$.
If $c_1$ is proper, then $u$ gives rise to an endomorphism 
$R\Gm_c(u):R\Gm_c(X,\C{F})\to R\Gm_c(X,\C{F})$. 

When $X$ is proper, the general Lefschetz-Verdier trace formula
\cite[Cor. 4.7]{Il} asserts that the trace $\Tr(R\Gm_c(u))$ equals the sum  
$\sum_{\beta\in\pi_0(Fix(c))}LT_{\beta}(u)$, where 
$Fix(c):=\{y\in C\,|\,c_1(y)=c_2(y)\}$ is the scheme of fixed points of $c$,
and $LT_{\beta}(u)$ is a so-called ``local term'' of $u$ at $\beta$.
This result has two drawbacks: it fails when $X$ is not proper, and the 
``local terms'' are very inexplicit. 

Deligne conjectured that the situation becomes better if $k$ is the algebraic 
closure of a finite field $\fq$, $X\overset{c_1}{\lla}C\overset{c_2}{\lra}X$ and 
$\C{F}$ are defined over
$\fq$,  $c_2$ is quasi-finite, and we twist $c_1$ by a sufficiently 
high power of the geometric Frobenius morphism. 
More precisely, he conjectured that in this case 
the Lefschetz-Verdier trace formula holds also for non-proper $X$'s, 
the scheme of fixed points $Fix(c)$ is finite, and for each $y\in c_1^{-1}(x)\cap c_2^{-1}(x)$, 
the local term $LT_{y}(u)$  equals the trace of the 
endomorphism $u_y:\C{F}_x\to\C{F}_x$, induced by $u$.

The conjecture was first proven by  
Pink \cite{Pi} assuming the resolution of singularities, and then by 
Fujiwara \cite{Fu} unconditionally. 

Theorem of Fujiwara has a fundamental importance for Langlands program. 
For example, it was crucially used by Flicker--Kazhdan, Harris--Taylor and Lafforgue. 

In a joint project \cite{KV} with David Kazhdan on the global Langlands correspondence 
over function fields, we needed a generalization of 
the above result to the case, when $c_1$ is not necessary proper but there exists 
an open subset $U\subset X$ such that $c_1|_{c_1^{-1}(U)}$ is proper,  
$\C{F}$ vanishes on $X\sm U$, and $X\sm U$ is ``locally $c$-invariant''. In this case, 
$u$ still gives rise to an endomorphism 
$R\Gm_c(u)$, and  the main result of this work asserts that the conclusion of 
Deligne's conjecture holds in this case.

The strategy of our proof is very similar to that of \cite{Pi} and \cite{Fu}: first we 
reduce the problem to vanishing of local terms $LT_{\beta}$, then we make 
the correspondence ``contracting'' by twisting it with a 
sufficiently high power of Frobenius, and finally we show the vanishing of local terms 
for ``contracting'' correspondences.  

Our approach differs from that of Fujiwara in two respects. First of all, 
our definition of a ``contracting'' correspondence is much simpler both to define 
and to use. 
Namely, we work with the most naive notion of an ``infinitesimally'' contracting correspondence, 
which has a simple geometric description in terms of a ``deformation to the normal cone''. 
As a result, our bound on the power of Frobenius is sharper and more explicit.

Secondly, to prove a generalization of Deligne's conjecture described above, 
we work ``locally''. More precisely, motivated by Pink \cite{Pi}, to show the 
vanishing of local terms, 
we first show the vanishing of so-called ``trace maps'', from which 
local terms are obtained by integration. Also our definition of 
local terms is more elementary. 


The paper is organized as follows. In the first section we introduce basic objects and 
constructions used later.
More precisely, in  Subsection 1.1 we define correspondences, cohomological correspondences
and basic operations on them. In Subsection 1.2 we introduce ``trace maps'' and 
show that the Lefschetz-Verdier trace formula follows from the commutativity of the 
trace maps with proper push-forwards. In Subsection 1.3 we define 
specialization of correspondences, while in Subsection 1.4 we restrict ourself to a particular case
called ``specialization (or deformation) to the normal cone''. 
We finish this section by Subsection 1.5, in which we introduce locally invariant subschemes
and study their properties.

In the second section we prove our main results, assuming certain functorial properties
of the trace maps. More precisely, in Subsection 2.1 we introduce the notion of correspondences 
``contracting near fixed points'' and show that for such correspondences, the local terms are 
equal to the ``naive local terms''. Our proof uses the additivity of the trace maps, the 
commutativity of the trace maps with specializations, and a theorem of Verdier \cite{Ve} 
asserting that specialization to the normal cone commutes with the restriction to the zero section. 

In Subsection 2.2 we study correspondences over finite fields and show that in our situation
the correspondence become ``contracting near fixed points'' after we twist $c_1$ by a 
sufficiently high power of the geometric Frobenius morphism. Finally, in 
Subsection 2.3 we prove the generalization of Deligne's conjecture, described above.

The third section is devoted to the proof of the theorem of Verdier \cite{Ve}, which 
is crucial for Subsection 2.1. Though our argument is almost identical to the original one, 
it is more detailed.

In the fourth section we prove that the trace maps commute with proper push-forwards and 
specializations. Though similar assertions were proven by Illusie and Fujiwara, respectively,
their results have unnecessary properness assumptions. Moreover, we find Fujiwara's proof 
rather sketchy. To make our proofs more structural, we prove a more general result asserting that  
the trace maps commute with so-called ``cohomological morphisms''.
We would like to observe that our proof reduces to a long sequence of rather 
straight-forward calculations. We believe that  it should be possible to replace 
these calculation by some ``categorical'' argument.

The fifth section is devoted to the proof of the additivity of trace maps used in Subsection 2.1. 
Following Pink, we deduce this property from the additivity of filtered trace maps, 
stated by Illusie in \cite[4.13]{Il}. Though the result is considered well-known among specialists,
we decided to include the proof for completeness. To prove the result, we first recall 
basic properties of filtered derived categories in Subsection 4.1, then we define 
filtered six operations in Subsection 4.2, and finally we show the 
additivity of filtered trace maps in Subsection 4.3.

In the appendix we show that our local terms coincide with those defined by 
Illusie and used by Pink and Fujiwara. Though we do not use this result anywhere 
in the paper, we decided to include it to avoid a confusion.

Finaly we would like to stress that though this paper is written in the language of schemes,  
our arguments apply word-by-word to the case of (compactifiable) algebraic spaces 
(and most likely also to the case Deligne--Mumford stacks).

A shorter exposition of our results appeared in \cite{Va}.

\section*{Acknowledgments}
This work would not be possible without David Kazhdan, who among other things 
explained to me how to define $R\Gm_c(u)$ in the case described above and 
suggested that an analog of Deligne's conjecture should hold in this situation.

Also the author thanks Alexander Beilinson, who gave a reference to \cite{Ve} and whose 
comments and suggestions helped to simplify and improve the exposition. Finally I thank 
Dennis Gaitsgory for his interest and valuable remarks. 

Part of the work was done while the author visited the University of Chicago and 
Northwestern University. The author thanks these institutions for hospitality and 
financial support.

\section*{Notation and conventions}
\begin{Emp'} \label{E:Zar}
For a scheme $X$, we will denote by $X_{red}$ the corresponding reduced scheme.
We will identify closed subsets of $X$ with the corresponding 
closed reduced subschemes. For a closed subscheme $Z\subset X$, denote by 
$\C{I}_Z\subset\C{O}_X$ the sheaf of ideals of $Z$. For a morphism of schemes $f:Y\to X$, 
we denote by $f^{\cdot}:\C{O}_X\to \C{O}_Y$
the pullback map of functions.  For a morphism of schemes $f:Y\to X$ and  
a closed subscheme $Z\subset X$, we denote by $f^{-1}(Z)$ the schematic inverse 
image of $Z$, i.e., the closed subscheme of $X'$ such that 
$\C{I}_{f^{-1}(Z)}=f^{\cdot}(\C{I}_{Z})\cdot\C{O}_{X'}\subset\C{O}_{X'}$.
\end{Emp'}

\begin{Emp'} \label{E:etale}
Throughout the paper all schemes will be  separated and of finite type over 
a fixed separably closed field $k$. We also fix a prime $l$,  invertible in $k$, and a 
commutative ring with identity $\La$, which is either finite and is annihilated by some 
power of $l$, or a  finite extension of $\B{Z}_l$ or $\B{Q}_l$. 

To each scheme $X$ as above, we associate a category $D_{ctf}^b(X,\La)$ of 
``complexes of finite tor-dimension 
with constructible cohomology'' (see \cite[Rapport 4.6]{SGA4.5} when $\La$ is finite and 
\cite[1.1.2-3]{De} in other cases). This category is known to be stable under the six 
operations  $f^*, f^!, f_*, f_!, \otimes$ and $\C{RHom}$ 
(see \cite[Th. finitude, 1.7]{SGA4.5}). In Section \ref{S:verdier} we will work with a 
larger category $D_{c}^b(X,\La)$ of ``complexes with constructible cohomology''. 

For each $X$ as above, we denote by $\pi_X:X\to\pt:=\Spec k$ the structure morphism, by 
$\La_X\in D_{ctf}^b(X,\La)$ the constant sheaf with fiber $\La$,  
by $K_X=\pi_X^!(\La_{\pt})$ the dualizing complex of $X$, and 
by $\B{D}=\B{D}_X:=\C{RHom}(\cdot,K_X)$ the Verdier duality. We will also write 
$R\Gm_c(X,\cdot)$ instead of $\pi_{X!}$ and $R\Gm(X,\cdot)$ instead of $\pi_{X*}$. 
For an embedding $i:Y\hra X$ and 
$\C{F}\in D_{ctf}^b(X,\La)$, we will often write $\C{F}|_Y$ instead of 
$i^*\C{F}$. 

Let $\B{F}$ be an algebraic closure of the finite field $\fq$. 
We say that an object $\C{X}$ over $\B{F}$ is {\em defined over} $\fq$, if it is 
a pullback of the corresponding object over $\fq$. 
\end{Emp'}

\begin{Emp'} \label{E:stmor}
(a) We will repeatedly use the fact that functor $f_*$ is right adjoint to $f^*$,
that  $f^!$ is right adjoint to $f_!$, and that $\C{RHom}$ is adjoint to $\otimes$.
We will denote by $adj$ the adjoint morphisms 
$\Id\to f_*f^*$, $f^*f_*\to\Id$, $\Id\to f^!f_!$,  $ f_!f^!\to\Id$, and by $ev$ the evaluation 
map $\C{A}\otimes \C{RHom}(\C{A},\C{B})\to \C{B}$. We also denote by $ev_{\C{F}}$ the evaluation
morphism $\B{D}\C{F}\otimes\C{F}\to K_X$. 

We will also freely use various base change morphisms 
(see, for example, \cite[XVII, 2.1.3 and XVIII, 3.1.12.3, 3.1.13.2, 3.1.14.2]{SGA4}), which we will denote by $BC$.

(b) For a morphism $f:X\to Y$, we denote by 
$t_f$ the morphism 
\begin{equation} \label{Eq:tf}
f_*\C{A}\otimes f_*\C{B}\to f_*(\C{A}\otimes\C{B}),
\end{equation}
adjoint to the composition 
$
f^*(f_*\C{A}\otimes f_*\C{B})=f^*f_*\C{A}\otimes f^*f_*\C{B}\overset{adj}{\lra}
\C{A}\otimes\C{B}$, by $proj$ the isomorphism 
$f_!(\C{A}\otimes f^*\C{B})\isom f_!\C{A}\otimes\C{B}$ (\cite[XVII, 5.2.9]{SGA4}), 
usually called the projection formula,
and by $t_{f^!}$ both the morphism 
\begin{equation} \label{Eq:tf!}
f^!\C{A}\otimes f^*\C{B}\to  f^!(\C{A}\otimes \C{B}),
\end{equation}
adjoint to the composition $f_!(f^!\C{A}\otimes f^*\C{B})\overset{proj}{\lra}f_!f^!\C{A}\otimes\C{B}
\overset{adj}{\lra}\C{A}\otimes\C{B}$
and the  isomorphism 
\begin{equation} \label{Eq:tf!2}
f^!\C{RHom}(\C{A},\C{B})\isom\C{RHom}(f^*\C{A},f^!\C{B})
\end{equation}
(see \cite[XVIII, 3.1.12.2]{SGA4}), adjoint to the composition
\[
f^!\C{RHom}(\C{A},\C{B})\otimes f^*\C{A} \overset{t_{f^!}}{\lra}
f^!(\C{RHom}(\C{A},\C{B})\otimes\C{A})\overset{ev}{\lra}f^!\C{B}.
\]
Also when $f$ is proper, we denote by $f_!$ the integration map
\[
H^0(X,K_X)=H^0(X,f^!K_Y)=H^0(Y,f_!f^!K_Y)\overset{adj}{\lra}H^0(Y,K_Y).
\]
\end{Emp'}

\begin{Emp'} \label{E:prod}
For a product $X_1\times X_2$, we denote by $p_i$ the projection 
$X_1\times X_2\to X_i$. For each $\C{A}\in D_{ctf}^b(X_1,\La)$ and 
$\C{B}\in D_{ctf}^b(X_2,\La)$, we set 
$\C{A}\pp\C{B}:=p_1^*\C{A}\otimes p_2^*\C{B}\in D_{ctf}^b(X_1\times X_2,\La)$.
In this setting Illusie \cite[(1.7.3) and (3.1.1)]{Il} constructed isomorphisms  
\begin{equation} \label{Eq:*}
\diam:K_{X_1}\pp\C{B}\isom p_2^!\C{B}\;\;\; \text{ and }\;\;\;
\B{D}\C{A}\pp\C{B}\isom \C{RHom}(p_1^*\C{A},p_2^!\C{B}),
\end{equation}
the first of which is the composition 
\[
p_1^*\pi_{X_2}^!\La_{\pt}\otimes p_2^*\C{B}\overset{BC}{\lra}
p_2^!\pi_{X_1}^*\La_{\pt}\otimes p_2^*\C{B}\overset{t_{p_2^!}}{\lra}p_2^!(\La_{X_1}\otimes\C{B}), 
\]
while the second is adjoint to the composition
\[
p_1^*\C{A}\otimes (\B{D}\C{A}\pp\C{B})=(\C{A}\otimes\B{D}\C{A})\pp\C{B}
\overset{ev}{\lra}K_{X_1}\pp\C{B}\overset{\diam}{\lra} p_2^!\C{B}.
\]
In particular, isomorphism $\diam$ induces an isomorphism $K_{X_1}\pp K_{X_2}\isom 
 K_{X_1\times X_2}$.
\end{Emp'}

\section{Basic constructions}
\subsection{Correspondences} \label{SS:corr}

\begin*
\vskip 8truept
\end*

\begin{Def} \label{D:corr} 
(a)  By a {\em correspondence}, we mean a morphism of schemes of the form  
$c=(c_1,c_2):C\to X_1\times X_2$.

(b) Let $c:C\to X_1\times X_2$ and $b:B\to Y_1\times Y_2$ be correspondences.
By a {\em morphism} from $c$ to $b$ we mean a triple $[f]=(f_1,f^{\diez},f_2)$ making the following 
diagram commutative 
\begin{equation} \label{Eq:funct}
\CD 
        X_1   @<{c_1}<<                    C        @>{c_2}>>         X_2\\
        @V{f_1}VV                        @V{f^{\diez}}VV                       @VV{f_2}V\\
        Y_1 @<{b_1}<<                    B    @>{b_2}>>            Y_2.
\endCD 
\end{equation}
\end{Def}

\begin{Not} \label{N:excor}
(a)  Denote by $c_{tr}$ the trivial correspondence $\pt\to\pt\times\pt$.
For an arbitrary correspondence $c:C\to X_1\times X_2$, denote by $[\pi]_c$ the structure morphism
$(\pi_{X_1},\pi_{C},\pi_{X_2})$ from $c$ to $c_{tr}$.

(b) We say that a morphism $[f]=(f_1,f^{\diez},f_2)$ between correspondences is
 {\em proper} (resp. {\em an open embedding}, resp. {\em a closed embedding}), 
if each component $f_1,f^{\diez}$ and $f_2$ satisfies this property. We say that
a correspondence $c$ is {\em proper} over $k$, if the structure morphism $[\pi]_c$ is proper.

(c) Let  $c:C\to X_1\times X_2$ be a correspondence. By a {\em compactification of $c$} we mean a 
proper correspondence $\ov{c}:\ov{C}\to\ov{X}_1\times\ov{X}_2$ over $k$ equipped with an open 
embedding $[j]=(j_1,j^{\diez},j_2)$ from $c$ to $\ov{c}$ such that 
$j_1,j^{\diez}$ and $j_2$ are dominant. 
\end{Not}
The following observation will be useful later.

\begin{Rem} \label{R:comp}
Let $c:C\to X_1\times X_2$ be a correspondence. Then every pair of compactifications 
$j_1:X_1\hra\ov{X}_1$ and $j_2:X_2\hra\ov{X}_2$ can be extended to a compactification 
$\ov{c}:\ov{C}\to \ov{X}_1\times \ov{X}_2$ of $c$. Indeed, choose first any compactification 
$j':C\hra C'$ of $C$, define $\ov{C}$ to be the closure of the image of 
\[
(j', j_1\circ c_1,j_2\circ c_2):C\hra C'\times\ov{X}_1\times\ov{X}_2,
\] 
and denote the
projection $\ov{C}\to\ov{X}_1\times\ov{X}_2$ by $\ov{c}$. 
\end{Rem}

\begin{Def} \label{D:cohcor}
 Let $c:C\to X_1\times X_2$ be a correspondence. By a {\em $c$-morphism}, we mean a morphism 
$u:c_{2!}c_1^*\C{F}_1\to\C{F}_2$ for some $\C{F}_1\in D^b_{ctf}(X_1,\La)$ and 
$\C{F}_2\in D^b_{ctf}(X_2,\La)$. 
\end{Def}

\begin{Rem}
A $c$-morphism is usually called a {\em cohomological correspondence lifting $c$}.
\end{Rem}


\begin{Emp} \label{E:pushf}
{\bf Push-forward of cohomological correspondences.}

(a) In the notation of \ref{D:corr} (b) assume any of the following three conditions:

\indent\indent(i) the left inner square of \form{funct} is Cartesian;

\indent\indent(ii) morphisms $f_1$ and $f^{\diez}$ are proper;

\indent\indent(iii) morphisms $c_1$ and $b_1$ are proper.

In all these cases, we have a base change morphism  
$BC:b_1^* f_{1!}\to f^{\diez}_! c_1^*$. Namely, in the cases (i) and (ii), $BC$ is the standard 
one, while in the case (iii), $BC$ is adjoint to the map
\[
f_{1!}\overset{adj}{\lra} f_{1!} c_{1*} c_1^*= f_{1!} c_{1!} c_1^*= b_{1!}f^{\diez}_! c_1^*=
b_{1*}f^{\diez}_! c_1^*.
\]
Every $c$-morphism $u:c_{2!}c_1^*\C{F}_1\to\C{F}_2$ gives rise to a 
$b$-morphism 
\[
[f]_!(u):b_{2!}b_1^*(f_{1!}\C{F}_1)\overset{BC}{\lra} b_{2!}f^{\diez}_! c_1^*\C{F}_1=
 f_{2!}c_{2!}c_1^*\C{F}_1\overset{u}{\lra} f_{2!}\C{F}_2.
\]

(b) Proper push-forwards from (a) are compatible with compositions. Namely, 
let  $[g]$ be a morphism from a correspondence $b$ to a third correspondence satisfying
the same condition among (i)-(iii) as $[f]$.  Then the push-forward 
$([g]\circ[f])_!$ defined in (a) equals the composition $[g]_!\circ [f]_!$.
Indeed, proper push-forwards are induced by base change morphisms, so the assertion 
follows from the fact that base change morphisms are compatible with compositions.
\end{Emp}

\begin{Emp} \label{E:rgm}
{\bf Functor $R\Gm_c$.}
Let $c:C\to X_1\times X_2$ be a correspondence with $c_1$ is proper, and let 
$u: c_{2!}c_1^*\C{F}_1\to\C{F}_2$ be a $c$-morphism.

(a) The structure 
morphism $[\pi]_c$ from \rn{excor} (a) satisfies assumption
(iii) of \re{pushf} (a). Hence $u$ gives rise to a morphism
\[
R\Gm_c(u):=[\pi]_{c!}(u):R\Gm_c(X_1,\C{F}_1)\to R\Gm_c(X_2,\C{F}_2).
\]
Explicitly, $R\Gm_c(u)$ is the composition 
\[
R\Gm_c(X_1,\C{F}_1)\overset{c_1^*}{\lra}R\Gm_c(C, c_1^*\C{F}_1)=
R\Gm_c(X_2,c_{2!}c_1^*\C{F}_1)\overset{u}{\lra}R\Gm_c(X_2,\C{F}_2).
\]

(b) Let $\ov{c}$  be a compactification of $c$, and let 
$[j]=(j_1,j^{\diez},j_2)$ be the corresponding open embedding of $c$ into $\ov{c}$.
Then $[j]$  satisfies assumption
(iii) of \re{pushf} (a), hence $u$ extends to a $\ov{c}$-morphism 
$\ov{u}:=[j]_!(u):\ov{c}_{2!}\ov{c}_1^*(j_{1!}\C{F}_1)\to j_{2!}\C{F}_2$. 

(c) In the notation of (b), we have $[\pi]_{c!}=[\pi]_{\ov{c}!}\circ [j]_!$ (use \re{pushf} (b)). 
Therefore natural isomorphisms $R\Gm_c(X_1,\C{F}_1)\cong R\Gm_c(\ov{X_1},j_{1!}\C{F}_1)$ and 
$R\Gm_c(X_2,\C{F}_2)\cong R\Gm_c(\ov{X_2},j_{2!}\C{F}_2)$
identify $R\Gm_c(u)$ with $R\Gm_c(\ov{u})$.
\end{Emp}

\begin{Emp} \label{E:pullb}
{\bf Pullback of cohomological correspondences.}

\noindent In the notation of \ref{D:corr} (b), assume any of the following conditions:

(i) the right inner square of \form{funct} is Cartesian ``up to nilpotents'',
i.e., the canonical 
\indent\;\;\;\: morphism $C_{red}\to (B\times_{Y_2} X_2)_{red}$ is an isomorphism;

(ii) morphisms $f^{\diez}$ and $f_2$ are \'etale.

\noindent In both cases, we have a base change morphism $BC:c_{2!}f^{\diez*}\to f_2^*b_{2!}$.
 Hence every $b$-morphism $u:b_{2!}b_1^*\C{F}_1\to\C{F}_2$ 
gives rise to a $c$-morphism

\[
[f]^*(u):c_{2!}c_1^*(f_{1}^*\C{F}_1)=c_{2!}f^{\diez*}b_1^*\C{F}_1 \overset{BC}{\lra}
 f_2^*b_{2!}b_1^*\C{F}_1\overset{f_2^*u}{\lra} f_{2}^*\C{F}_2.
\]
\end{Emp}

\begin{Emp} \label{E:restr}
{\bf Restriction of correspondences.}
Let $c:C\to X\times X$ be a correspondence, and  
$u:c_{2!}c_1^*\C{F}\to\C{F}$ a $c$-morphism.

 For open subsets $W\subset C$ and $U\subset X$, 
denote by $c|_W:W\to X\times X$ and $c|_U:c_1^{-1}(U)\cap c_2^{-1}(U)\to U\times U$
the restrictions of $c$, and let $[j_W]$ and $[j_U]$ be the open embeddings of $c|_W$ 
 and $c|_U$ into $c$, respectively. Then $[j_W]$ and $[j_U]$ satisfy assumption (ii) of  
\re{pullb}, so $u$ restricts to a $c|_W$-morphism $u|_W:=[j_W]^*(u)$
and a  $c|_U$-morphism $u|_U:=[j_U]^*(u)$.
\end{Emp}

\subsection{Trace maps} \label{SS:trace}

\begin*
\vskip 8truept
\end*

\begin{Not} \label{N:partcase}
For a correspondence $c:C\to X\times X$, let $\Dt$ be the diagonal map 
$X\hra X\times X$, put $Fix(c)$ be the fiber product $X\times_{X\times X} C$, and denote 
by $\Dt': Fix(c)\hra C$ and $c':Fix(c)\to X$
the inclusion map and the restriction of $c$, respectively.
We call $Fix(c)$ {\em the scheme of fixed points of $c$}.
\end{Not}

The following construction is crucial for this work.

\begin{Emp} \label{E:locterms}
{\bf Trace maps and local terms.} Fix  a correspondence $c:C\to X\times X$ 
and $\C{F}\in D^b_{ctf}(X,\La)$.

(a) Denote by $\un{\C{Tr}}:\C{RHom}(c_1^*\C{F},c_2^!\C{F})\to\Dt'_*K_{Fix(c)}$
the composition
\[
\C{RHom}(c_1^*\C{F},c_2^!\C{F})\isom c^!(\B{D}\C{F}\pp\C{F})\to c^!\Dt_*K_X
\isom\Dt'_* c'^!K_X=\Dt'_*K_{Fix(c)},
\]
where the first map is the composition of the inverses of \form{tf!2} and \form{*}
\[
\C{RHom}(c_1^*\C{F},c_2^!\C{F})\isom c^!\C{RHom}(p_1^*\C{F},p_2^!\C{F})\isom
c^!(\B{D}\C{F}\pp\C{F}),
\]
the second one is induced by the map $\B{D}\C{F}\pp\C{F}\to\Dt_*K_{X}$, 
adjoint to the evaluation map $\Dt^*(\B{D}\C{F}\pp\C{F})=\B{D}\C{F}\otimes\C{F}\to K_{X}$,
and the third one is the base change isomorphism  
$c^!\Dt_*\isom \Dt'_* c'^!$.

 Using identifications
\[
H^0(C,\C{RHom}(c_1^*\C{F},c_2^!\C{F}))=\Hom(c_1^*\C{F},c_2^!\C{F})\isom \Hom(c_{2!}c_1^*\C{F},\C{F}),
\]
where the last isomorphism is obtained by adjointness, the map  $H^0(C,\un{\C{Tr}})$ gives rise to
the map
\begin{equation} \label{Eq:lterm}
\C{Tr}:\Hom(c_{2!}c_1^*\C{F},\C{F})\to H^0(C,\Dt'_*K_{Fix(c)})=H^0(Fix(c),K_{Fix(c)}).
\end{equation}

Whenever necessary, we will also use notation $\un{\C{Tr}}_c$, $\un{\C{Tr}}_{\C{F}}$,
$\C{Tr}_c$ or $\C{Tr}_{\C{F}}$ to emphasize that the trace maps $\un{\C{Tr}}$ and 
$\C{Tr}$  depend on $c$ and $\C{F}$.

(b) For an open subset $\beta$ of $Fix(c)$, we denote by 
\begin{equation} \label{Eq:lterm2}
\C{Tr}_{\beta}: Hom(c_{2!}c_1^*\C{F},\C{F})\to H^0(\beta, K_{\beta})
\end{equation}
the composition of $\C{Tr}$ and the restriction map 
$H^0(Fix(c),K_{Fix(c)})\to  H^0(\beta, K_{\beta})$.
If moreover, $\beta$ is proper over $k$, we denote by 
\begin{equation} \label{Eq:lterm3}
LT_{\beta}:Hom(c_{2!}c_1^*\C{F},\C{F})\to\La
\end{equation}
the composition of $\C{Tr}_{\beta}$ and the integration map 
$\pi_{\beta!}:H^0(\beta, K_{\beta})\to \La$.

When $\beta$ is a connected component of $Fix(c)$ which is proper over $k$, 
then $LT_{\beta}(u)$ is usually 
called the {\em local term} of $u$ at $\beta$.
\end{Emp}

\begin{Emp} \label{E:extr}
{\bf Example.} If $c=c_{tr}$ (see \rn{excor} (a)), then $Fix(c)=\pt$, hence
$H^0(Fix(c),K_{Fix(c)})=\La$. Moreover, in this case $\C{F}$ is just a bounded complex of 
finitely generated free $\La$-modules (modulo homotopy), 
and the trace map $\C{Tr}_{c_{tr}}$ coincides with the usual trace map 
$\Hom(\C{F},\C{F})\to\La$.
\end{Emp} 


\begin{Rem}
We will show in the appendix that our local terms are equivalent to those 
defined by Illusie \cite{Il} and later used by Pink \cite{Pi} and 
Fujiwara \cite{Fu}. However, our definition is more elementary.
\end{Rem}

The following proposition, whose proof will be given in \rs{funct} (see \re{proofs}), 
asserts that the trace maps commute with proper push-forwards. Though the result resembles 
\cite[Cor 4.5]{Il}, we do not assume that morphisms $c:C\to X\times X$ and 
$b:B\to Y\times Y$ are proper. 

\begin{Prop} \label{P:pushf} 
Let $[f]=(f,f^{\diez},f)$ be a proper morphism from a correspondence 
$c:C\to X\times X$ to $b:B\to Y\times Y$. Then the morphism $f':Fix(c)\to Fix(b)$ induced by $[f]$ 
is proper as well, and for every $\C{F}\in D^b_{ctf}(X,\La)$, the following diagram commutes
\begin{equation} \label{Eq:pushf}
\CD 
       Hom(c_{2!}c_1^*\C{F},\C{F})   @>\C{Tr}_c>> H^0(Fix(c), K_{Fix(c)})\\
      @V{[f]_!}VV            @V{f'_!}VV\\
       Hom(b_{2!}b_1^* (f_!\C{F}),f_!\C{F}) @>\C{Tr}_b>>  H^0(Fix(b), K_{Fix(b)}).
\endCD 
\end{equation}
\end{Prop}

As a particular case, we deduce the well-known Lefschetz-Verdier trace formula 
(\cite[Cor. 4.7]{Il}).

\begin{Cor} \label{C:LTF}
Let  $c:C\to X\times X$ be a proper correspondence over $k$. 
Then for every $c$-morphism 
 $u: c_{2!}c_1^*\C{F}\to\C{F}$, we have an equality
\begin{equation} \label{Eq:ltf} 
\Tr(R\Gm_c(u))=\sum_{\beta\in \pi_0(Fix(c))}LT_{\beta}(u).
\end{equation}
\end{Cor}

\begin{proof}
Consider the structure morphism $[\pi]=[\pi]_c$ from \rn{excor} (a). 
Then the right-hand side of \form{ltf} equals
$\pi_{Fix(c)!}(\C{Tr}_c(u))$, while the left-hand side equals
$\Tr(R\Gm_c(u))=\Tr([\pi]_!(u))$ (see \re{rgm} (a) and \re{extr}). Hence the equality 
\form{ltf} 
follows from the commutativity of \form{pushf}.
\end{proof}

\subsection{Specialization} \label{SS:spec}

\begin*
\vskip 8truept
\end*

\begin{Emp} \label{E:hens}
{\bf Set up.}
Let $k$ be a separably closed field, $R$ a discrete valuation domain over $k$, whose residue
field is $k$. Let $K$ be the fraction field of $R$,
$\ov{R}$ the integral closure of $R$ in $K^{sep}$, and $R^h$  the (strict) henselization of 
$R$. 

Set $\C{D}:=\Spec R$, $\ov{\C{D}}:=\Spec \ov{R}$ and $\C{D}^h:=\Spec R^h$, denote by 
$\eta$, $\ov{\eta}$, and $\eta^h$ the generic points of $\C{D}$,  $\ov{\C{D}}$ and 
$\C{D}^h$, respectively, and by $s$ the special points of each $\C{D}$, $\ov{\C{D}}$ and 
$\C{D}^h$.

(a) For a scheme $\wt{X}$ over $\C{D}$ and $?=s,\eta,\ov{\eta},\eta^h,\ov{\C{D}},\C{D}^h$, 
we denote by $\wt{X}_?$, the corresponding scheme over $?$.
Let $\ov{i}:\wt{X}_s\hra\wt{X}_{\ov{\C{D}}}$ and 
$\ov{j}:\wt{X}_{\ov{\eta}}\hra\wt{X}_{\ov{\C{D}}}$ be 
the canonical closed an open embeddings.

(b) For a morphism $\wt{f}:\wt{X}\to\wt{Y}$ of schemes over $\C{D}$ and  
$?=s,\eta,\ov{\eta},\ov{\C{D}}$, we 
denote by $\wt{f}_{?}$ the corresponding morphism $\wt{X}_?\to\wt{Y}_?$.

(c) For a scheme $X$ over $k$, set $X_{\C{D}}:=X\times \C{D}$. For a morphism $f:X\to Y$ 
of schemes over $k$, set $f_{\C{D}}:=f\times\Id_{\C{D}}:X_{\C{D}}\to Y_{\C{D}}$.
For each $?=s,\eta,\ov{\eta},\ov{\C{D}}$, we will write $X_{?}$ instead of $(X_{\C{D}})_{?}$
in the notation of (a) and $f_{?}$ instead of $(X_{\C{D}})_{?}$ in the notation of (b).
We will also identify $X_s=(X_{\C{D}})_s$ with $X$.

(d) For a scheme $X$ over $k$ and an object $\C{F}\in D_{ctf}^b(X,\La)$, we denote
by $\C{F}_{\eta}$ and $\C{F}_{\ov{\C{D}}}$ its pullbacks to 
$D_{ctf}^b(X_{\eta},\La)$ and $D_{ctf}^b(X_{\ov{\C{D}}},\La)$. 
For a scheme $\wt{X}_{\eta}$ over $\eta$ and an object 
$\C{F}_{\eta}\in D_{ctf}^b(\wt{X}_{\eta},\La)$, denote by 
$\C{F}_{\ov{\eta}}$ and $\C{F}_{\eta^h}$ its pullbacks to $D_{ctf}^b(\wt{X}_{\ov{\eta}},\La)$
and  $D_{ctf}^b(\wt{X}_{\eta^h},\La)$.

(e) Denote by $\Psi_{\wt{X}}: D_{ctf}^b(\wt{X}_{\eta},\La)\to D_{ctf}^b(\wt{X}_{s},\La)$ 
the corresponding functor of nearby cycles. Explicitly, 
$\Psi_{\wt{X}}(\C{F}_{\eta}):=\ov{i}^*\ov{j}_*\C{F}_{\ov{\eta}}$, if $R$ is strictly henselian 
(compare, for example, \cite[$\S$4]{Il3}), 
and $\Psi_{\wt{X}}(\C{F}_{\eta}):=\Psi_{\wt{X}_{\C{D}^h}}(\C{F}_{\eta^h})$, in general.
\end{Emp}

\begin{Emp} \label{E:specfunc}
{\bf Specialization functor.}

(a) We say that a separated scheme $\wt{X}$ of finite type over $\C{D}$ {\em lifts} 
a scheme $X$ over $k$, if it is equipped with a morphism 
$\varphi=\varphi_X:\wt{X}\to X_{\C{D}}$ such that $\varphi_{\eta}:\wt{X}_{\eta}\to X_{\eta}$
is an isomorphism. In this case, we define a functor 
\[
sp_{\wt{X}}:=sp_{\wt{X},\varphi}:D_{ctf}^b(X,\La)\to D_{ctf}^b(\wt{X}_{s},\La) 
\]
by the rule $sp_{\wt{X}}(\C{F}):=\Psi_{\wt{X}}(\varphi_{\eta}^*\C{F}_{\eta})$.

(b) We say that a morphism $\wt{f}:\wt{X}\to\wt{Y}$ of schemes over $\C{D}$
{\em lifts} a  morphism $f:X\to Y$ of schemes over $k$, if $\wt{X}$ lifts $X$, 
$\wt{Y}$ lifts $Y$, and $\varphi_Y\circ\wt{f}=f_{\C{D}}\circ\varphi_X$. 
In this case, we have base change morphisms
$BC^*:\wt{f}_s^* sp_{\wt{Y}}\to sp_{\wt{X}} f^*$, 
$BC_*: sp_{\wt{Y}} f_*\to \wt{f}_{s*} sp_{\wt{X}}$, 
$BC^!: sp_{\wt{X}} f^!\to \wt{f}_s^! sp_{\wt{Y}}$ and
$BC_!:\wt{f}_{s!} sp_{\wt{X}}\to sp_{\wt{Y}} f_!$, defined as follows.

Since the map $\C{D}^h\to\C{D}$ is etale, the functor $\C{F}_{\eta}\mapsto\C{F}_{\eta^h}$ commutes 
with all the operations. Thus it remains to define the base change functors in the case when  
$R$ is strictly henselian. In this case, $BC^*$ is by definition the composition 
\[
\wt{f}_s^*sp_{\wt{Y}}\C{F}=\wt{f}_s^*\ov{i}^*\ov{j}_*\varphi^*_Y\C{F}_{\ov{\eta}}=
\ov{i}^*\wt{f}_{\ov{\C{D}}}^*\ov{j}_*\varphi^*_Y\C{F}_{\ov{\eta}}\overset{BC}{\lra} 
\ov{i}^*\ov{j}_*\wt{f}_{\ov{\eta}}^*\varphi^*_Y\C{F}_{\ov{\eta}}=\ov{i}^*\ov{j}_*
\varphi^*_X(f^*\C{F})_{\ov{\eta}}=sp_{\wt{X}}f^*\C{F},
\]
 $BC_!$ is the composition 
\[
\wt{f}_{s!}sp_{\wt{X}}\C{F}=\wt{f}_{s!}\ov{i}^*\ov{j}_*\varphi^*_X\C{F}_{\ov{\eta}}\cong
\ov{i}^*\wt{f}_{\ov{\C{D}}!}\ov{j}_*\varphi^*_X\C{F}_{\ov{\eta}}\overset{BC}{\lra} 
\ov{i}^*\ov{j}_*\wt{f}_{\ov{\eta}!}\varphi^*_X\C{F}_{\ov{\eta}}\cong\ov{i}^*\ov{j}_*
\varphi^*_Y(f_!\C{F})_{\ov{\eta}}=sp_{\wt{X}}f_!\C{F},
\]
while  $BC_*$ and $BC^!$ are defined by adjointness.

(c) If $\wt{f}$ is proper, it follows from the proper base change theorem that the 
base change morphism $BC_!$ is an isomorphism, and $BC_*$ is its inverse. Similarly, the 
base change morphism $BC^*$ is an isomorphism, if  $\wt{f}$ is smooth. 
\end{Emp}

\begin{Emp} \label{E:trivdef}
{\bf Examples.}
(a) If $\wt{X}=X_{\C{D}}$ and $\varphi$ is the identity, then $\wt{X}_s=X$, and the natural
morphism of functors $\C{F}=\ov{i}^*\C{F}_{\ov{\C{D}}}\overset{adj}{\lra}\ov{i}^*
\ov{j}_*\ov{j}^*\C{F}_{\ov{\C{D}}}= sp_{\wt{X}}\C{F}$ is an isomorphism (see \rr{trivdef} below).

(b) Let $f:X\to\pt$ and $\wt{f}:\wt{X}\to\C{D}$ be the structure morphisms. 
Then the composition 
\[
BC^!\circ BC_*:f_* f^!\La_{\pt}=sp_{\C{D}}f_* f^!\La_{\pt}\to 
\wt{f}_{s*}\wt{f}_s^! sp_{\C{D}}\La_{\pt}=\wt{f}_{s*}\wt{f}_s^!\La_{\pt}
\]
defines a homomorphism $H^0(X,K_X)\to H^0(\wt{X}_s,K_{\wt{X}_s})$, 
which we will denote simply by $sp_{\wt{X}}$.
\end{Emp}

\begin{Emp} \label{E:cor}
{\bf Specialization of correspondences.}
Let  $c:C\to X\times X$ be a correspondence over $k$, $\wt{c}:\wt{C}\to \wt{X}\times \wt{X}$ a 
correspondence over $\C{D}$ lifting $c$, and 
$\wt{c}_{s}$ the special fiber of $\wt{c}$. Then every $c$-morphism 
$u:c_{2!}c_1^*\C{F}_1\to\C{F}_2$ gives rise to a  $\wt{c}_{s}$-morphism
\[
sp_{\wt{c}}(u):\wt{c}_{s2!}\wt{c}_{s1}^* sp_{\wt{X}}\C{F}_1
\overset{BC^*}{\lra}\wt{c}_{s2!} sp_{\wt{C}} c_{1}^*\C{F}_1
\overset{BC_!}{\lra} sp_{\wt{X}}c_{2!}c_1^*\C{F}_1
\overset{u}{\lra}
 sp_{\wt{X}}\C{F}_2.
\] 
\end{Emp}

The following proposition,  whose proof will be given in \rs{funct} (see \re{proofs}), 
asserts that the trace maps commute with specialization. Though the result resembles 
\cite[Prop. 1.7.1]{Fu}, we do not assume that $c$ and $\wt{c}$ are proper.

\begin{Prop} \label{P:spec}
Let $c:C\to X\times X$ be a  correspondence over $k$ and let 
$\wt{c}:\wt{C}\to \wt{X}\times \wt{X}$ be a 
correspondence over $\C{D}$ lifting $c$. 
Then for each $\C{F}\in D_{ctf}^b(X,\La)$, the following diagram 
is commutative
\begin{equation} \label{Eq:spec}
\CD 
      \Hom(c_{2!}c_{1}^*\C{F},\C{F})
  @>\C{Tr}_c>> H^0(Fix(c), K_{Fix(c)})\\
      @V{sp_{\wt{c}}}VV            @V{sp_{Fix(\wt{c})}}VV\\
      \Hom(\wt{c}_{s2!}\wt{c}_{s1}^* sp_{\wt{X}}\C{F},
 sp_{\wt{X}}\C{F})
      @>\C{Tr}_{\wt{c}_s}>>  H^0(Fix(\wt{c}_s), K_{Fix(\wt{c}_s)}).
\endCD 
\end{equation}
\end{Prop}

\subsection{ Deformation to the normal cone} \label{SS:normcone}

\begin*
\vskip 8truept
\end*


We will apply the construction of the previous subsection in the following particular case.
We fix any $R$ as in \ref{SS:spec}, for example, $R=k[t]_{(t)}$ (the localization of $k[t]$
at $(0)$)  or $R=k[[t]]$.

\begin{Not} \label{N:deform}
Let $X$ be a scheme over $k$, and $Z\subset X$ a closed subscheme.
We will identify $X$ with the special fiber $(X_{\C{D}})_s$ of $X_{\C{D}}$.

(a) Denote by $\wt{X}_{Z}$ the affine scheme 
$\C{Spec}(\C{O}_{X_{\C{D}}}[\frac{\C{I}_{Z_{\C{D}}}}{t}])$ over $X_{\C{D}}$, where
$t$ is any uniformizer of $R$, and  
$\C{O}_{X_{\C{D}}}[\frac{\C{I}_{Z_{\C{D}}}}{t}]$ is an $\C{O}_{X_{\C{D}}}$-subalgebra of 
$\C{O}_{X_{\C{D}}}[\frac{1}{t}]=\C{O}_{X_{\eta}}$. 

(b) The embedding $\C{O}_{X_{\C{D}}}\hra \C{O}_{X_{\C{D}}}[\frac{\C{I}_Z}{t}]$ gives rise to the 
birational projection $\wt{\varphi}:\wt{X}_{Z}\to X_{\C{D}}$, which is an isomorphism over 
$\eta$. Thus $\wt{X}_{Z}$ is a lift of $X$ in the sense of \re{specfunc} (a).

(c) The special fiber $(\wt{X}_{Z})_s=\C{Spec}(\bigoplus_{n=0}^{\infty}
(\C{I}_Z)^n/(\C{I}_Z)^{n+1})$ of $\wt{X}_{Z}$ is the normal cone of $X$ to $Z$, 
which we will denote by $N_{Z}(X)$. 

(d) The projection $\C{O}_{X_{\C{D}}}[\frac{\C{I}_{Z_{\C{D}}}}{t}]\to
\C{O}_{X_{\C{D}}}/\C{I}_{Z_{\C{D}}}=\C{O}_{Z_{\C{D}}}$ defines a closed embedding 
$\wt{i}:Z_{\C{D}}\hra\wt{X}_{Z}$. The special fiber
$\wt{i}_s:Z\hra N_{Z}(X)$ corresponds to the projection to the zero's component 
$\bigoplus_{n=0}^{\infty}(\C{I}_Z)^n/(\C{I}_Z)^{n+1}\to \C{O}_X/\C{I}_Z=\C{O}_Z$,
thus identifies $Z$ with the zero section of  $N_{Z}(X)$. 

(e) The composition of the map $BC^*$ from \re{specfunc} (b), corresponding to the embedding 
$\wt{i}$ from (d), and  the isomorphism of \re{trivdef} (a) defines a morphism

\begin{equation} \label{Eq:verdier}
sp_{\wt{X}_{Z}}(\C{F})|_Z\to sp_{Z_{\C{D}}}(\C{F}|_Z)\isom\C{F}|_Z.
\end{equation} 
\end{Not}

The following result was proven by Verdier in \cite[$\S$8, (SP5)]{Ve}. 
Since this result is crucial for our argument, 
while Verdier's argument is slightly sketchy and contains a small gap, 
we will present its proof in 
Section \ref{S:verdier}. 

\begin{Prop} \label{P:verdier}
For each $\C{F}\in D_{ctf}^b(X,\La)$, the map \form{verdier}
is an isomorphism.
\end{Prop}


The following functorial properties of $\wt{X}_Z$ will be used later.

\begin{Lem} \label{L:normfunc}
Let  $f:X'\to X$ be a morphism of schemes over $k$, $Z\subset X$ 
a closed subscheme, and $Z'$ a closed subscheme of $f^{-1}(Z)$, that is, 
$f^{\cdot}(\C{I}_{Z})\subset\C{I}_{Z'}$. 

(a) The morphism $f$ has a unique lift $\wt{f}:\wt{X'}_{Z'}\to\wt{X}_{Z}$ over 
$\C{D}$.

(b) The image $\wt{f}_s({N}_{Z'}(X'))$ is supported set-theoretically at 
the zero section $Z\subset N_{Z}(X)$ if and only if there exists $n\in \B{N}$ such that 
 $f^{\cdot}(\C{I}_{Z})^n\subset (\C{I}_{Z'})^{n+1}$.

(c) If $Z'=f^{-1}(Z)$, then the map $\wt{X'}_{Z'}\to\wt{X}_{Z}\times_{X}X'$ induced by $\wt{f}$ 
is a closed embedding. Moreover, the latter map is an isomorphism, if $f$ is flat. 

(d) If $Z'=f^{-1}(Z)$, then the preimage $\wt{f}_s^{-1}(Z)$ 
of the zero section $Z\subset N_{Z}(X)$ equals $Z'$.

(e)  If $Z'=f^{-1}(Z)$ and $f$ is a closed embedding (resp. proper, resp. smooth), then 
$\wt{f}$ is a closed embedding  (resp. proper, resp. smooth) as well. 
\end{Lem}
\begin{proof}
(a) Note that $f_{\eta}:X'_{\eta}\to X_{\eta}$ corresponds to the morphism of sheaves
 $\wt{f}^{\cdot}:\C{O}_{X_{\C{D}}}[\frac{1}{t}]\to\C{O}_{X'_{\C{D}}}[\frac{1}{t}]$, induced by 
$f^{\cdot}:\C{O}_{X}\to \C{O}_{X'}$. 
Since $f^{\cdot}(\C{I}_{Z})\subset \C{I}_{Z'}$, we obtain the inclusion 
$\wt{f}^{\cdot}(\C{O}_{X_{\C{D}}}[\frac{\C{I}_{Z_{\C{D}}}}{t}])
\subset\C{O}_{X'_{\C{D}}}[\frac{\C{I}_{Z'_{\C{D}}}}{t}]$, which implies 
that $f_{\eta}$ extends uniquely to the morphism  $\wt{f}:\wt{X'}_{Z'}\to\wt{X}_{Z}$, 
as claimed. 

(b) Since $\wt{f}_s:{N}_{Z'}(X')\to {N}_{Z}(X)$ corresponds to 
the map 
$\bigoplus_{n=0}^{\infty}(\C{I}_{Z})^n/(\C{I}_{Z})^{n+1}\to 
\bigoplus_{n=0}^{\infty}(\C{I}_{Z'})^n/(\C{I}_{Z'})^{n+1}$, induced by 
$f^{\cdot}$, the assertion follows. 

(c) follows from the fact that 
$\C{I}_{f^{-1}(Z)}\subset \C{O}_{X'}$ equals by definition the image
of $\C{I}_{Z}\otimes_{\C{O}_{X'}}\C{O}_{X'}$.

(d) By (c), $\wt{f}_s^{-1}(Z)$ is a closed subscheme of $Z\times_{X} X'=Z'$.
Since the opposite inclusion always holds, we obtain the assertion.

(e) follows immediately from (c).
\end{proof}

\begin{Emp} \label{E:defcor}
{\bf Deformation of correspondences.}
Let $c:C\to X\times X$ be a correspondence, and $Z\subset X$ a closed subscheme.
By \rl{normfunc} (a), $c$ lifts to a correspondence
$\wt{c}_{Z}:\wt{C}_{c^{-1}(Z\times Z)}\to\wt{X}_{Z}\times \wt{X}_{Z}$ over $\C{D}$.
\end{Emp}

\begin{Cor} \label{C:diag}
In the notation of \re{defcor}, $\wt{Fix(c)}_{c'^{-1}(Z)}$ is naturally
a closed subscheme of $Fix(\wt{c}_{Z})$.
\end{Cor}
\begin{proof}
Since $Fix(c)$ is a closed subscheme of $C$, and $c'^{-1}(Z)\subset Fix(c)$ is the 
schematic preimage of $c^{-1}(Z\times Z)\subset C$, we conclude from \rl{normfunc} (e)
that $\wt{Fix(c)}_{c'^{-1}(Z)}$ is naturally
a closed subscheme of $\wt{N}_{c^{-1}(Z\times Z)}(C)$. Since by 
\rl{normfunc} (a) the image $\wt{c}_{Z}(\wt{Fix(c)}_{c'^{-1}(Z)})$ is 
(schematically) contained in 
$\Dt(\wt{X}_{Z})\subset\wt{X}_{Z}\times \wt{X}_{Z}$, we obtain the assertion.
\end{proof}

\subsection{Locally invariant subsets} \label{SS:locinv}

\begin*
\vskip 8truept
\end*

\begin{Def}
Let $c:C\to X\times X$ be a correspondence, and $Z\subset X$ a closed subset.

(a) We say that $Z$ is {\em $c$-invariant}, if  
$c_1(c_2^{-1}(Z))\subset X$ is set-theoretically contained in  $Z$.

(b) We say that $Z$ is {\em $c$-invariant in a neighborhood of fixed points},
if there exists an open neighborhood $W\subset C$ of $Fix(c)$ such that
$Z$ is $c|_W$-invariant.

(c) We say that $Z$ is {\em locally $c$-invariant},
if for each $x\in Z$ there exists an open neighborhood $U\subset X$ of $x$ such that
$c_1 (c_2^{-1}(Z\cap U))$ is contained in $Z\cup (X\sm U)$, or, equivalently, such that
$Z\cap U\subset U$ is $c|_U$-invariant.
\end{Def}

\begin{Emp} \label{E:qfinite} 
{\bf Example.} If $x\in X$ is a closed point such that $c_2^{-1}(x)$ is finite, then 
 $\{x\}$ is locally $c$-invariant. Indeed, $U:=X\sm [c_1(c_2^{-1}(x))\sm x]$ is the required 
open set.
\end{Emp}

\begin{Lem} \label{L:invneib}
Let $c:C\to X\times X$ be a correspondence, and $Z\subset X$ a closed subset.

(a) If  $Z$ is $c$-invariant, then $Z$ is locally $c$-invariant.

(b) The open subset $W(Z):=C\sm [\ov{c_2^{-1}(Z)\sm c_1^{-1}(Z)}]$, where 
$\ov{c_2^{-1}(Z)\sm c_1^{-1}(Z)}$ is the closure of $c_2^{-1}(Z)\sm c_1^{-1}(Z)$ in $C$, 
is the largest open subset $W$ of $C$ such that $Z$ is $c|_W$-invariant.

(c) $Z$ is $c$-invariant in a neighborhood of fixed points if and only if the closure
of $c_2^{-1}(Z)\sm c_1^{-1}(Z)$ in $C$ does not intersect with $Fix(c)$.

(d) $Z$ is locally $c$-invariant if and only if for each irreducible component 
$S$ of  $c_2^{-1}(Z)\sm c_1^{-1}(Z)$, the closures of $c_1(S)$ and  $c_2(S)$ in $X$
do not intersect.

(e) If $Z$ is locally $c$-invariant, then $Z$ is  $c$-invariant in a neighborhood of fixed points.
\end{Lem}

\begin{proof}
(a) and (b) are clear.

(c)  $Z$ is $c$-invariant in a neighborhood of fixed points if and only if the subset
$W(Z)$ from (b) contains $Fix(c)$. This implies the assertion.

(d) For an open subset $U\subset X$, the inclusion $c_1(c_2^{-1}(Z\cap U))\subset Z\cup (X\sm U)$ 
is equivalent to the inclusion 
$c_2^{-1}(Z)\sm c_1^{-1}(Z)\subset c_1^{-1}(X\sm U)\cup c_2^{-1}(X\sm U)$.
The last condition holds if and only if  
for each irreducible component $S$ of  $c_2^{-1}(Z)\sm c_1^{-1}(Z)$, we have either
$c_1(S)\subset X\sm U$ or $c_2(S)\subset X\sm U$. Thus $Z$ is locally $c$-invariant if and only if
for each $x\in Z$ and each $S$ as above, we have either $z\notin \ov{c_1(S)}$ or  
$z\notin \ov{c_2(S)}$. In other words, $Z$ is locally $c$-invariant if and only if
each $\ov{c_1(S)}\cap\ov{c_2(S)}$ does not intersect with $Z$. Since $c_2(S)$ is contained in $Z$, 
this happens if and only if each $\ov{c_1(S)}\cap\ov{c_2(S)}$ is empty.

(e) Assume that $Z$ is not  $c$-invariant in a neighborhood of fixed points. Then by (c) 
there exists an irreducible component $S$ of  $c_2^{-1}(Z)\sm c_1^{-1}(Z)$ whose 
closure $\ov{S}\subset C$ has a non-empty intersection with $Fix(c)$. Hence 
$c_1(\ov{S})\cap c_2(\ov{S})\neq\emptyset$, thus $Z$ is not
locally $c$-invariant by (d).
\end{proof}

By the following lemma, the notion of local invariance well behaves after a compactification. 

\begin{Lem} \label{L:comp}
Let $c:C\to X\times X$ be a correspondence and $U\subset X$ an open subset such that
$c_1^{-1}(U)$ is dense in $C$, $c_1|_{c_1^{-1}(U)}$ is proper, and  
$X\sm U$ is locally $c$-invariant.
Then there exists a compactification  $\ov{c}:\ov{C}\to \ov{X}\times \ov{X}$ of 
$c$ such that $\ov{X}\sm U$ is locally $\ov{c}$-invariant.
\end{Lem}
\begin{proof}
Set $Z:=X\sm U$. First we claim that there exists a compactification 
$\ov{X}$ of $X$ such that for each irreducible component $S$ 
of $c_2^{-1}(Z)\sm  c_1^{-1}(Z)$, the closures of $c_1(S)$ and $c_2(S)$ in 
$\ov{X}$ do not intersect. 

To show it, choose any compactification $\ov{X}$ of $X$. Assume that there exists $S$ such that
 $\ov{c_1(S)}\cap \ov{c_2(S)}\neq\emptyset$.
By \rl{invneib} (d), 
the intersection $\ov{c_1(S)}\cap \ov{c_2(S)}$ is contained in $\ov{X}\sm X$.
Replacing $\ov{X}$ by the blow-up $Bl_{\ov{c_1(S)}\cap \ov{c_2(S)}}(\ov{X})$,
we reach the situation where the closures of $c_1(S)$ and $c_2(S)$ in $\ov{X}$ have an empty 
intersection. Since the number of irreducible components of $c_2^{-1}(Z)\sm c_1^{-1}(Z)$ 
is finite, we achieve the goal after a finite number of blow-ups. 
 
Choose $\ov{X}$ as above and extend it to a compactification 
$\ov{c}:\ov{C}\to\ov{X}\times\ov{X}$ of $c$ (use \rr{comp}). 
We claim that this compactification satisfies the required property.

Note that $\ov{c}_1^{-1}(U)=c_1^{-1}(U)\subset C$. Indeed, since
$c_1^{-1}(U)$ is dense in $C$ and $C$ is dense in $\ov{C}$, the subset   
$c_1^{-1}(U)$ is open and dense in  $\ov{c}_1^{-1}(U)$. Therefore the required equality 
$\ov{c}_1^{-1}(U)=c_1^{-1}(U)$ follows from the fact that 
$\ov{c}_1|_{c_1^{-1}(U)}=c_1|_{c_1^{-1}(U)}$ is proper and $\ov{c}_1|_{\ov{c}_1^{-1}(U)}$ is 
separated.

By what we proved above, $\ov{c}_2^{-1}(\ov{X}\sm U)\sm\ov{c}_1^{-1}(\ov{X}\sm U)=
\ov{c}_2^{-1}(\ov{X}\sm U)\cap\ov{c}_1^{-1}(U)$ equals
$c_2^{-1}({X}\sm U)\cap c_1^{-1}(U)=c_2^{-1}(Z)\sm c_1^{-1}(Z)$. 
Hence by our assumption on $\ov{X}$ and \rl{invneib} (d),
we obtain that $\ov{X}\sm U$ is locally $\ov{c}$-invariant, as claimed. 
\end{proof}

\begin{Rem} \label{R:inv}
An important particular case of the lemma is $U=X$.
In this case, $X\sm U=\emptyset$, $\ov{X}$ can be chosen arbitrary, 
and $\ov{c}_2^{-1}(\ov{X}\sm U)\sm\ov{c}_1^{-1}(\ov{X}\sm U)=\emptyset$, thus
$\ov{X}\sm U$ is $\ov{c}$-invariant.
\end{Rem}

\begin{Emp} \label{E:restr2}
{\bf Restriction of correspondences.}
Let $c:C\to X\times X$ be a correspondence, 
$u:c_{2!}c_1^*\C{F}\to\C{F}$ a $c$-morphism, and  $Z\subset X$ 
a closed subset.

(a) Assume first that $Z$ is $c$-invariant. In this case, we denote by 
$c|_Z:c_2^{-1}(Z)_{red}\to Z\times Z$
the restriction of $c$, and by $[i_Z]$ the closed embedding of $c|_Z$ into $c$. 
Then $[i_Z]$ satisfies assumption (i) of \re{pullb}, so $u$ restricts to a $c|_Z$-morphism
$u|_Z:=[i_Z]^*(u)$.

(b) For a general $Z$, let $W=W(Z)\subset C$ be as in \rl{invneib} (b). 
Then $Z$ is $c|_W$-invariant, and we denote the correspondence 
$(c|_W)|_Z$ defined in (a) simply by $c|_Z$. Moreover, by \re{restr}, 
$u$ restricts to a $c|_W$-morphism $u|_W$, hence by (a), to a $c|_Z$-morphism 
$u|_Z:=(u|_W)|_Z$.
\end{Emp}

\begin{Emp} \label{E:naive}
{\bf Example.}
Assume that $c_2$ is quasi-finite. Then for each closed point $x$ of $X$, the correspondence 
$c|_x$ equals $c_1^{-1}(x)\cap c_2^{-1}(x)\to\{x\}\times\{x\}$. 
Also $c_1^{-1}(x)\cap c_2^{-1}(x)=Fix(c|_x)$ is finite.

For each $y\in c_1^{-1}(x)\cap c_2^{-1}(x)$, the restriction of $u|_x$ to 
$\{y\}\to\{x\}\times \{x\}$ gives an endomorphism $u_y:\C{F}_x\to\C{F}_x$.
Explicitly, $u_y$ equals the restriction of 
$u_x:(c_{2!}c_1^*\C{F})_x=\oplus_{y\in c_2^{-1}(x)}(c_1^*\C{F})_y\to \C{F}_x$
to $(c_1^*\C{F})_y=\C{F}_x$. 
In particular, the local term $LT_y(u|_x)$ equals $\Tr(u_y)$ (see \re{extr}). 
\end{Emp}

\begin{Rem}
The trace $\Tr(u_y)$ is usually called {\em the naive local term} of $u$ at $y$.
\end{Rem}

\begin{Emp} \label{E:clop}
{\bf Notation.} Let $c:C\to X\times X$ be a correspondence, $u:c_{2!}c_1^*\C{F}\to\C{F}$ 
a $c$-morphism, and $Z\subset X$ a closed  $c$-invariant subset. Then the morphism $[i_Z]$ from 
\re{restr2} is a closed embedding. Hence by \re{pushf} (a) (ii), $u|_Z$ extends to a 
$c$-morphism $[i_Z]_{!}(u|_Z)$.

Put $U:=X\sm Z$, then $c_1^{-1}(U)\subset c_2^{-1}(U)$. Hence the open embedding 
 $[j_U]$ from \re{restr} satisfies assumption (i) of \re{pushf} (a), therefore  $u|_U$ extends to a $c$-morphism $[j_U]_{!}(u|_U)$.
\end{Emp}

The following proposition, whose proof will be given in \rs{additivity}, resembles
\cite[Prop. 2.4.3]{Pi}. 

\begin{Prop} \label{P:add}
In the notation of \re{clop}, we have an equality
\[
\C{Tr}_c(u)= \C{Tr}_c([i_Z]_{!}(u|_Z))+\C{Tr}_c([j_U]_{!}(u|_U)).
\]  
\end{Prop}

\section{Main results}
\subsection{Contracting correspondences}

\begin*
\vskip 8truept
\end*

The following definition is crucial for this work.

\begin{Def} \label{D:contr}
Let $c:C\to X\times X$ be a correspondence, and $Z\subset X$ a closed subscheme.

(a) We say that $c$ {\em stabilizes} $Z$,  if 
$c_1(c_2^{-1}(Z))\subset X$ is scheme-theoretically contained in $Z$, i.e.,
if $c_1^{\cdot}(\C{I}_{Z})\subset c_2^{\cdot}(\C{I}_{Z})\cdot\C{O}_C$.

(b) We say that $c$ is {\em contracting near} $Z$, if $c$ stabilizes $Z$ 
and there exists $n\in \B{N}$ such that 
$c_1^{\cdot}(\C{I}_{Z})^n\subset c_2^{\cdot}(\C{I}_{Z})^{n+1}\cdot\C{O}_{C}$.

(c) We say that $c$ is {\em contracting near $Z$ in a neighborhood of fixed points},
if there exists an open neighborhood $W\subset C$ of $Fix(c)$ such that 
 $c|_W$ is contracting near $Z$.
\end{Def}

\begin{Rem} \label{R:contr}
(a) If  $c$ stabilizes $Z$, then the closed subset $Z_{red}\subset X$ is $c$-invariant.

(b) A correspondence $c:C\to X\times X$ stabilizes a closed subscheme $Z\subset X$ 
if and only if a closed subscheme $c^{-1}(Z\times Z)\subset C$ equals $c_2^{-1}(Z)$. 
Therefore by \rl{normfunc} (b), $c$ is contracting near $Z$ if and only if $c$ 
stabilizes $Z$ and the image of $(\wt{c}_Z)_{1s}$ is 
supported set-theoretically at the zero section $Z\subset N_{Z}(X)$. 

(c) If a correspondence $c$ is contracting near $Z$, then the corresponding rigid 
correspondence $c^{rig}$ is contracting near $Z^{rig}$ in the sense of \cite[Def. 3.1.1]{Fu}.
Furthermore, it is likely that the two notions are equivalent.
\end{Rem} 
 
Now we are ready to formulate our first main result.

\begin{Thm} \label{T:locterms}
Let $c:C\to X\times X$ be a correspondence contracting near a closed subscheme 
$Z\subset X$ in a neighborhood of fixed points, and let $\beta$ be an open connected subset of 
$Fix(c)$ such that $c'(\beta)\cap Z\neq\emptyset$. Then 

(a) $\beta$ is contained set-theoretically in $c'^{-1}(Z)$, hence $\beta$ 
is an open connected subset of $Fix(c|_Z)$.

(b) For every $c$-morphism $u:c_{2!}c_1^*\C{F}\to\C{F}$, we have
$\C{Tr}_{\beta}(u)=\C{Tr}_{\beta}(u|_Z)$. In particular, if $\beta$ is proper over $k$, 
then $LT_{\beta}(u)=LT_{\beta}(u|_Z)$.
\end{Thm}
\begin{proof}
Choose an open neighborhood $W\subset C$ of $Fix(c)$ such that $c|_W$ is contracting near $Z$.
Then $Fix(c|_W)=Fix(c)$, therefore neither assertion of the theorem 
will change if we replace $c$ by $c|_W$ 
and $u$ by its restriction $u|_W$. Hence we can assume that  $c$ is contracting near $Z$.
Moreover, replacing $C$ by an open subset $C\sm [Fix(c)\sm \beta]$, we can assume that
$Fix(c)=\beta$, hence $\C{Tr}_{\beta}=\C{Tr}_c$. For the proof, we will apply the 
construction of \re{defcor}.

(a) Since $c$ is contracting near $Z$, the image of $(\wt{c}_Z)_{1s}$ is supported at  
the zero section $Z\subset N_{Z}(X)$ (see \rr{contr} (b)).  Therefore, by \rco{diag}, 
the image of $\wt{c}'_s:N_{c'^{-1}(Z)}(\beta)\to  N_{Z}(X)$ is supported at $Z$ as well.
Hence by the ``only if'' assertion of \rl{normfunc} (b), there exists 
$n\in\B{N}$ such that $(c')^{\cdot}(\C{I}_Z)^n\subset (\C{I}_{c'^{-1}(Z)})^{n+1}$, or, equivalently,
$(\C{I}_{c'^{-1}(Z)})^{n}=(\C{I}_{c'^{-1}(Z)})^{n+1}\subset\C{O}_{\beta}$. 
Since $\beta$ is Noetherian and connected, we conclude from this that  
$(\C{I}_{c'^{-1}(Z)})^{n}=0$. Thus $\beta_{red}=c'^{-1}(Z)_{red}$,
as claimed.

(b) Put $U:=X\sm Z$. By \rp{add}, the trace 
$\C{Tr}_c(u)$ equals the sum $\C{Tr}_c([i_Z]_{!}(u|_Z))+\C{Tr}_c([j_U]_{!}(u|_U))$. By 
\rp{pushf} for the closed embedding $[i_Z]$, we have 
$\C{Tr}_c([i]_{Z!}(u|_Z))=i'_{Z!}\C{Tr}_{c|_Z}(u|_Z)$. Since $(i'_{Z})_{red}$ is the identity map
(by (a)), the map $i'_{Z!}$ is the identity as well. Thus  
it remains to show that $\C{Tr}_c([j_U]_{!}(u|_U))=0$. 
For this, we may replace $u$ by $[j_U]_{!}(u|_U)$, thus assuming that $\C{F}|_Z=0$. 
We will show that in this case, $\C{Tr}_c(u)=0$.

By \rp{spec} for correspondences $c$ and $\wt{c}_Z$, we have an equality
\[
sp_{Fix(\wt{c}_Z)}(\C{Tr}_c(u))=\C{Tr}_{(\wt{c}_Z)_s}(sp_{\wt{c}_Z}(u)).
\]
To prove the vanishing of $\C{Tr}_c(u)$ we will show that the map $sp_{\wt{c}_Z}$
vanishes, while the map $sp_{Fix(\wt{c}_Z)}$ is an isomorphism. 

As we already mentioned above, the image of $(\wt{c}_Z)_{1s}$ is supported at  
$Z\subset N_{Z}(X)$. On the other hand, by \rp{verdier} we obtain that 
$sp_{\wt{X}_Z}(\C{F})|_Z\cong\C{F}|_Z=0$. Therefore $(\wt{c}_Z)_{1s}^* sp_{\wt{X}_Z}\C{F}=0$, 
hence the specialization map $sp_{\wt{c}_{Z}}$ vanishes, as claimed.

By \re{trivdef} (a), to show that $sp_{Fix(\wt{c}_Z)}$ is an isomorphism, it will suffice to 
check that $Fix(\wt{c}_Z)_{red}$ is isomorphic to $(Fix(c)_{\C{D}})_{red}$. Recall that by 
$c'^{-1}(Z)_{red}=Fix(c)_{red}$ (by (a)), hence 
$(Fix(c)_{\C{D}})_{red}=(\wt{Fix(c)}_{c'^{-1}(Z)})_{red}$. 
Therefore by \rco{diag}, $(Fix(c)_{\C{D}})_{red}$ is a closed 
subscheme of  $Fix(\wt{c}_Z)_{red}$. Since the generic fibers of $Fix(c)_{\C{D}}$ and 
$Fix(\wt{c}_Z)$ are equals, it remains to prove that the special fiber 
$Fix((\wt{c}_Z)_s)$ is set-theoretically supported at $Fix(c)=(Fix(c)_{\C{D}})_s$.

Since $(\wt{c}_Z)_{1s}(Fix((\wt{c}_Z)_s))$  is set-theoretically supported at $Z$,
the same is true for $(\wt{c}_Z)_{2s}(Fix((\wt{c}_Z)_s))$. Hence, by
\rl{normfunc} (d), $Fix((\wt{c}_Z)_s)$ is supported at  
$c_2^{-1}(Z)\subset N_{c_2^{-1}(Z)}(C)$. But the restriction of $(\wt{c}_Z)_s$ to 
$c_2^{-1}(Z)$ equals the restriction $c|_Z$ of $c$.  
Therefore $Fix((\wt{c}_Z)_s)$ is supported at  
$Fix(c)$, as claimed.

This completes the proof of the theorem.
\end{proof} 

\subsection{Correspondences over finite fields} \label{SS:finite}

\begin*
\vskip 8truept
\end*

\begin{Emp} \label{E:frtw}
{\bf Twisting of correspondences.}
(a) For a  scheme $X$ over $\B{F}$ defined over $\fq$, we 
denote by $\Fr_{X,q}$ the geometric Frobenius morphism $X\to X$ over $\fq$.

(b) For a correspondence $c:C\to X_1\times X_2$ defined over $\fq$ and 
$n\in\B{N}$, we denote by $c^{(n)}$ the correspondence $(c_1^{(n)},c_2):C\to X_1\times X_2$,
where $c_1^{(n)}$ is $\Fr_{X_1,q}^n\circ c_1=c_1\circ \Fr_{C,q}^n$.
\end{Emp}

\begin{Not} \label{N:degree}
Let $f:Y\to X$ be a morphism of Noetherian schemes over $k$.

(a) Let $Z$ be a closed subset of $X$. Then  
$\C{I}_{(f^{-1}(Z))_{red}}$ is the radical of $\C{I}_{f^{-1}(Z)}$. Since $Y$ is Noetherian,
there exists a positive integer $m$ such that 
$(\C{I}_{(f^{-1}(Z))_{red}})^m\subset\C{I}_{f^{-1}(Z)}$.
The smallest such an $m$ we call {\em the ramification of $f$ at $Z$} and denote by 
$ram(f,Z)$.

(b) Assume that $f$ is quasi-finite. Then for each closed point $x\in X$, the ramification 
degree $ram(f,x)$ is at most $\dim_k k[f^{-1}(x)]$. In particular, the set $\{ram(f,x)\}_x$ 
is bounded. We denote  by $ram(f)$ the maximum of the $ram(f,x)$'s and call it 
{\em the ramification degree of $f$}.
\end{Not}

\begin{Lem} \label{L:contrac}
Let $c:C\to X\times X$ be a correspondence over $\B{F}$ defined over $\fq$, let $Z$ be 
a closed subset of $X$, and let $n\in\B{N}$ be such that $q^n>ram(c_2,Z)$ and 
$Z$ is $c^{(n)}$-invariant. Then the correspondence $c^{(n)}$ is contracting near $Z$.
\end{Lem}

\begin{proof}
Denote $ram(c_2,Z)$ by $d$, and let $\varphi_{X,q}$ be the arithmetic Frobenius isomorphism 
$X\isom X$ over $\fq$. Then for every section $f$ of $\C{O}_{X}$, we have 
$\Fr_{X,q}^{\cdot}(f)=(\varphi_{X,q}^{-1})^{\cdot}(f)^q$. 
Therefore $(c_1^{(n)})^{\cdot}(\C{I}_Z)$ equals 
$(\varphi^{-n}_{X,q})^{\cdot} c_1^{\cdot}(\C{I}_Z)^{q^n}$. 

Since $Z$ is $c^{(n)}$-invariant, the sheaf 
$(\varphi^{-n}_{X,q})^{\cdot} c_1^{\cdot}(\C{I}_Z)^{q^n}$ 
is therefore contained in $\C{I}_{c_2^{-1}(Z)_{red}}$. Hence 
$(\varphi^{-n}_{X,q})^{\cdot} c_1^{\cdot}(\C{I}_Z)
\subset\C{I}_{c_2^{-1}(Z)_{red}}$, thus 
$(c_1^{(n)})^{\cdot}(\C{I}_Z)\subset (\C{I}_{c_2^{-1}(Z)_{red}})^{q^n}$.
As  $q^n\geq d+1$, we conclude that 
$(c_1^{(n)})^{\cdot}(\C{I}_Z)\subset (\C{I}_{c_2^{-1}(Z)_{red}})^{d+1}
\subset \C{I}_{c_2^{-1}(Z)}$,
implying that $c^{(n)}$ stabilizes $Z$. Furthermore, $(c_1^{(n)})^{\cdot}(\C{I}_Z)^d$ is contained in 
$(\C{I}_{c_2^{-1}(Z)_{red}})^{d(d+1)}
\subset(\C{I}_{c_2^{-1}(Z)})^{d+1}$, hence $c^{(n)}$ is contracting near $Z$, as claimed.
\end{proof}

\begin{Cor} \label{C:twisting}
Let $c:C\to X\times X$ be a correspondence over $\B{F}$ defined over $\fq$.

(a) Let $Z$ be a closed  locally $c$-invariant subset of $X$ defined over $\fq$. 
Then for each  $n\in\B{N}$ with $q^n>ram(c_2,Z)$ the correspondence $c^{(n)}$ is contracting near 
$Z$ in a neighborhood of fixed points.

(b) If $c_2$ is quasi-finite, then for each  $n\in\B{N}$ with $q^n>ram(c_2)$, 
the correspondence $c^{(n)}$ is contracting near every closed point $x$ of $X$ 
in a neighborhood of fixed points.
\end{Cor}

\begin{proof}

(a) Since $Z$ is  locally $c$-invariant and defined over $\fq$, it is 
locally $c^{(n)}$-invariant. Hence by \rl{invneib} (c) and (e), 
the open subset $W=W(Z)\subset C$ from \rl{invneib} (b) contains $Fix(c^{(n)})$, 
and $Z$ is $c^{(n)}|_W$-invariant. Then $c^{(n)}|_W=(c|_W)^{(n)}$ and 
$ram(c_2|_W,Z)\leq ram(c_2,Z)$, hence the assertion follows from 
\rl{contrac} for the correspondence $c|_W$.
 
(b) As $c_2$  is quasi-finite, every closed point $x$ of $X$ is 
locally $c^{(n)}$-invariant (see \re{qfinite}).
Thus, as in (a), the assertion follows from \rl{contrac}.
\end{proof}

\subsection{Generalization of a  theorem of Fujiwara} \label{SS:fuj}

\begin*
\vskip 8truept
\end*

To formulate our main result, we need the following generalization of construction of 
\re{rgm} (a).

\begin{Emp} \label{E:rgm'}
Let $c:C\to X_1\times X_2$ be a correspondence, $u:c_{2!}c_1^*\C{F}_1\to\C{F}_2$ a 
$c$-morphism, and $U_1\subset X_1$ an open subset such that 
$c_1|_{c_1^{-1}(U_1)}:c_1^{-1}(U_1)\to U_1$ is proper and $\C{F}_1|_{X_1\sm U_1}=0$.
Then $u$ gives rise to a morphism 
\[
R\Gm_c(u):R\Gm_c(X_1,\C{F}_1)\to R\Gm_c(X_2,\C{F}_2),
\]
defined as follows.
 
Denote by $c^0:c_1^{-1}(U_1)\to U_1\times X_2$ the restriction of $c$, and by  
$u^0$ the restriction of $u$ to $c^0$. Since $c^0_1$ is proper, the construction of 
\re{rgm} (a) gives rise to a homomorphism 
\[
R\Gm_c(u^0): R\Gm_c(U_1,\C{F}_1|_{U_1})\to R\Gm_c(X_2,\C{F}_2).
\]
As  $\C{F}_1|_{X_1\sm U_1}=0$, the canonical map 
$j_!:R\Gm_c(U_1,\C{F}_1|_{U_1})\to R\Gm_c(X_1,\C{F}_1)$ is an isomorphism. We define
$R\Gm_c(u)$ to be the composition $R\Gm_c(u^0)\circ(j_!)^{-1}$.
\end{Emp}

Now we are ready to formulate a generalization of the theorem of Fujiwara, 
which was suggested to me by David Kazhdan and is crucial for our work \cite{KV}.

\begin{Thm} \label{T:Del} 
Let $c:C\to X\times X$ be a correspondence over $\B{F}$ defined over $\fq$.

(a) Assume that $c_2$ is quasi-finite. Then for every $n\in\B{N}$ with $q^n>ram(c_2)$, 
the set $Fix(c^{(n)})$ is finite. 

(b) Let  $U\subset X$ be an open subset defined over $\fq$ such that 
$c_1|_{c_1^{-1}(U)}$ is proper, $c_2|_{c_2^{-1}(U)}$ is quasi-finite,
and $X\sm U$ is locally $c$-invariant.

Then there exists a positive integer $d\geq ram(c_2|_{c_2^{-1}(U)})$ with 
the following property: For every  
$\C{F}\in D_{ctf}^b(X,\La)$ with $\C{F}|_{X\sm U}=0$, 
every $n\in\B{N}$ with $q^n>d$ and every  $c^{(n)}$-morphism 
$u: c_{2!}(c^{(n)}_1)^*\C{F}\to\C{F}$, we have an equality
\begin{equation} \label{Eq:fuj}
\Tr(R\Gm_c(u))=\sum_{y\in Fix(c^{(n)})\cap c_2^{-1}(U)}\Tr(u_y).
\end{equation}

(c) In the notation of (b), assume that $X$ and $C$ are proper over $\B{F}$. Then \\  
$d:=max\{ram(c_2|_{c_2^{-1}(U)}),ram(c_2,X\sm U)\}$ satisfies the conclusion
of (b). 
\end{Thm}

\begin{Rem} 
(a) Both sides of \form{fuj} are well-defined. Namely, 
$R\Gm_c(u)$ was defined in \re{rgm'}, $u_y$ was defined in \re{naive}, 
and the sum is finite by \ref{T:Del} (a).

(b) The constant $d$ in \ref{T:Del} (b) can be explicitly estimated. 
Namely, one can see from the proof that the picture can be compactified, and then 
\ref{T:Del} (c) gives an estimate for $d$.

(c) In the notation of \ref{T:Del} (b), assume that $\C{F}\in D_{ctf}^b(X,\La)$ is 
equipped with a morphism $\psi:\Fr_{X,q}^*\C{F}\to\C{F}$ (say, $\C{F}$ is a Weil sheaf) 
and with a $c$-morphism $u:c_{2!}c_1^*\C{F}\to\C{F}$.  Then for each $n\in\B{N}$,
 $\C{F}$ is equipped with a $c^{(n)}$-morphism 
$u^{(n)}:=u\circ \psi^n:c_{2!}(c^{(n)}_1)^*\C{F}\to\C{F}$, so one can apply formula 
\form{fuj}. In the case $U=X$, the assertion thus reduces to Deligne's conjecture proven by 
Fujiwara \cite{Fu}.
\end{Rem}

\begin{proof}
(a) We have to show that each connected component $\beta\in \pi_0(Fix(c^{(n)}))$ is a point.
Pick a closed point $x\in (c^{(n)})'(\beta)$. By \rco{twisting} (b), the correspondence 
$c^{(n)}$ is contracting 
near $x$ in a neighborhood of fixed points. Hence by \rt{locterms} (a),
$\beta$ is a connected component of a finite scheme $Fix(c|_x)$ (see \re{naive}).  
Therefore $\beta$ is a point, as claimed. 

(c) Fix $n\in\B{N}$ with $q^n>d$. By the Lefschetz-Verdier trace formula (\rco{LTF}), 
we have an equality
\begin{equation} \label{Eq:lef}
\Tr(R\Gm_c(u))=\sum_{\beta\in\pi_0(Fix(c^{(n)}))}LT_{\beta}(u).
\end{equation}
Pick any $\beta\in \pi_0(Fix(c^{(n)}))$ such that $c_2(\beta)\subset X\sm U$. 
By \rco{twisting} (a), the correspondence 
$c^{(n)}$ is contracting near $X\sm U$ in a neighborhood of fixed points. Therefore by 
\rt{locterms}, $\beta$ is a connected component of $Fix(c|_{X\sm U}^{(n)})$, and  
$LT_{\beta}(u)$ equals $LT_{\beta}(u|_{X\sm U})$. 
However, $\C{F}|_{X\sm U}=0$, hence $u|_{X\sm U}$ vanishes. Thus
$LT_{\beta}(u)=LT_{\beta}(u|_{X\sm U})$ vanishes as well.

Pick now any  $\beta\in \pi_0(Fix(c^{(n)}))$ such that 
$c_2(\beta)\cap U\neq\emptyset$. By (a), such $\beta$ is simply a point 
$y\in Fix(c^{(n)})\cap c_2^{-1}(U)$, 
while by \rco{twisting} (b),  the correspondence 
$c^{(n)}$ is contracting near $c_2(y)$ in a neighborhood of fixed points. Hence 
by \rt{locterms} (b),  $LT_{\beta}(u)$ equals  $LT_y(u|_x)$. Thus by 
\re{naive}, it equals $\Tr(u_y)$.
This shows that the right-hand side of \form{lef} is equal to that of \form{fuj}, as claimed.

(b) For the proof we can replace $c$ and $u$ by their restrictions to  $c_1^{-1}(U)$. 
Then the assumptions of \rl{comp} are satisfied, hence there exists a compactification 
$\ov{c}:\ov{C}\to \ov{X}\times \ov{X}$ of $c$ such that $\ov{X}\sm U$ is 
locally $\ov{c}$-invariant. 

Denote by $c^0:c_1^{-1}(U)\to U\times X$ the restriction of $c$, by $u^0$ the restriction of $u$ 
to $c^0$, and by $[j]=(j^0,j_{C^0},j)$ the inclusion map of 
$c^0$ into $\ov{c}$. Since $c^0_1$ is proper, $[j]$ satisfies assumption (iii) of \re{pushf} (a). 
Therefore $u^0$ extends to a $\ov{c}$-morphism 
$\ov{u}:=[j]_!(u^0):\ov{c}_{2!}\ov{c}_1^*\ov{\C{F}}\to\ov{\C{F}}$, where 
$\ov{\C{F}}:=j_!\C{F}=j^0_!(\C{F}|_U)$. 

We claim that the equality \form{fuj} for $c, U$ and $u$ is equivalent to that 
for $\ov{c}, U$ and $\ov{u}$. Indeed, the equality 
$\Tr(R\Gm_c(\ov{u}))=\Tr(R\Gm_c(u))$ follows from \re{rgm} (c), while 
the equality $\Tr(\ov{u}_y)=\Tr(u_y)$ for each 
$y\in Fix(c^{(n)})\cap c_2^{-1}(U)$ follows from the definition.  
This shows that assertion (b) is a consequence of (c). 

This completes the proof of the theorem.
\end{proof}
\section{Theorem of Verdier} \label{S:verdier}

The goal of this section is to prove \rp{verdier}.  
Set $U:=X\sm Z$, let $i:Z\hra X$ and $j:U\hra X$ be the natural embeddings, and 
denote by $C_{X,Z}(\C{F})$ the cone of the morphism 
$sp_{\wt{X}_Z}(\C{F})|_Z\to\C{F}|_Z$ from \form{verdier}. 
We claim that $C_{X,Z}(\C{F})=0$ for all $X$, $Z$ and all 
$\C{F}\in D^b_c(X,\La)$.

\subsection{Contractible case} \label{SS:partcase}

\begin*
\vskip 8truept
\end*

In this subsection we will prove \rp{verdier} in the case 
$X=\B{A}^1$, $Z:=\{x\in\B{A}^1|x^k=0\}$ and $\C{F}=\La_X$. 
Though this case can be easily done by direct calculation, we will deduce 
it from the fact that $X$ is contractible to $Z$.
(This approach was suggested to us by David Kazhdan).

\begin{Def} \label{D:contract}
We say that a scheme $X$ is {\em contractible to its closed subcheme $Z$},  
if there exists a morphism $H:X\times\B{A}^1\to X$ such that 

(i) $H|_{X\times\{1\}}:X\to X$ is the identity;

(ii) $H(X\times\{0\})$ and  $H(Z\times\B{A}^1)$ are scheme-theoretically contained in $Z$; 

(iii) the restriction $H|_{Z_{red}\times\B{A}^1}$ is the projection 
$\pi:Z_{red}\times \B{A}^1\to Z_{red}\subset X$.
\end{Def}

\begin{Ex} \label{Ex:basic}
Let $X=\B{A}^1$ and $Z:=\{x\in\B{A}^1|x^k=0\}$. Hence $X$ is contractible to $Z$ via  
the contraction map $H:(x,a)\mapsto xa$.
\end{Ex}

\begin{Lem} \label{L:contract}
\rp{verdier} holds, if $X$ is contractible to $Z$ and 
$\C{F}=\La$.
\end{Lem}

\begin{proof}
For each $a\in\B{A}^1$, denote by $i_a$ the inclusion 
$X\cong X\times\{a\}\hra X\times\B{A}^1$, and set $H_a:=H\circ i_a:X\to X$.
Then $H_1=\Id_X$, and $H_0$ factors through $i:Z\hra X$.  Since 
$H_a(Z)$ is schematically contained in $Z$ (by (ii)), $H_a$ lifts to a morphism
$\wt{H}_a:\wt{X}_Z\to \wt{X}_Z$ (by \rl{normfunc} (a)). Since 
$(\wt{H}_a)_s|_{Z_{red}}$ is the identity (by (iii)), the base change morphism $BC^*$ of 
\re{specfunc} (b) corresponding to $H_a$ induces a morphism 
\[
\phi_a:sp_{\wt{X}_Z}(\La)|_Z=[(\wt{H}_a)_s^*sp_{\wt{X}_Z}(\La)]|_Z \to sp_{\wt{X}_Z}(\La)|_Z.
\]
Then $\phi_1$ is the identity map, while $\phi_0$ factors through the morphism 
$sp_{\wt{X}_Z}(\La)|_Z\to\La_Z$ from \form{verdier}. Thus it will suffice to show that $\phi_a$ is independent of $a\in\B{A}^1$.

Consider the endomorphism 
\[
\phi:\pi^* (sp_{\wt{X}_Z}(\La)|_Z)=[H^*sp_{\wt{X}_Z}(\La)]|_{Z\times \B{A}^1}
\to [sp_{\wt{X\times\B{A}^1}_{Z\times\B{A}^1}}(\La)]|_{Z\times \B{A}^1}
\isom\pi^* (sp_{\wt{X}_Z}(\La)|_Z)
\]
of $\pi^*(sp_{\wt{X}_Z}(\La)|_Z)\in D_{c}^{b}(Z,\La)$, obtained as a composition of the 
the base change morphism  $BC^*$ of \re{specfunc} (b) 
corresponding to $H$ and the inverse of the base change morphism  $BC^*$ corresponding to
the (smooth) projection $\pi:X\times\B{A}^1\to X$. 

Since $H_a=H\circ i_a$, $\pi\circ i_a=\Id_X$, and base change morphisms are 
compatible with compositions, the fiber of $\phi$ over each $a\in\B{A}^1$
equals $\phi_a$. Therefore the assertion follows from \rl{families} below.
\end{proof}

\begin{Lem} \label{L:families}
Let $X$ and $Y$ be schemes over $k$ such that $X$ is connected. 
For every $\C{A},\C{B}\in D_{c}^b(Y,\La)$ and every  
$\phi\in \Hom(\La\pp\C{A},\La\pp\C{B})$, its restriction 
$\phi_x\in\Hom(\C{A},\C{B})$ to $\{x\}\times X=X$ is independent of $x\in X$.
\end{Lem}  

\begin{proof}
Note that $\Hom(\La\pp\C{A},\La\pp\C{B})$ is canonically isomorphic to 
\[
H^0(R\Gm(X\times Y,\C{RHom}(\La\pp\C{A},\La\pp\C{B}))\cong 
H^0(R\Gm(X,\La)\otimes R\Gm(Y,\C{RHom}(\C{A},\C{B}))),
\]
while $\Hom(\C{A},\C{B})\cong H^0(R\Gm(Y,\C{RHom}(\C{A},\C{B})))$. Moreover, the restriction
$\phi\mapsto\phi_x$ is induced by the pullback $i_x^*:R\Gm(X,\La)\to R\Gm(\{x\},\La)=\La$
corresponding to the embedding $i_x:\{x\}\hra X$. Since $X$ is connected, the pullback 
$i_x^*$ is independent of $x\in X$, as claimed.
\end{proof}

\subsection{General case} \label{SS:gencase}

\begin*
\vskip 8truept
\end*

To deduce the general case from the one considered in Subsection \ref{SS:partcase}, 
we follow Verdier \cite[pp. 354--355]{Ve}.
 
\begin{Lem} \label{L:red}
 
Let $f:X'\to X$ be a morphism of schemes over $k$, $Z\subset X$ 
a closed subscheme, $Z'=f^{-1}(Z)$, and let 
$f_Z:Z'\to Z$  be the restriction of $f$.

(a) If $f$ is proper, then for each $\C{F}'\in D_c^b(X',\La)$, we have 
$C_{X,Z}(f_!\C{F}')\cong f_{Z!}(C_{X',Z'}(\C{F}'))$.

(b) If $f$ is smooth, then for each $\C{F}\in D_c^b(X,\La)$, we have
$C_{X',Z'}(f^*\C{F})\cong f^*_{Z}(C_{X,Z}(\C{F}))$.
\end{Lem}

\begin{proof} 
Let $\wt{f}:\wt{X'}_{Z'}\to\wt{X}_{Z}$ the morphism lifting $f$ (use \rl{normfunc} (a)).  
If $f$ is proper (resp. smooth), then $\wt{f}$ is proper  (resp. smooth) as well 
(by \rl{normfunc} (e)). Then the assertion  follows from the proper  (resp. smooth) 
base change theorem.
\end{proof}

\begin{Emp} \label{E:step1}
{\bf Reduction steps}. 

(I) By construction, the morphism \form{verdier} is an isomorphism when $Z=X$ 
(compare \rr{trivdef} below). Hence by \rl{red} (a) for the closed embedding $i:Z\hra X$, 
we obtain $C_{X,Z}(i_!i^*\C{F})=C_{Z,Z}(i^*\C{F})=0$. Therefore 
for each $\C{F}\in  D^b_c(X,\La)$, the canonical morphism $C(j_!j^*\C{F})\to C(\C{F})$
is an isomorphism. Hence it will suffice to show the assertion under the additional 
assumption $\C{F}|_Z=0$.

(II) Let $f$ be the blow-up $X':=Bl_Z(X)\to X$. Since $f$ induces an isomorphism  
$X'\sm f^{-1}(Z)\isom X\sm Z$, while $\C{F}|_Z=0$, we have $\C{F}=f_!f^*\C{F}$. 
Thus by \rl{red} (a) for $f$, the assertion for $(X,Z,\C{F})$ follows from that for 
$(X',f^{-1}(Z),f^*\C{F})$. In particular, we can assume that $Z\subset X$ is a 
Cartier divisor.

(III) Since the assertion is local, we can assume that $X$ is affine and 
$Z\subset X$ is given by one equation (by (II)). Then there exists a closed embedding 
$f:X\hra\B{A}^n$ such that $Z:=f^{-1}(\B{A}^{n-1}\times\{0\})$. 
By \rl{red} (a) for $f$, 
the assertion for $(X,Z,\C{F})$ is equivalent to that for 
$(\B{A}^n,\B{A}^{n-1}\times\{0\},f_!\C{F})$. Thus it is enough to prove the assertion for 
$X=\B{A}^n$ and $Z=\B{A}^{n-1}\times\{0\}$.
\end{Emp}

\begin{Cl} \label{C:generic}
Assume that $X$ is normal, 
$Z\subset X$ is given by one equation $f=0$ and $\C{F}|_Z=0$. Then 
$C_{X,Z}(\C{F})$ vanishes at each generic point of $Z$.   
\end{Cl}
\begin{proof}
Since the assertion is local, we may assume that both $X$ and $Z$ are irreducible.

 Assume first that  $\C{F}=j_!\La_{U}$.
Since $X$ is normal, we can replace it by open subset (containing the generic point of $Z$)
such that $Z_{red}$ is smooth. Further shrinking $X$, we can assume that exists a 
smooth morphism $f:X\to \B{A}^1$ such that $Z_{red}=f^{-1}(0)$. Then there exists 
$k\in\B{N}$ such that $Z=f^{-1}(\{x^k=0\})$. The vanishing of $C_{X,Z}(\La_X)$ then 
follows from \rl{red} (b), \rl{contract} and Example \ref{Ex:basic}.  Hence 
$C_{X,Z}(j_!\La_{U})=C_{X,Z}(\La_X)=0$ by \re{step1} (I).

 For the general case, we may assume that $\La$ is finite. 
Next, since $\C{F}\mapsto C(\C{F})$ is a triangulated functor, we can assume that 
$\C{F}$ is a sheaf (and not a complex of sheaves). Then replacing $X$ by an open subset 
(containing the generic point of $Z$) we may assume that 
$\C{F}=j_!\C{G}$  for some local system $\C{G}$ on $U$. We will show by 
induction on $k$ that the cohomology sheaf $\C{H}^k(C_{X,Z}(\C{F}))$ vanishes generically.

When $k<0$, the assertion is clear. For the induction step, choose a finite 
\'etale covering $f:U'\to U$ such that the local system $f^*\C{G}$ is trivial, and let
$\ov{f}:X'\to X$ be the normalization of $X$ extending $f$. 
Denote $\C{G}'$ the cokernel of the natural embedding $\C{G}\hra f_!f^*\C{G}$. Then
$\C{G}'$ is a local system and $j_!\C{G}'$ is the cokernel of the embedding
$\C{F}\hra \ov{f}_!\ov{f}^*\C{F}$. 

By  what we proved above and \rl{red} (a), 
$C_{X,Z}(\ov{f}_!\ov{f}^*\C{F})$ vanishes generically. 
Hence by the long exact sequence for the cohomology we see that 
$\C{H}^k(C_{X,Z}(\C{F}))$ is generically isomorphic to   $\C{H}^{k-1}(C_{X,Z}(j_!\C{G}'))$.
Hence it vanishes by the induction hypothesis.
\end{proof}

\begin{Emp} \label{E:compl}
{\bf Completion of the proof}.
By \re{step1}, it will suffice to prove the assertion for 
$X=\B{A}^n$ and $Z=\B{A}^{n-1}\times\{0\}$. We will prove the assertion by induction on $n$.

If $n=1$, then $Z$ is a point, so the assertion follows from \rcl{generic}.

Assume now that $n>1$. Let $Y\subset Z$ be the closure of the support of 
$C_{X,Z}(\C{F})$. By \rcl{generic}, $Y\neq Z$. Hence by Noether normalization theorem, 
there exists a line $l\subset Z\subset X$ such that the restriction to $Y$ of the 
linear projection $q:X\to X':=X/l$ is finite. 

Put $Z':=Z/l\subset X'$ , let $\ov{q}:\ov{X}\to X'$ be the 
compactification of $q$, set $\ov{Z}:=\ov{q}^{-1}(Z')$ and let 
$\ov{\C{F}}\in D_c^{b}(\ov{X},\La)$ be the extension of $\C{F}$ by zero. 
Then it follows from \rl{red} (a) and the induction hypothesis that 
$\ov{q}_!(C_{\ov{X},\ov{Z}}(\ov{\C{F}}))\cong C_{X',Z'}(\ov{q}_!\ov{\C{F}})=0$. 
On the other hand,  $C_{\ov{X},\ov{Z}}(\ov{\C{F}})$ is supported on 
$\ov{Y}:=Y\cup(\ov{Z}\sm Z)$, while the restriction $\ov{q}|_{\ov{Y}}$ is finite.
Therefore $C_{\ov{X},\ov{Z}}(\ov{\C{F}})=0$, hence 
$C_{X,Z}(\C{F})=0$.
\end{Emp}

\begin{Rem} \label{R:trivdef}
Similar (but simpler) arguments can be used to prove that  for all $\C{F}\in D_c^b(X,\La)$, 
the map $\C{F}\to sp_{\wt{X}}\C{F}$ from \re{trivdef} (a) is an isomorphism.
(This fact was implicitly used in the definition of the map \form{verdier}). 

Denote the cone of $\C{F}\to sp_{\wt{X}}\C{F}$ by $C_X(\C{F})$. For a morphism $f:X'\to X$
we have $C_{X}\circ f_!\cong f_{!}\circ C_{X'}$, if $f$ is proper, and 
$C_{X'}\circ f^*\cong f^*\circ C_X$, if $f$ is smooth. In particular, as in 
\re{step1} (III), we reduce to the case $X=\B{A}^n$.

To show that  $C_X(\C{F})$ vanishes generically, we may assume that $\C{F}$ is a sheaf.
Next passing to an open subset of $X$, we may assume that $X$ is smooth and $\C{F}$ is a local system.
Then passing to an \'etale cover of $X$, we may assume that $\C{F}=\La_X$. 
Since $X$ is smooth, we are reduced to the case $X=\pt$, in which case it is standard.
The rest of the proof goes as in \re{compl} word-by-word. 
\end{Rem}

\section{Functorial properties of trace maps} \label{S:funct}

In this section we will prove a result that generalizes both \rp{pushf}
and \rp{spec}. 

\subsection{Cohomological morphisms} \label{SS:cohmor}

\begin*
\vskip 8truept
\end*

\begin{Emp} \label{E:cohdatum}
{\bf Cohomological pre-morphisms.} (a) Let $c:C\to X\times X$ be a correspondence.
Then $c$ gives rise to a diagram 
\begin{equation} \label{Eq:diag}
\CD 
    Fix(c) @>c'>> X\\
      @V{\Dt'}VV @V{\Dt}VV\\
      C@>c>>X\times X@>p_i, i=1,2>>X@>\pi_X>>\pt(=\Spec k), 
\endCD 
\end{equation}
and we denote by  $Ob(c)$ and $Mor(c)$ the sets 
$\{X,X\times X,C,Fix(c),\pt\}$ and 
\[
\{(\Id_Z,\pi_Z)_{Z\in \Ob(c)}; c,c',\Dt,\Dt';p_i, 
c_i:=p_i\circ c, (i=1,2); c'':=c\circ\Dt'=\Dt\circ c\},
\] 
respectively. In other words, 
 $Ob(c)$ and $Mor(c)$ are the sets of objects and morphisms ``appearing in \form{diag}''.

(b) Let  $c:C\to X\times X$ and $\ov{c}:\ov{C}\to \ov{X}\times\ov{X}$ be 
a pair of correspondences. We denote by $\bar{\cdot}$ the natural 
bijections $Ob(c)\isom Ob(\ov{c})$ and  $Mor(c)\isom Mor(\ov{c})$.
By a {\em cohomological pre-morphism} from $c$ to $\ov{c}$ we mean a collection 
\[
(\{f_Z,t_Z\}_{Z\in Ob(c)};\{BC^*_{g},BC^!_{g}\}_{g:Z_1\to Z_2\in Mor(c)};\iota),
\] 
where

\indent\indent\textbullet\: {$f_Z$} is a functor $D_{ctf}^b(Z,\La)\to D_{ctf}^b(\ov{Z},\La)$;

\indent\indent\textbullet\:  {$t_Z$} is a morphism of functors $f_Z\C{A}\otimes f_Z\C{B}\to 
f_Z(\C{A}\otimes \C{B})$;

\indent\indent\textbullet\:  {$BC^*_{g}$} is a morphism of functors $\ov{g}^* f_{Z_2}\to f_{Z_1}g^*$;

\indent\indent\textbullet\:  {$BC^!_{g}$} is a morphism of functors $f_{Z_1}g^!\to\ov{g}^! f_{Z_2}$;

\indent\indent\textbullet\: {$\iota$}  is an isomorphism  $f_{\pt}\La_{\pt}\isom\La_{\ov{\pt}}$.

(c) For  each $Z\in\Ob(c)$ and  $g:Z_1\to Z_2\in\Mor(c)$, 
cohomological pre-morphism $(f_Z,t_Z;BC^*_{g},BC^!_{g};\iota)$ 
gives rise to:

\indent\indent(i) morphism 
$\wt{\pi}_Z:f_Z K_Z=f_Z(\pi_Z^!\La_{\pt})\overset{BC^!_{g}}{\lra}\pi_{\ov{Z}}^!f_{\pt}\La_{\pt}
\overset{\iota}{\lra} \pi_{\ov{Z}}^!\La_{\ov{\pt}}=K_{\ov{Z}}$;

\indent\indent(ii) morphism $r_Z:f_Z\C{RHom}(\C{A},\C{B})\to \C{RHom}(f_Z\C{A},f_Z\C{B})$, 
adjoint to 
\[
f_Z\C{RHom}(\C{A},\C{B})\otimes f_Z\C{A}\overset{t_Z}{\lra}
f_Z(\C{RHom}(\C{A},\C{B})\otimes \C{A}) \overset{ev}{\lra}f_Z\C{B};
\]

\indent\indent(iii) morphism  $d_Z:f_Z\B{D}\C{A}\overset{r_{Z}}{\lra}\C{RHom}(f_Z\C{A},f_Z K_Z)
\overset{\wt{\pi}_Z}{\lra}\B{D} f_Z\C{A}$;

\indent\indent(iv) morphism
$BC_{*g}: f_{Z_2} g_*\to \ov{g}_* f_{Z_1}$, adjoint to the composition
\[
\ov{g}^* f_{Z_2} g_*\overset{BC^*_g}{\lra} f_{Z_1}g^*g_*\overset{adj}{\lra} f_{Z_1};
\]

\indent\indent(v) morphism $BC_{!g}:\ov{g}_! f_{Z_2}\to f_{Z_1} g_!$, adjoint to the composition
\[
f_{Z_2}\overset{adj}{\lra}f_{Z_2}g^!g_!\overset{BC^!_g}{\lra} \ov{g}^!f_{Z_1} g_!.
\]
(d) If $(f_Z,t_Z;BC^*_{g},BC^!_{g};\iota)$ is a cohomological pre-morphism 
from $c$ to $\ov{c}$, and\\ $(f_{\ov{Z}},t_{\ov{Z}};BC^*_{\ov{g}},BC^!_{\ov{g}};\ov{\iota})$ 
is a cohomological pre-morphism from $c$ to $\ov{c}$, then the composition
 $(f_{\ov{Z}}\circ f_Z,t_Z\circ t_{\ov{Z}}; BC^*_{g}\circ BC^*_{\ov{g}}, BC^!_{\ov{g}}\circ BC^!_{g};
\ov{\iota}\circ\iota)$ is a cohomological pre-morphism 
from $c$ to $\ov{\ov{c}}$.

\end{Emp}

\begin{Not}
(a) Let $c^0$ and $\ov{c}^0$ be open subcorrespondences of $c$ and $\ov{c}$, respectively.
For each $Z\in\Ob(c)$ (resp. $g\in\Mor(c)$) we denote by $Z^0$ (resp $g^0$) the 
corresponding element of $\Ob(c^0)$ (resp. $\Mor(c^0)$).

(b) We say that a cohomological pre-morphism $(f_Z,t_Z;BC^*_{g},BC^!_{g};\iota)$  
from $c$ to $\ov{c}$ {\em extends} a cohomological pre-morphism
 $(f_{Z^0},t_{Z^0};BC^*_{g^0},BC^!_{g^0};\iota)$ from $c^0$ to $\ov{c}^0$,  
if 

\indent\indent(i) for each $Z\in\Ob(c)$ and $\C{A}\in D^b_{ctf}(Z,\La)$, 
we are given an isomorphism 
\indent\indent\indent\;$f_{Z^0}(\C{A}|_{Z^0})\isom f_{Z}(\C{A})|_{\ov{Z}^0}$;

\indent\indent(ii) the isomorphisms of (i) identify $t_{Z^0}$ with  
$t_Z|_{\ov{Z}^0}$, $BC^*_{g^0}$ with $BC^*_{g}|_{\ov{Z}_1^0}$, and  
\indent\indent\;\;\;\;$BC^!_{g^0}$ with $BC^!_{g}|_{\ov{Z}_1^0}$.
\end{Not}

\begin{Def} \label{D:cohmor}
(a) A cohomological pre-morphism $(f_Z,t_Z;BC^*_{g},BC^!_{g};\iota)$ from  $c$ to $\ov{c}$
is called a {\em compactifiable cohomological morphism}, if it satisfies the following axioms:

{(I)} Each $t_Z$ is commutative and associative.

{(II)} The composition 
$f_{\pt}\La_{\pt}\otimes f_{\pt}\La_{\pt}\overset{t_{\pt}}{\lra} 
f_{\pt}(\La_{\pt}\otimes\La_{\pt})=f_{\pt}\La_{\pt}
\overset{\iota}{\lra}\La_{\ov{\pt}}$ equals 
$\iota\otimes\iota:f_{\pt}\La_{\pt}\otimes f_{\pt}\La_{\pt}\to
\La_{\ov{\pt}}\otimes\La_{\ov{\pt}}=\La_{\ov{\pt}}$.

{(III)} Each 
$BC^*_{\Id_Z}=\Id_{f_Z}$, and for each $g_1:Z_1\to Z_2$ and $g_2:Z_2\to Z_3$
in $\Mor(c)$, the morphism $BC^*_{g_2\circ g_1}$ decomposes as
$
\ov{g}^*_{1}\ov{g}^*_{2}f_{Z_3}\overset{BC^*_{g_2}}{\lra}\ov{g}^*_1f_{Z_2}g^*_2
\overset{BC^*_{g_1}}{\lra}f_{Z_1}g^*_{1}g^*_{2}.
$

{(IV)} For each $g:Z_1\to Z_2\in\Mor(c)$ and $\C{A},\C{B}\in  D_{ctf}^b(Z_2,\La)$,
the following diagram is commutative
\[
\CD
\ov{g}^*f_{Z_2}\C{A}\otimes \ov{g}^*f_{Z_2}\C{B}  @>t_{Z_2}>>
\ov{g}^*f_{Z_2}(\C{A}\otimes \C{B})\\
@VBC^*_{g}\otimes BC^*_{g}VV @VBC^*_{g}VV\\
f_{Z_1}g^*\C{A}\otimes f_{Z_1}g^*\C{B} 
@>t_{Z_1}>>
f_{Z_1}g^*(\C{A}\otimes \C{B}).
\endCD
\]

{(V)} Each morphism of functors $d_Z:f_Z\B{D}\to \B{D} f_Z$ is an isomorphism.

(VI) For each $g\in Mor(c)$ such that $g$ and $\ov{g}$ are proper, the base change 
morphisms $BC_{*g}$ and $BC_{!g}$ are isomorphisms, and $BC_{!g}$ is inverse to $BC_{*g}$.

(VII) Cohomological pre-morphism $(f_Z,t_Z;BC^*_{g},BC^!_{g};\iota)$ can be extended 
to a cohomological pre-morphism between compactifications of $c$ and $\ov{c}$ satisfying axioms 
(I)-(VI) above.

(b)  A cohomological pre-morphism $(f_Z,t_Z;BC^*_{g},BC^!_{g};\iota)$ from  $c$ to $\ov{c}$
is called a {\em cohomological morphism}, if it decomposes as a composition of finitely many 
compactifiable cohomological morphisms.
\end{Def}

\begin{Emp} \label{E:becohmor}
{\bf Basic examples.}

(a) {\bf Proper push-forward.} Let  $[f]=(f,f^{\diez},f)$ be a proper morphism from a correspondence 
$c:C\to X\times X$ to $\ov{c}:\ov{C}\to \ov{X}\times\ov{X}$. Then $[f]$ gives rise to 
a compactifiable cohomological morphism defined as follows.
 
For each $Z\in Ob(c)$, $[f]$ defines a proper morphism $[f]_Z:Z\to\ov{Z}$ such that $[f]_{\pt}=\Id$
and $\ov{g}\circ[f]_{Z_1}=[f]_{Z_2}\circ g$ for all $g:Z_1\to Z_2\in Mor(c)$. Hence $[f]$ gives 
rise to a cohomological pre-morphism $(f_Z,t_Z;BC^*_{g},BC^!_{g};\iota)$ from $c$ to $\ov{c}$, 
where $f_Z$ is $[f]_{Z*}=[f]_{Z!}$, $t_Z$ is the morphism 
$t_{[f]_Z}$ of \form{tf}, $\iota$ is the identity, $BC^*_g$ and $BC^!_g$ are the usual 
base change morphisms $\ov{g}^*[f]_{Z_2*}\to [f]_{Z_1*}g^*$ and   
$[f]_{Z_1!}g^!\to \ov{g}^![f]_{Z_2!}$, while  $BC_{*g}$ and $BC_{!g}$ are isomorphisms
$ [f]_{Z_2*}g_*\isom \ov{g}_*[f]_{Z_1*}$ and 
$\ov{g}_![f]_{Z_1!}\isom [f]_{Z_2!}g_!$.

We claim that the pre-morphism defined above is actually a compactifiable cohomological morphism. 
Indeed, axioms (I)-(IV) easily follow by adjointness, while axiom (VI) is clear. Next 
axiom (V) follows from fact that $d_Z$ coincides with the isomorphism
$[f]_{Z*}\C{RHom}(\C{A},[f]^!_{Z} K_{\ov{Z}})\isom\C{RHom}([f]_{Z!}\C{A},K_{\ov{Z}})$, 
obtained by the sheafification of the adjointness $\Hom(\C{A},[f]^!_{Z} K_{\ov{Z}})\isom
\Hom([f]_{Z!}\C{A},K_{\ov{Z}})$. 

Finally, to prove axiom (VII) it remains to show that $[f]$ extends to a morphism 
$[f^{\star}]$ between compactifications of $c$  and $\ov{c}$. To construct $[f^{\star}]$, 
we argue as in \rr{comp}. Namely, first choose a compactification 
$\ov{c}^{\star}:\ov{C}^{\star}\to\ov{X}^{\star}\times\ov{X}^{\star}$ of $\ov{c}$, 
next choose a compactification $X^{\star}$ of $X$ such that $f$ extends to a 
morphism $f^{\star}:X^{\star}\to \ov{X}^{\star}$, and finally choose 
a  compactification $C^{\star}$ of $C$ such that 
$(f^{\diez},c):C\to \ov{C}\times (X\times X)$ extends to a morphism
$(f^{\diez,\star},c^{\star}):C^{\star}\to \ov{C}^{\star}\times (X^{\star}\times X^{\star})$.
Then $[f^{\star}]=(f^{\star}f^{\diez,\star},f^{\star})$ is a required morphism between 
compactifications $c^{\star}$ and $\ov{c}^{\star}$.

(b) {\bf Specialization.} Let $c:C\to X\times X$ be a correspondence over $k$, 
and let $\wt{c}:\wt{C}\to\wt{X}\times\wt{X}$ be a correspondence over $\C{D}$ 
lifting $c$. Then $\wt{c}$ gives rise to a cohomological morphism
from $c$ to $\ov{c}:=\wt{c}_s$ defined as follows. 

Replacing $\wt{c}$ by its pullback to 
$\C{D}^h$, we may assume that $R$ is strictly henselian.
For each $Z\in Ob(c)$, $\wt{c}$ defines a lift $\wt{Z}$ of $Z$, while for  
each $g:Z_1\to Z_2\in Mor(c)$, $\wt{c}$ defines a lift 
$\wt{g}:\wt{Z}_1\to\wt{Z}_2$ of $g$, whose special fiber $\wt{g}_s$ is the corresponding morphism 
$\ov{g}:\ov{Z}_1\to \ov{Z}_2\in Mor(\ov{c})$.   Hence $\wt{c}$ gives rise to a 
cohomological pre-morphism $(f_Z,t_Z;BC^*_{g},BC^!_{g};\iota)$ from $c$ to $\ov{c}$, 
where $f_Z$ is $sp_{\wt{Z}}$, (thus  $f_{\pt}$ is $sp_{\C{D}}=\Id$),   
$\iota$ is the identity, $BC_g^!, BC_g^*, BC_{!g}$ and $BC_{*g}$ are the maps of 
\re{specfunc} (b), and finally 
$t_Z:sp_{\wt{Z}}\C{A}\otimes sp_{\wt{Z}}\C{B}\to sp_{\wt{Z}}(\C{A}\otimes\C{B})$ is 
the morphism 
\[
\ov{i}^*(\ov{j}_*\varphi^*\C{A}_{\ov{\eta}}\otimes\ov{j}_*\varphi^*\C{B}_{\ov{\eta}})
\overset{t_{\ov{j}}}{\lra} \ov{i}^*\ov{j}_*\varphi^*(\C{A}\otimes\C{B})_{\ov{\eta}}.
\]

We claim that this pre-morphism satisfies axioms (I)-(VI) of \rd{cohmor}.
Indeed, axioms (I)-(IV) follow by adjointness, axiom (VI) 
follows from the proper base change theorem (implicitly used in the definition of $BC_{!g}$), 
while axiom (V) follows from the corresponding 
assertion for the functor $\Psi_{\wt{Z}}$, proven in \cite[Thm. 4.2]{Il3}. 

Assume first that $\wt{c}$ extends to a correspondence proper over $\C{D}$, 
lifting a compactification of $c$. In this case, axiom (VII) clearly holds, thus our cohomological 
pre-morphism is actually a compactifiable cohomological morphism. 

In general, it follows from \rco{comp} below, that there exists a correspondence $\wt{c}':\wt{C}'\to\wt{X}'\times\wt{X}'$ over $\C{D}$ lifting $c$ and
a proper morphism $\wt{[f]}$ from $\wt{c}'$ to $\wt{c}$ such that $\wt{c}'$ extends 
to a correspondence proper over $\C{D}$, lifting a compactification of $c$.
Therefore our cohomological pre-morphism decomposes as a composition of two 
compactifiable cohomological morphisms: the one from $c$ to $\wt{c}'_s$, defined by $\wt{c}'$, 
and the one from $\wt{c}'_s$ to $\wt{c}_s$, defined by $\wt{[f]}_s$. 
This proves the assertion.
\end{Emp} 

\begin{Lem} \label{L:comp1}
Given a scheme $X$ over $k$, a lift $\wt{X}$ of $X$ over $\C{D}$, and 
a compactification $\ov{X}$ of $X$, there exists a lift $\wt{\ov{X}}$ of $\ov{X}$, 
proper over $\C{D}$ and an open 
subscheme $\wt{\ov{X}}^0\subset\wt{\ov{X}}$ lifting $X$ and proper over $\wt{X}$.
\end{Lem}

\begin{proof}
Choose any compactification $j:\wt{X}\hra\wt{X}''$ over $\C{D}$, define $\wt{X}'$ be the closure of 
the image of $(j,\varphi):\wt{X}\hra\wt{X}''\times X_{\C{D}}$, and define 
$\varphi':\wt{X}'\to  X_{\C{D}}$ be the 
projection. By the construction, $\wt{X}'$ is a lift of $X$, proper over $X_{\C{D}}$, and $\wt{X}\subset\wt{X}'$ is an open 
subscheme. 

It follows from the explicit form of Chow's lemma 
(see, for example, \cite[Thm 2.1]{Co}), that there exists a closed subscheme $Z\subset X_{\C{D}}$,
 supported set-theoretically at $X$, such that the blow-up 
$B_Z(X_{\C{D}})\to X_{\C{D}}$ factors through $\wt{X}'\to X_{\C{D}}$. If we denote by 
 $\ov{Z}\subset \ov{X}$ the schematic closure of $Z$, then the assertion of lemma holds with 
$\wt{\ov{X}}$ equal the blow-up $Bl_{\ov{Z}}(\ov{X}_{\C{D}})$ and $\wt{\ov{X}}^0$ equal the 
preimage of $\wt{X}$ in $Bl_Z(X_{\C{D}})$. 
\end{proof}

\begin{Cor} \label{C:comp}
In the notation of \re{becohmor} (b), there exists a correspondence 
$\wt{c}':\wt{C}'\to\wt{X}'\times\wt{X}'$ over $\C{D}$ such that $\wt{c}'$ 
lifts $c$,  $\wt{c}'$ is proper over $\wt{c}$ and for every compactification 
$\ov{c}:\ov{C}\to\ov{X}\times\ov{X}$ of $c$, 
$\wt{c}'$ extends to a correspondence $\wt{\ov{c}}$ which is proper over $\C{D}$ and 
lifts $\ov{c}$.
\end{Cor}
\begin{proof}
Let $\wt{\ov{X}}^0\subset \wt{\ov{X}}$ and $\wt{\ov{C}}^0\subset \wt{\ov{C}}$ 
correspond by the lemma to triples $(X,\wt{X},\ov{X})$ and $(C,\wt{C},\ov{C})$, respectively.
Denote by $\wt{\ov{C}}'$ and  $\wt{\ov{C}}'^0 $ the closures of the embeddings 
$\ov{C}_{\eta}\hra\wt{\ov{C}}\times\wt{\ov{X}}\times\wt{\ov{X}}$ and 
$C_{\eta}\hra\wt{\ov{C}}^0\times\wt{\ov{X}}^0\times\wt{\ov{X}}^0$, induced by $\ov{c}$ and $c$, 
and let $\wt{\ov{c}}:\wt{\ov{C}}'\to \wt{\ov{X}}\times\wt{\ov{X}}$ and 
$\wt{\ov{c}}^0:\wt{\ov{C}}'^0\to \wt{\ov{X}}^0\times\wt{\ov{X}}^0$ be the projection maps. 
Then the correspondences $\wt{\ov{c}}^0$ and $\wt{\ov{c}}$ satisfy all the required properties.
\end{proof}

\begin{Rem}
If $\wt{c}$ is the correspondence $\wt{c}_Z$ from \re{defcor}, then it extends to a  
correspondence proper over $\C{D}$ and lifting $\ov{c}$. Therefore in this case 
the cohomological morphism, defined in \re{becohmor} (b) is compactifiable.
However, we do not know whether this is true in general. 
\end{Rem}
\subsection{Properties of compactifiable cohomological morphisms} \label{SS:cohmor2}

\begin*
\vskip 8truept
\end*

Fix a compactifiable cohomological morphism $(f_Z,t_Z;BC^*_g,BC^!_g;\iota)$
from $c$ to $\ov{c}$. 

\begin{Lem} \label{L:basechange}

(a) For each $Z\in\Ob(c)$, we have $BC^!_{\Id_Z}=\Id_{f_Z}$.
For for each $g_1:Z_1\to Z_2$ and $g_2:Z_3\to Z_3$ in $\Mor(c)$, we have
$BC^!_{g_2\circ g_1}= BC^!_{g_2}\circ BC^!_{g_1}$. 

(b) For each $g:Z_1\to Z_2\in\Mor(c)$, the morphism $\wt{\pi}_{Z_1}:f_{Z_1} K_{Z_1}\to K_{\ov{Z}_1}$
decomposes as $f_{Z_1}g^!K_{Z_2}\overset{BC^!_{g}}{\lra}\ov{g}^!f_{Z_2} K_{Z_2}
\overset{\wt{\pi}_{Z_2}}{\lra}\ov{g}^!K_{\ov{Z}_2}$.

(c) For each commutative diagram 
\[
\CD
Z_1 @>g_1>> Z_2\\
@Vg_3VV @Vg_4VV\\
Z_3 @>g_2>> Z_4
\endCD
\]
in $\Mor(c)$, the following diagrams are commutative
\[(c1)\;\;\;\;\;\;\;\;\;\;\;\;\;\;\;\;\;\;\;\;\;\;\;\;
\CD
f_{Z_3}g^*_2 g_{4*} @<BC^*_{g_2}<< \ov{g}^*_{2}f_{Z_4}g_{4*}  @>BC_{*g_4}>> 
 \ov{g}^*_{2}\ov{g}_{4*}f_{Z_2}\\
@VBCVV @. @VBCVV\\
f_{Z_3}g_{3*}g_1^* @>BC_{* g_3}>> \ov{g}_{3*}f_{Z_1}g_1^*  @<BC^*_{g_1}<< 
\ov{g}_{3*}\ov{g}_1^*f_{Z_2}
\endCD\;\;\;\;\;\;\;\;\;\;\;\;\;\;\;\;\;\;\;\;\;\;\;\;
\]

\[(c2)\;\;\;\;\;\;\;\;\;\;\;\;\;\;\;\;\;\;\;\;\;\;\;\;
\CD
f_{Z_3}g_{3!}g_1^! @<BC_{!g_3}<< \ov{g}_{3!}f_{Z_1}g_1^! @>BC^!_{g_1}>> 
 \ov{g}_{3!}\ov{g}_1^!f_{Z_2}\\
@VBCVV @. @VBCVV\\
f_{Z_3}g^!_2 g_{4!} @>BC^!_{g_2}>>\ov{g}^!_2 f_{Z_4}g_{4!} @<BC_{!g_4}<< 
\ov{g}^!_2 \ov{g}_{4!}f_{Z_2}
\endCD\;\;\;\;\;\;\;\;\;\;\;\;\;\;\;\;\;\;\;\;\;\;\;\;
\]
\end{Lem}

\begin{proof}

(a) The assertion is local, so we can replace $c$ and $\ov{c}$ by their compactifications, 
thus assuming that $g_1$, $g_2$, $\ov{g}_1$ and $\ov{g}_2$ are proper. 
Now the assertion follows from axioms (III) and (VI). Namely, axiom (III) implies 
the corresponding assertion for $BC_*$ (by adjointness), hence 
for  $BC_!$ (by axiom (VI)), and finally for $BC^!$ (by adjointness).  

(b) Since $\wt{\pi}_{Z_i}$ is the composition of $BC^!_{\pi_{Z_i}}$ and $\iota$, the assertion 
follows from the equality $BC^!_{\pi_{Z_1}}= BC^!_{\pi_{Z_2}}\circ BC^!_{g}$, proven in (a).

(c) we will prove the commutativity of $(c1)$, while that of $(c2)$ is similar.

Put $g:=g_2\circ g_3=g_4\circ g_1$. Then by axiom (III), we have
$BC^*_{g_3}\circ BC^*_{g_2}=BC^*_{g}=BC^*_{g_1}\circ BC^*_{g_4}$. 
Since the composition $\ov{g_3}^* f_{Z_3}g_{3*}\overset{BC^*_{g_3}}{\lra} 
f_{Z_1} g_3^* g_{3*}\overset{adj}{\lra}f_{Z_1}$ equals
$\ov{g_3}^* f_{Z_3}g_{3*}\overset{BC_{*g_3}}{\lra}\ov{g_3}^*\ov{g}_{3*} f_{Z_1}
\overset{adj}{\lra}f_{Z_1}$,
while the composition $\ov{g_4}^* f_{Z_4}g_{4*}\overset{BC^*_{g_4}}{\lra} 
f_{Z_2} g_4^* g_{4*}\overset{adj}{\lra}f_{Z_2}$ equals
$\ov{g_4}^* f_{Z_4}g_{4*}\overset{BC_{*g_4}}{\lra}\ov{g_4}^*\ov{g}_{4*} f_{Z_2}
\overset{adj}{\lra}f_{Z_2}$, 
both morphisms $\ov{g}^*_{2}f_{Z_4}g_{4*}\to \ov{g}_{3*}f_{Z_1}g_1^*$ are
adjoint to the composition 
\[
\ov{g}^*f_{Z_4}g_{4*}\overset{BC^*_{g}}{\lra} 
f_{Z_1} g^* g_{4*}=f_{Z_1} g_1^*g^*_4 g_{4*}\overset{adj}{\lra}f_{Z_1}{g_1}^*. 
\]
Hence they are equal. 
\end{proof}

\begin{Lem} \label{L:proper}
For each $g:Z_1\to Z_2\in\Mor(c)$, the following diagrams are commutative
\begin{equation*}(a)\;\;\;\;\;\;\;\;\;\;\;\;\;\;\;
\CD
f_{Z_1} g^!\C{A}\otimes \ov{g}^*f_{Z_2}\C{B}  @>BC^*_{g}>>
f_{Z_1} g^!\C{A}\otimes f_{Z_1} g^*\C{B}  @>t_{g^!}\circ t_{Z_1}>>
f_{Z_1} g^!(\C{A}\otimes \C{B})\\
@VBC^!_{g}VV @. @VBC^!_{g}VV\\
\ov{g}^! f_{Z_2}\C{A}\otimes \ov{g}^*f_{Z_2}\C{B}  @>t_{\ov{g}^!}>>
\ov{g}^!(f_{Z_2}\C{A}\otimes f_{Z_2}\C{B})  @> t_{Z_2}>>
\ov{g}^!f_{Z_2}(\C{A}\otimes\C{B});
\endCD\;\;\;\;\;\;\;\;\;\;\;\;\;\;\;
\end{equation*}
\[(b)
\CD
\!\!\!\!\!\!f_{Z_1} g^!\C{RHom}(\C{A},\!\C{B}) @>BC^!_g>> 
\!\!\!\!\!\! \ov{g}^!f_{Z_2}\C{RHom}(\C{A},\C{B}) @>r_{Z_2}>>
\!\! \!\!\!\!\ov{g}^!\C{RHom}(f_{Z_2}\C{A},f_{Z_2}\C{B})\\
@Vt_{g^!}VV @. @Vt_{\ov{g}^!}VV\\
\!\!\!\!\!\! f_{Z_1}\C{RHom}(g^*\C{A},\!g^!\C{B})\!\!\! @>r_{Z_1}>> 
\!\!\!\!\C{RHom}(f_{Z_1}g^*\C{A},\!f_{Z_1}g^!\C{B})\!\!\!@>(BC^*_g, BC^!_g)>>
\!\!\!\!\C{RHom}(\ov{g}^*f_{Z_2}\C{A},\!\ov{g}^!f_{Z_2}\C{B}).
\endCD
\]
\end{Lem}

\begin{proof}
(a) The assertion is local, so we can replace $c$ and $\ov{c}$ by their compactifications. 
Thus we can 
assume that $g$ and $\ov{g}$ are proper.
Now the assertion follows from  axioms (IV) and (VI) by adjointness. 
Namely, the diagram adjoint to (a) extends to the diagram 
\begin{equation} \label{Eq:a}
\CD
\ov{g}_!(f_{Z_1} g^!\C{A}\otimes \ov{g}^*f_{Z_2}\C{B})  @> t_{Z_1}\circ BC^*_{g}>>
\ov{g}_! f_{Z_1}(g^!\C{A}\otimes g^*\C{B})  @>t_{g^!}>>
\ov{g}_! f_{Z_1} g^!(\C{A}\otimes \C{B})\\
@VprojVV @VBC_{!g}VV @VBC_{!g}VV\\
\ov{g}_!f_{Z_1} g^!\C{A}\otimes f_{Z_2}\C{B}  @. 
f_{Z_2}g_!(g^!\C{A}\otimes g^*\C{B})  @>t_{g^!}>>
f_{Z_2} g_!g^!(\C{A}\otimes \C{B})\\
 @VBC_{!g}VV @VprojVV @VadjVV\\
f_{Z_2} g_!g^!\C{A}\otimes f_{Z_2}\C{B}  @>t_{Z_2}>> 
f_{Z_2}(g_!g^!\C{A}\otimes \C{B}) @>adj>> f_{Z_2}(\C{A}\otimes\C{B}).
\endCD
\end{equation} 
The top right inner square of \form{a} is commutative by functoriality, while 
the bottom right inner square is commutative by the definition of $t_{g^!}$.

For the left inner square  of \form{a}, recall that $g$ and $\ov{g}$ are assumed 
to be proper, therefore
by axiom (VI) all the vertical morphisms are isomorphisms. 
Therefore we can replace $g_!$ by $g_*$,  $\ov{g}_!$ by $\ov{g}_*$ 
and all vertical morphisms by their inverses.
Then by adjointness, we get a diagram 
\begin{equation} \label{Eq:a'}
\CD
f_{Z_1} g^!\C{A}\otimes f_{Z_1}g^*\C{B}  @> t_{Z_1}>>
f_{Z_1}(g^!\C{A}\otimes g^*\C{B})\\
@AadjAA @AadjAA\\
f_{Z_1}g^* g_* g^!\C{A}\otimes f_{Z_1} g^*\C{B} @>t_{Z_1}>> 
f_{Z_1}(g^* g_* g^!\C{A}\otimes g^* \C{B})\\
 @ABC^*_{g}\otimes BC^*_{g}AA  @ABC^*_{g}AA\\
\ov{g}^* f_{Z_2} g_* g^!\C{A}\otimes \ov{g}^* f_{Z_2}\C{B}  @>t_{Z_2}>> 
\ov{g}^* f_{Z_2}(g_* g^!\C{A}\otimes \C{B}).
\endCD
\end{equation}
As the bottom inner square of \form{a'} is commutative by axiom (IV), and the top one 
is commutative by functoriality, we obtain that diagram \form{a'} and hence
the left inner square of \form{a} is commutative. This completes the proof.

(b) The assertion follows from (a) by adjointness. Namely,
consider diagram 
\begin{equation} \label{Eq:b}
\CD
f_{Z_1} g^!\C{RHom}(\C{A},\C{B})\otimes \ov{g}^*f_{Z_2}\C{A}
\!@>>>\! f_{Z_1}g^!(\C{RHom}(\C{A},\C{B})\otimes\C{A})\!@>ev>>\!f_{Z_1}g^!\C{B}\\
@VBC^!_gVV       @VBC^!_gVV @VBC^!_gVV\\
\ov{g}^!f_{Z_2}\C{RHom}(\C{A},\C{B})\otimes \ov{g}^*f_{Z_2}\C{A}
\!@>>>\! \ov{g}^!f_{Z_2}(\C{RHom}(\C{A},\C{B})\otimes\C{A})\!@>ev>>\!\ov{g}^!f_{Z_2}\C{B},
\endCD
\end{equation}
where the left inner square is the commutative diagram of (a). Since the right inner square 
of \form{b} is commutative by functoriality, the exterior square of \form{b} is commutative. 
Since the composition maps 
\[
f_{Z_1} g^!\C{RHom}(\C{A},\C{B})\otimes 
\ov{g}^*f_{Z_2}\C{A}\to \ov{g}^!f_{Z_2}\C{B},
\]
of the exterior square of \form{b}, are adjoint to the composition maps 
\[
f_{Z_1} g^!\C{RHom}(\C{A},\C{B})\to \C{RHom}(\ov{g}^*f_{Z_2}\C{A},\ov{g}^!f_{Z_2}\C{B})
\]
of the diagram (b), we obtain the assertion.
\end{proof}

\begin{Not} \label{N:auxmor}

(i)  For $Z\in\Ob(c)$, denote by $i_Z:\La_{\ov{Z}}\to f_Z\La_{{Z}}$
the composition 
\[
\pi_{\ov{Z}}^* \La_{\ov{\pt}}\overset{\iota}{\lra}\pi_{\ov{Z}}^*f_{\pt}\La_{\pt}
\overset{BC^*_{\pi_Z}}{\lra}f_Z(\pi_Z^*\La_{{Z}})
\]
and by $h^0_Z:H^0(Z,\C{A})\to H^0(\ov{Z},f_Z\C{A})$ the composition 
\[
\Hom(\La_Z,\C{A})\overset{f_Z}{\lra} \Hom(f_Z\La_Z,f_Z\C{A})\overset{i_Z}{\lra}
\Hom(\La_{\ov{Z}},f_Z\C{A}).
\]

(ii) For $\C{A},\C{B}\in  D_{ctf}^b(X,\La)$, denote by 
$\Pi_{\C{A},\C{B}}:f_X\C{A}\pp f_X\C{B}\to f_{X\times X}(\C{A}\pp \C{B})$ the 
composition
\[
\ov{p}_1^*f_X\C{A}\otimes\ov{p}_2^* f_X\C{B}\overset{BC^*_{p_1}\otimes BC^*_{p_2}}{\lra}
f_{X\times X}p_1^*\C{A}\otimes f_{X\times X}p_2^*\C{B}\overset{t_{X\times X}}{\lra}
f_{X\times X} (p_1^*\C{A}\otimes p_2^*\C{B}).
\]
\end{Not}

\begin{Lem} \label{L:iZ}
(a) For each $g:Z_1\to Z_2\in\Mor(c)$, the morphism $i_{Z_1}:\La_{\ov{Z}_1}\to f_{Z_1}\La_{Z_1}$
decomposes as 
$\ov{g}^*\La_{\ov{Z}_2}\overset{i_{Z_2}}{\lra}\ov{g}^* f_{Z_2}\La_{Z_2}
\overset{BC^*_{g}}{\lra}f_{Z_1}g^*\La_{Z_2}$.

(b) The identity morphism $\Id_{f_Z\C{A}}$ decomposes as 
\[
t_Z\circ i_Z: f_Z\C{A}=f_Z\C{A}\otimes\La_{\ov{Z}}\overset{i_Z}{\lra}
f_Z\C{A}\otimes f_Z\La_{{Z}}\overset{t_{Z}}{\lra}
f_Z(\C{A}\otimes\La_{{Z}})=f_Z\C{A}.
\]

(c) Each morphism $f_Z:\Hom(\C{A},\C{B})\to\Hom(f_Z\C{A},f_Z\C{B})$ decomposes as 
\[
H^0(Z,\C{RHom}(\C{A},\C{B}))\overset{h^0_Z}{\lra}
H^0(\ov{Z},f_Z(\C{RHom}(\C{A},\C{B})))\overset{r_Z}{\lra}
H^0(\ov{Z},\C{RHom}(f_Z\C{A},f_Z\C{B})).
\]

(d) For each $g:Z_1\to Z_2\in\Mor(c)$, the composition
\[
H^0(Z_1,\C{A})=H^0(Z_2,g_*\C{A})\overset{h^0_{Z_2}}{\lra}H^0(\ov{Z_2},f_{Z_2}g_*\C{A})
\overset{BC_{*g}}{\lra}H^0(\ov{Z_2},\ov{g}_*f_{Z_1}\C{A})=H^0(\ov{Z_1},f_{Z_1}\C{A})
\]
equals $h^0_{Z_1}:H^0(Z_1,\C{A})\to H^0(\ov{Z_1},f_{Z_1}\C{A})$.

(e) Morphism $t_X:f_X\C{A}\otimes f_X\C{B}\to f_{X}(\C{A}\otimes\C{B})$ decomposes as
\[
\ov{\Dt}^*(f_X\C{A}\pp f_X\C{B})\overset{\Pi_{\C{A},\C{B}}}{\lra}\ov{\Dt}^* 
f_{X\times X}(\C{A}\pp \C{B})\overset{BC^*_{\Dt}}{\lra} f_{X}\Dt^*(\C{A}\pp \C{B}).
\]

(f) The diagram
\[
\CD
\ov{p}^*_1 f_X K_X @>\wt{\pi}_X>> \ov{p}^*_1 K_{\ov{X}} @>BC>>  \ov{p}_2^!\La_{\ov{X}}\\
@VBC_{p_1}^*VV @. @Vi_XVV\\
 f_{X\times X}p^*_1 K_X @>BC>>  f_{X\times X}p^!_2\La_X @>BC_{p_2}^!>>
\ov{p}_2^!f_X\La_{X},
\endCD
\]
where we denote by $BC$ base change morphisms 
$\ov{p}^*_1\pi_{\ov{X}}^!\La_{\ov{\pt}}\to\ov{p}_2^!\pi_{\ov{X}}^*\La_{\ov{\pt}}$ and
 $p^*_1\pi_{X}^!\La_{\pt}\to p_2^!\pi_{X}^*\La_{\pt}$, is commutative. 
\end{Lem}
\begin{proof}
(a)  Since $i_{Z_i}$ is the composition of  $\iota$ and $BC^*_{\pi_{Z_i}}$, the assertion 
follow from the equality $BC^*_{\pi_{Z_1}}=BC^*_{g}\circ BC^*_{\pi_{Z_2}}$
(axiom (III)).

(b) Since the composition 
\begin{equation} \label{Eq:can}
f_Z\C{A}\otimes  f_Z\B{D}\C{A}\overset{t_Z}{\lra}f_Z(\C{A}\otimes\B{D}\C{A})
\overset{ev}{\lra}f_Z K_Z\overset{\wt{\pi}_Z}{\lra} K_{\ov{Z}},
\end{equation}
induces an isomorphism  $f_Z\B{D}\C{A}\isom \B{D}f_Z\C{A}$ (axiom (V)), it will suffice 
to check that the composition 
\begin{equation} \label{Eq:can'}
f_Z\C{A}\otimes  f_Z\B{D}\C{A}\overset{(t_Z\circ i_Z)\otimes\Id}{\lra}
f_Z\C{A}\otimes f_Z\B{D}\C{A}\overset{\form{can}}{\lra} K_{\ov{Z}}
\end{equation}
equals \form{can}. 
By associativity and commutativity of $t_Z$
(axiom (I)), the map \form{can'} decomposes as a composition 
\[
f_Z\C{A}\otimes  f_Z\B{D}\C{A}\overset{ev\circ t_Z}{\lra}
f_Z K_Z\overset{i_Z}{\lra}f_Z K_Z\otimes f_Z\La_Z
\overset{t_Z}{\lra} f_Z K_Z
\overset{\wt{\pi}_Z}{\lra} K_{\ov{Z}}.
\]
Thus we have to show that the composition
\[
f_Z K_Z\overset{i_Z}{\lra}f_Z K_Z\otimes f_Z\La_Z
\overset{t_Z}{\lra} f_Z K_Z\overset{\wt{\pi}_Z}{\lra}  K_{\ov{Z}}
\]
equals $\wt{\pi}_Z$. Using axiom (II), the assertion reduces to the commutativity of diagram (a) of 
\rl{proper} in the case $g=\pi_Z$ and $\C{A}=\C{B}=\La_{\ov{Z}}$.

(c) For each $y\in\Hom(\C{A},\C{B})$, the element 
$r_Z\circ h^0_Z(y)\in\Hom(f_Z\C{A},f_Z\C{B})$ equals the composition
$f_Z\C{A}\overset{i_{Z}}{\lra}f_Z\C{A}\otimes f_Z\La_{{Z}}\overset{t_{Z}}{\lra}
f_Z\C{A}\overset{f_{Z}(y)}{\lra}f_Z\C{B}$.
Hence by (a) it equals $f_{Z}(y)$, as claimed.

(d) Since $BC_{*g}$ is defined by adjointness from $BC^*_g$, 
the composition sends morphism 
$(g^*\La_{Z_2}\overset{h}{\lra}\C{A})\in\Hom(g^*\La_{Z_2},\C{A})=H^0(Z_1,\C{A})$ to
the morphism 
\[
(\ov{g}^*\La_{\ov{Z}_2}\overset{i_{Z_2}}{\lra}\ov{g}^* f_{Z_2}\La_{Z_2}
\overset{BC^*_{g}}{\lra}f_{Z_1}g^*\La_{Z_2}\overset{h}{\lra} f_{Z_1}\C{A})\in
\Hom(\ov{g}^*\La_{\ov{Z}_2},f_{Z_1}\C{A})=H^0(\ov{Z}_1,f_{Z_1}\C{A}).
\]
Therefore the assertion follows from (a).

(e) Using identities $p_i\circ\Dt=\Id_X$ and  $\ov{p}_i\circ\ov{\Dt}=\Id_{\ov{X}}$, we can form  
a diagram
\[
\CD
\!\!\ov{\Dt}^*\ov{p}_1^*f_X\C{A}\otimes \ov{\Dt}^*\ov{p}_2^*f_X\C{B} 
\!\!@>BC^*_{p_1}\otimes BC^*_{p_2}>>\!\!
\ov{\Dt}^*f_{X\times X}{p}_1^*\C{A}\otimes \ov{\Dt}^*f_{X\times X}{p}_2^*\C{B}
\!\! @>BC^*_{\Dt}\otimes BC^*_{\Dt}>>\!\! f_{X}\C{A}\otimes  f_{X}\C{A}\\
@| @Vt_{X\times X}VV @Vt_XVV\\
f_X\C{A}\otimes f_X\C{B} @>\Pi_{\C{A},\C{B}}>> 
\ov{\Dt}^*f_{X\times X}({p}_1^*\C{A}\otimes {p}_2^*\C{B})
 @>BC^*_{\Dt}>>  f_{X}(\C{A}\otimes\C{B}).
\endCD
\] 
Its left inner square is commutative by the definition of $\Pi_{\C{A},\C{B}}$,
while the right inner square is commutative by axiom (IV). 
Since the composition of the top horisontal morphisms is the identity by axiom (III), 
the assertion follows.

(f) The assertion is local, so we can replace $c$ and $\ov{c}$ by their compactifications,
hence we can assume that all morphisms in question are proper. By adjointness, we have to show 
the commutativity of the diagram
\begin{equation} \label{Eq:c}
\CD
\ov{p}_{2!}\ov{p}^*_1 f_X K_X @>BC>> \ov{\pi}^*_X\ov{\pi}_{X!}f_X K_X @>(1)>> 
\ov{\pi}^*_X\La_{\ov{\pt}}=\La_{\ov{X}}\\
@VBC_{! p_2}\circ BC_{p_1}^*VV @VBC^*_{\pi_X}\circ BC_{!\pi_X}^*VV @Vi_XVV\\
 f_{X\times X}p_{2!}p^*_1 K_X @>BC>>  f_{X\times X}\pi_X^*\pi_{X!} K_X 
@>adj>> f_X\La_{X},
\endCD
\end{equation}
where $(1)$ is induced by the morphism $\ov{\pi}_{X!}f_X K_X\to \La_{\ov{\pt}}$
adjoint to $\wt{\pi}_X$. 
Since all morphisms are assumed to be proper,
the commutativity of the left inner square of \form{c} follows from \rl{basechange} (c1) 
and axiom (VI). 
Next  by \rl{basechange} (b) applied to $\pi_X$, the map 
$(1)$ decomposes as a composition 
\[
\ov{\pi}_{X!}f_X K_X\overset{BC_{!\pi_X}}{\lra}f_{\pt}\pi_{X!} K_X
\overset{adj}{\lra}f_{\pt}\La_{\pt}\overset{\iota}{\lra}\La_{\ov{\pt}}.
\]
Since $i_X$ decomposes as $\pi_{\ov{X}}^*\La_{\ov{\pt}}\overset{\iota^{-1}}{\lra}
\pi_X^* f_{\pt}\La_{\pt}\overset{BC_{\pi_X}^*}{\lra}f_X\pi_X^*\La_{\pt}$, 
the commutativity of the right inner square of \form{c} follows.
\end{proof}

\subsection{Compatibility with trace maps} \label{SS:comp}

\begin*
\vskip 8truept
\end*

Now we formulate the main result of this section, whose proof will be given in 
Subsection \ref{SS:pfofcd}.

\begin{Thm} \label{T:cohmor}
For every cohomological morphism 
$(f_Z,t_Z;BC^*_g,BC^!_{g};\iota)$ from $c$ to $\ov{c}$ and every $\C{F}\in D_{ctf}^b(X,\La)$,
the following diagram is commutative
\begin{equation} \label{Eq:cohmor}
\CD 
    f_C\C{RHom}(c_1^*\C{F},c_2^!\C{F}) @>{f_C(\un{\C{Tr}}_c)}>> 
    f_C\Dt'_* K_{Fix(c)}\\
      @V{(BC_{c_1}^*, BC_{c_2}^!)\circ r_C}VV  @V{\wt{\pi}_{Fix(c)}\circ BC_{*\Dt'}}VV\\
     \C{RHom}(\ov{c}_1^*f_X\C{F},\ov{c}_2^!f_X\C{F}) @>\un{\C{Tr}}_{\ov{c}}>> 
\ov{\Dt}'_* K_{Fix(\ov{c})}.
\endCD 
\end{equation}
\end{Thm}
\begin{Cor} \label{C:cohmor}
For every cohomological morphism 
$(f_Z,t_Z;BC^*_g, BC^!_{g};\iota)$ from $c$ to $\ov{c}$ and every $\C{F}\in D_{ctf}^b(X,\La)$, 
the following diagram is commutative
\begin{equation} \label{Eq:cohmor'}
\CD 
    \Hom(c_{2!}c_1^*\C{F},\C{F}) @>{\C{Tr}_c}>> H^0(Fix(c), K_{Fix(c)})\\
      @V{BC_{c_1}^* \circ BC_{!c_2}\circ f_X}VV     @V\wt{\pi}_{Fix(c)}\circ h^0_{Fix(c)}VV\\
     \Hom(\ov{c}_{2!}\ov{c}_1^*f_X\C{F},f_X\C{F}) @>\C{Tr}_{\ov{c}}>> 
H^0(Fix(\ov{c}),K_{Fix(\ov{c})}).
\endCD 
\end{equation}
\end{Cor}

\begin{Rem}
Since the assertions of the theorem and the corollary are compatible with compositions,
we can assume that cohomological morphism  $(f_Z,t_Z;BC^*_g, BC^!_{g};\iota)$ is 
compactifiable. Thus we can apply the results of Subsection \ref{SS:cohmor2}. 
\end{Rem}

\begin{proof}[Proof of the Corollary]
By \rt{cohmor}, the diagram
\[
\CD 
    H^0(C,\C{RHom}(c_1^*\C{F},c_2^!\C{F})) @>{\un{\C{Tr}}_c}>> 
     H^0(C,\Dt'_* K_{Fix(c)})\\
      @V{(BC_{c_1}^*, BC_{c_2}^!)\circ r_C\circ h^0_C}VV  
@V{\wt{\pi}_{Fix(c)}\circ BC_{*\Dt'}\circ h^0_C}VV\\
      H^0(\ov{C},\C{RHom}(\ov{c}_1^*f_X\C{F},\ov{c}_2^!f_X\C{F})) @>\un{\C{Tr}}_{\ov{c}}>> 
 H^0(\ov{C},\ov{\Dt}'_* K_{Fix(\ov{c})})
\endCD
\]
is commutative. Therefore it will suffice to show 
the commutativity of 
the following two diagrams: 
\begin{equation} \label{Eq:cohmor''}
\CD 
   \!\!  \Hom(c_{2!}c_1^*\C{F}_1,\C{F}_2) @>adj>> \Hom(c_1^*\C{F}_1,c^!_{2}\C{F}_2)@=
H^0(C, \C{RHom}(c_1^*\C{F}_1,c_2^!\C{F}_2))\\
      @V{f_X}VV  @VV{f_C}V 
@V{r_C \circ h^0_C}VV\\    
 \!\!\Hom(f_X c_{2!}c_1^*\C{F}_1,f_X \C{F}_2) @.\!\!\!\!
\Hom(f_C c_1^*\C{F}_1,f_C c^!_{2}\C{F}_2)\!\! @=  
\!\!\! H^0(\ov{C}, \!\C{RHom}(f_C c_1^*\C{F}_1,\!f_C c_2^!\C{F}_2))\\
      @V{BC_{c_1}^* \circ BC_{!c_2}}VV  @V{(BC_{c_1}^*, BC_{c_2}^!)}VV 
@V{(BC_{c_1}^*, BC_{c_2}^!)}VV \\
   \!\! \!\!\!\! \!\!\Hom(\ov{c}_{2!}\ov{c}_1^*f_X\C{F}_1,\!f_X\C{F}_2)\!\!\!\!\!\!\!\! \!\! 
@>adj>>\!\!\!  
\Hom(\ov{c}_1^*f_X\C{F}_1,\!\ov{c}^!_{2}f_X\C{F}_2)\!\! @= 
\!\! H^0(\ov{C},\!\C{RHom}(\ov{c}_1^*f_X\C{F}_1,\!\ov{c}_2^!f_X\C{F}_2))
\endCD 
\end{equation}

\begin{equation} \label{Eq:cohmor'''}
\CD 
 H^0(C,\Dt'_* K_{Fix(c)}) @= H^0(Fix(c),K_{Fix(c)})\\
 @VBC_{*\Dt}\circ h^0_CVV  @V{h^0_{Fix(c)}}VV\\  
H^0(\ov{C},\Dt'_* f_{Fix(c)} K_{Fix(c)}) @= H^0(Fix(\ov{c}),f_{Fix(c)} K_{Fix(c)})\\
 @V\wt{\pi}_{Fix(c)}VV  @V{\wt{\pi}_{Fix(c)}}VV\\ 
 H^0(\ov{C},\ov{\Dt}'_* K_{Fix(\ov{c})}) @= H^0(Fix(\ov{c}), K_{Fix(\ov{c})}).
\endCD 
\end{equation}
The top right inner square of \form{cohmor''} is commutative by \rl{iZ} (c).
The bottom right inner square of \form{cohmor''} is commutative by functoriality.
The commutativity of the left inner square of \form{cohmor''} follows from 
the adjointness of $BC_{!c_2}$ and $BC^!_{c_2}$.
Finally, the bottom inner square of \form{cohmor'''} is commutative by functoriality,
while the commutativity of the top inner square of \form{cohmor'''} follows from \rl{iZ} (d).
\end{proof}

\begin{Emp} \label{E:proofs}
{\bf Proof of Propositions \ref{P:pushf} and \ref{P:spec}.}

Both results follow immediately from \rco{cohmor} and 
\re{becohmor}. Namely, in the case of cohomological morphism of \re{becohmor} (a), 
diagram \form{cohmor'} specializes to diagram \form{pushf}, while in the case of 
cohomological morphism of \re{becohmor} (b), diagram \form{cohmor'} specializes to 
diagram \form{spec}.
\end{Emp}

\subsection{Proof of \rt{cohmor}} \label{SS:pfofcd}

\begin*
\vskip 8truept
\end*

\begin{Emp}
{\bf Particular case.}
First we will prove the assertion in the particular case 
$c=\Id_{X\times X}$ and $\ov{c}=\Id_{\ov{X}\times \ov{X}}$. 
By the definition of the trace maps, in this case 
diagram \form{cohmor} extends to diagram 
\begin{equation} \label{Eq:cohmor1}
\CD 
    f_{X\times X}\C{RHom}(p_1^*\C{F},p_2^!\C{F}) @<<<  
f_{X\times X}(\B{D}\C{F}\pp\C{F})  @>>> f_{X\times X}\Dt_* K_{X}\\
      @VV{(BC^*_{p_1},BC^!_{p_2})\circ r_{X\times X}}V  @AA{\Pi_{\B{D}\C{F},\C{F}}}A                         @V{\wt{\pi}_X\circ BC_{*\Dt}}VV\\
     \C{RHom}(\ov{p}_1^*f_X\C{F},\ov{p}_2^!f_X\C{F})  @<<<
 f_X\B{D}\C{F}\pp f_X\C{F}  @>>>\ov{\Dt}_* K_{\ov{X}}.
\endCD 
\end{equation}
We claim that both inner squares of \form{cohmor1} are commutative.
\end{Emp}

\begin{Emp}
{\bf Proof of the commutativity of the right inner square of \form{cohmor1}.}

Consider diagram
\begin{equation} \label{Eq:cohmor2}
\CD 
    \ov{\Dt}^* f_{X\times X}(\B{D}\C{F}\pp\C{F})  @>BC^*_{\Dt}>>  
 f_{X}\Dt^*(\B{D}\C{F}\pp\C{F})  @>ev>>
f_{X} K_{X}\\
   @AA{\Dt^*\Pi_{\B{D}\C{F},\C{F}}}A      @|    @|\\
f_X\B{D}\C{F}\otimes f_X\C{F}@>t_X>> f_X(\B{D}\C{F}\otimes \C{F})
 @>ev >>  f_X K_X\\
 @| & &  @VV{\wt{\pi}_X}V\\
f_X\B{D}\C{F}\otimes f_X\C{F}
@>d_X>>\B{D}f_X\C{F}\otimes f_X\C{F}@>ev>>K_{\ov{X}}.
\endCD 
\end{equation}
We claim that all three inner squares of \form{cohmor2} are commutative. Indeed, 
the bottom inner square  of \form{cohmor2} is commutative by the definition of $d_X$.
The top left inner square of \form{cohmor2} is commutative by \rl{iZ} (e),
while the commutativity of the top right inner square  of \form{cohmor2} is clear. 
Since the exterior square of \form{cohmor2} is adjoint to the right inner square
of \form{cohmor1}, the assertion follows.
\end{Emp}

\begin{Emp}
{\bf Proof of the commutativity of the left inner square of \form{cohmor1}.}

We have to check the equality of the two morphisms 
\[
f_X\B{D}\C{F}\pp f_X\C{F}\to\C{RHom}(\ov{p}_1^*f_X\C{F},\ov{p}_2^!f_X\C{F})
\]
from the  left inner square of \form{cohmor1}, or, equivalently, of the corresponding 
morphisms 
\[
\ov{p}_1^*f_X\C{F}\otimes\ov{p}_1^* f_X\B{D}\C{F}\otimes \ov{p}_2^*f_X\C{F}\to 
\ov{p}_2^!f_X\C{F}.
\]  
The assertion reduces to the commutativity of the following diagram

\begin{equation} \label{Eq:cohmor3}
\CD 
[\ov{p}_1^*f_X\C{F}\otimes\ov{p}_1^*f_X\B{D}\C{F}]\otimes \ov{p}_2^* f_X\C{F}
@>ev\circ t_X>> \ov{p}_1^* f_X K_X\otimes \ov{p}_2^*f_X\C{F} @>\diam\circ\wt{\pi}_X>> 
\ov{p}_2^!f_X\C{F}\\
@VV{BC_{p_1}^*\otimes BC_{p_1}^*\otimes BC_{p_2}^*}V  @V{BC_{p_1}^*\otimes BC_{p_2}^*}VV  
@A{BC_{p_2}^!}AA\\
\!\!\!\![f_{X\times X}p_1^*\C{F}\!\otimes\!  f_{X\times X}p_1^*\B{D}\C{F}]
\!\otimes\! f_{X\times X}p_2^* \C{F}\!\!\!@>ev\circ t_{X\times X}>>\!\!\!
f_{X\times X}p_1^* K_X\!\otimes\!
f_{X\times X}p_2^* \C{F}\!\!\! @>\diam\circ t_{X\times X}>>\!\!\! f_{X\times X}p_2^!\C{F},
\endCD 
\end{equation}
where $\diam$ is the isomorphism, defined in \form{*}.  
The left inner square of \form{cohmor3} commutes by axiom (IV) for the morphism $p_1$. 

Using \rl{iZ} (b), the right inner square of \form{cohmor3} extends to the diagram
\begin{equation} \label{Eq:cohmor4}
\CD 
\ov{p}_1^* f_X K_X\otimes \ov{p}_2^*f_X\C{F}\!\!\! @>(1)>>\!\!\! 
\ov{p}_2^! f_X\La_X\!\otimes\!\ov{p}_2^*f_X\C{F}\!\!\! @>t_X\circ t_{\ov{p}_2^!}>> \!\!\!\ov{p}_2^!f_X\C{F}\\
@VV{BC_{p_1}^*}V  @A{BC_{p_2}^!}AA @|\\
f_{X\times X}p_1^* K_X\!\otimes\!\ov{p}_2^*f_X\C{F}\!\!\! 
@>BC>>\!\!\! f_{X\times X}p_2^!\La_X\!\otimes\!\ov{p}_2^*f_X\C{F}@. \!\!\!\ov{p}_2^!f_X\C{F}\\ 
@VV{BC_{p_2}^*}V  @V{BC_{p_2}^*}VV  @A{BC_{p_2}^!}AA\\
f_{X\times X}p_1^* K_X\!\otimes\! f_{X\times X}p_2^* \C{F}\!\!\! @>BC>>\!\!\!  
f_{X\times X}p_2^!\La_X\!\otimes\! f_{X\times X}p_2^* \C{F} 
\!\!\!@>t_{p_2^!}\circ t_{X\times X}>>\!\!\! f_{X\times X}p_2^!\C{F},
\endCD 
\end{equation}
where $(1)$ is induced by the top morphism $\ov{p}_1^* f_X K_X\to \ov{p}_2^! f_X\La_X$
of \rl{iZ} (f). The  top left inner square of \form{cohmor4} is commutative by 
\rl{iZ} (f). The  bottom left inner square of \form{cohmor4} is commutative
by functoriality. Finally, the right 
inner square of \form{cohmor4} is commutative by \rl{proper} (a) for the morphism $p_2$.
\end{Emp}

This completes the proof of \rt{cohmor} in the case $c=\Id_{X\times X}$
and   $\ov{c}=\Id_{\ov{X}\times \ov{X}}$.
\begin{Emp}
{\bf The general case.}
Diagram \form{cohmor} extends to the diagram 
\begin{equation} \label{Eq:cohmor7}
\CD 
   \!\!\! \!\!\! f_C\C{RHom}(c_1^*\C{F}\!,\!c_2^!\C{F})\!\!\!\!\!\!\! @<t_{c^!}<<\!\!\!\!\!\!\!\!\!\!\!\!\!\!
f_C c^!\C{RHom}(p_1^*\C{F}\!,\!p_2^!\C{F})\!\!\!\!\!\!\!\!\!\!\!@>(3)>>\!\!\!
    f_C c^! \Dt_* K_{X} \!\!\! @<BC<<\!\!\! f_C\Dt'_* K_{Fix(c)}\\
      @VV{(BC_{c_1}^*, BC_{c_2}^!)\circ r_C}V  @V(1)\circ BC^!_c VV @V(2)\circ BC^!_cVV 
@V{\wt{\pi}_{Fix(c)}\circ BC_{*\Dt'}}VV\\
    \!\!\! \C{RHom}(\ov{c}_1^*f_X\C{F}\!,\!\ov{c}_2^!f_X\C{F})\!\!\! 
 @<t_{\ov{c}^!}<<\!\!\!\ov{c}^!\C{RHom}(\ov{p}_1^*f_X\C{F}\!,\!\ov{p}_2^!f_X\C{F})
\!\!\!@>(4)>>\!\!\!\!\!
\ov{c}^!\ov{\Dt}_* K_{\ov{X}}
\!\!\!\!\! @<BC<<\!\!\!\!\!\ov{\Dt}'_* K_{Fix(\ov{c})}, 
\endCD 
\end{equation}
where $(1)$ and $(2)$ are the left and the right vertical maps of
\form{cohmor1}, respectively, $(3)=f_Cc^!\un{Tr}_{c}$ and $(4)=\ov{c}^!\un{Tr}_{\ov{c}}$.  
In particular, the commutativity of the middle inner 
square of \form{cohmor7} follows from that of  \form{cohmor1}. 
Thus it will suffice to show the commutativity of the left and the right inner squares 
of \form{cohmor7}.
\end{Emp}

\begin{Emp}
{\bf Proof of the commutativity of the left inner square of \form{cohmor7}.}

By axiom (III) and \rl{basechange} (a), the base change morphism $BC^*_{c_1}$ decomposes as 
$BC^*_{c}\circ BC^*_{p_1}$ and  $BC^!_{c_2}$ decomposes as 
$BC^!_{p_2}\circ BC^!_{c}$. Therefore the  left inner
square of \form{cohmor7} extends to the diagram

\begin{equation} \label{Eq:cohmor8}
\CD
    f_C\C{RHom}(c^* p_1^*\C{F},c^! p_2^!\C{F}) @<t_{c^!}<<f_C c^!\C{RHom}(p_1^*\C{F},p_2^!\C{F})\\
     @V{(BC_{c}^*, BC_{c}^!)\circ r_C}VV @V{r_{X\times X}\circ BC^!_c}VV \\
\C{RHom}(\ov{c}^*f_{X\times X} p_1^*\C{F},\ov{c}^! f_{X\times X} p_2^!\C{F}) 
@<t_{\ov{c}^!}<<\ov{c}^!\C{RHom}(f_{X\times X} p_1^*\C{F},f_{X\times X} p_2^!\C{F}) \\
      @V(BC^*_{p_1},BC^!_{p_2})VV   @V(BC^*_{p_1},BC^!_{p_2})VV\\
     \C{RHom}(\ov{c}^* \ov{p}_1^*f_X\C{F},\ov{c}^! \ov{p}_2^!f_X\C{F}) 
 @<t_{\ov{c}^!}<<\ov{c}^!\C{RHom}(\ov{p}_1^*f_X\C{F},\ov{p}_2^!f_X\C{F}).
\endCD
\end{equation}
The bottom inner square of \form{cohmor8} is commutative by functoriality, 
while the top inner square is commutative by \rl{proper} (b) for $g=c$, $\C{A}=p_1^*\C{F}$
and  $\C{B}=p_2^!\C{F}$.
\end{Emp}

\begin{Emp}
{\bf Proof of the commutativity of the right inner square of \form{cohmor7}.}

The right inner square of \form{cohmor7} extends to the diagram

\begin{equation} \label{Eq:cohmor9}
\CD 
   f_C c^! \Dt_* K_{X}  @<BC<<  f_C \Dt'_* c'^! K_{X} @= f_C\Dt'_* K_{Fix(c)}\\
      @VV BC_{*\Dt}\circ BC^!_{c}V  @VVBC^!_{c'}\circ BC_{*\Dt'}V @VVBC_{*\Dt'}V\\
 \ov{c}^! \ov{\Dt}_* f_X K_{X}  @<BC<< \ov{\Dt}'_*  \ov{c}'^! f_X K_{X}       
@<BC^!_{c'}<< \ov{\Dt}'_* f_{Fix(c)} K_{Fix(c)}\\ 
@VV{\wt{\pi}_X}V  @VV\wt{\pi}_XV @VV\wt{\pi}_{Fix(c)}V\\   
\ov{c}^!\ov{\Dt}_* K_{\ov{X}}
 @<BC<<\ov{\Dt}'_* \ov{c}'^! K_{\ov{X}} @=\ov{\Dt}'_* K_{Fix(\ov{c})}. 
\endCD 
\end{equation}
Since $\Dt$ and $\ov{\Dt}$ are proper, the top left inner square of \form{cohmor9} is commutative 
by \rl{basechange} (c2) and axiom (VI). 
The bottom left inner square of \form{cohmor9} is commutative by functoriality.
The top right inner square of \form{cohmor9} is commutative by definition, while 
the bottom right inner square is commutative by  \rl{basechange} (b).
\end{Emp}

This completes the proof of \rt{cohmor} and hence also of 
Propositions \ref{P:pushf} and \ref{P:spec}.

\section{Additivity of trace maps} \label{S:additivity}

In this subsection we will prove \rp{add}. Following Pink \cite[Prop. 2.4.3]{Pi},   
we will deduce it from the additivity of filtered trace maps.
We start from preliminaries.

\subsection{Filtered derived categories} \label{SS:filtered}
\begin*
\vskip 8truept
\end*

\begin{Emp} \label{E:fcom}
{\bf Filtered complexes} (compare \cite[V, $\S$1]{Il2}). Let  $\un{A}$ be an abelian category.

(a) Denote by $C(\un{A})$ the category of complexes in $\un{A}$. 
Denote by $CF(\un{A})$ the category, whose objects are complexes $\C{A}\in C(\un{A})$
 equipped with a decreasing filtration $F^i\C{A}, i\in\B{Z}$ such that 
$F^i\C{A}=\C{A}$ for $i<<0$ and $F^i\C{A}=0$ for $i>>0$, and whose morphisms are morphisms 
in $C(\un{A})$ preserving the filtrations.

(b) For each pair $-\infty\leq a\leq b\leq+\infty$, denote by
$CF^{[a,b]}(\un{A})$ the full subcategory consisting of complexes $\C{A}$
such that $F^i\C{A}=\C{A}$  for $i\leq a$ and $F^i\C{A}=0$ for $i>b$.
Denote by $\frak{t}^{[a,b]}$ the functor  $CF(\un{A})\to CF^{[a,b]}(\un{A})$ such that
$\frak{t}^{[a,b]}\C{A}:=F^a\C{A}/F^{b+1}\C{A}$ with the induced filtration
$F^i(\frak{t}^{[a,b]}\C{A}):=F^i\C{A}/F^{b+1}\C{A}$ for all $a\leq i\leq b$.  

(c) When  $b=+\infty$ (resp.   $a=-\infty$) we will write $\geq a$ (resp. $\leq b$) 
instead of $[a,b]$. Then  $\frak{t}^{\geq a}$ (resp. $\frak{t}^{\leq b}$) is a left (resp. right) 
adjoint functor to the embedding $CF^{\geq a}(\un{A})\hra CF(\un{A})$ (resp.  
$CF^{\leq b}(\un{A})\hra CF(\un{A})$).
Also each  $\frak{t}^{[a,b]}$ can be written as a composition 
$\frak{t}^{\geq a}\circ\frak{t}^{\leq b}=\frak{t}^{\leq b}\circ\frak{t}^{\geq a}$.

(d) Denote by $\om$ the forgetful functor $CF(\un{A})\to C(\un{A})$.
Then for each $i\in \B{Z}$ the functor  
$gr^i:=\om\circ\frak{t}^{[i,i]}:CF(\un{A})\to C(\un{A})$ associates to $\C{A}$ its graded
piece $F^i\C{A}/F^{i+1}\C{A}$.

(e) Everything said above can be generalized to the case of 
bounded derived categories $D^*F(\un{A})$, where $*=+,-,b$.
\end{Emp}

\begin{Emp}
{\bf Localization.}

(a) In the notation of \re{fcom}, denote by  
$DF(\un{A})$ (resp. $DF^{[a,b]}(\un{A})$) the localization of
$CF(\un{A})$ (resp. $CF^{[a,b]}(\un{A})$) by  filtered quasi-isomorphisms, that is, by morphisms 
$f$ in $CF(\un{A})$ such that each $gr^i(f)$ is a quasi-ismorphism in $C(\un{A})$ 
(see \cite[V, 1.2]{Il2}).
Each  $DF^{[a,b]}(\un{A})$ is a full subcategory of $DF(\un{A})$ (\cite[V, Lem. 1.2.7.1]{Il2}).

(b) Functors from \re{fcom} descend to the corresponding 
functors $\frak{t}^{[a,b]}, \frak{t}^{\geq a}, \frak{t}^{\leq b}, \om, gr^i$ between derived 
categories. Moreover, the adjointness properties of 
$\frak{t}^{\geq a}$ and $\frak{t}^{\leq b}$ still hold.

(c) For every $\C{A}\in DF(\un{A})$ and
$i\in\B{Z}$, there exists a distinguished triangle 
$\frak{t}^{\geq i}\C{A}\to \C{A}\to \frak{t}^{\leq i-1}\C{A}\to$. Moreover, by the adjointness 
properties of $\frak{t}^{\geq i}$ and $\frak{t}^{\leq i-1}$, this is the unique 
distinguished triangle $\C{A}^{\geq i}\to\C{A}\to\C{A}^{\leq i-1}\to$ with
$\C{A}^{\geq i}\in DF^{\geq i}(\un{A})$ and $\C{A}^{\leq i-1}\in DF^{\leq i-1}(\un{A})$.
\end{Emp}

\begin{Emp} 
{\bf Filtered (bi)functors.} Let $\un{A}$, $\un{B}$ and $\un{C}$ be three abelian categories.

(a) We say that a triangulated functor 
$\wt{G}:DF(\un{A})\to DF(\un{B})$ (resp. $\wt{G}:DF(\un{A})^{op}\to DF(\un{B})$)
is {\em filtered}, if it satisfies
$\wt{G}(DF^{[a,b]}(\un{A}))\subset DF^{[a,b]}(\un{B})$
(resp.  $\wt{G}(DF^{[a,b]}(\un{A})^{op})\subset DF^{[-b,-a]}(\un{B})$).

(b) By a {\em filtered lift} of  a triangulated functor  
$G:D(\un{A})\to D(\un{B})$ (resp. $G:D(\un{A})^{op}\to D(\un{B})$) 
we mean a pair $(\wt{G},\phi_G)$, where $\wt{G}$ is a filtered functor as in (a),  
and $\phi_G$ is an isomorphism of functors $\om\circ\wt{G}\isom G\circ\om$.

(c) We say that a triangulated bifunctor 
$\wt{G}:DF(\un{A})\times DF(\un{B})\to  DF(\un{C})$ 
(resp. $\wt{G}:DF(\un{A})^{op}\times DF(\un{B})\to  DF(\un{C})$)
is {\em bifiltered}, if it satisfies
$\wt{G}(DF^{[a,b]}(\un{A})\times DF^{[c,d]}(\un{B}))\subset DF^{[a+c,b+d]}(\un{C})$
(resp.  $\wt{G}(DF^{[a,b]}(\un{A})^{op}\times DF^{[c,d]}(\un{B}))\subset 
DF^{[-b+c,-a+d]}(\un{C})$).

(d)  By a {\em filtered lift} of  a triangulated bifunctor  
$G:D(\un{A})\times D(\un{B})\to D(\un{C})$ (resp. 
$G:D(\un{A})^{op}\times D(\un{B})\to D(\un{C})$) 
we mean a pair $(\wt{G},\phi_G)$, where $\wt{G}$ is a filtered bifunctor as in (c),  
and $\phi_G$ is an isomorphism of bifunctors $\om\circ\wt{G}\isom G\circ(\om,\om)$.

(e) Assume that we are in the situation of (b) (resp. (d)). By a {\em filtered lift}
of a morphism  $H:G_1\to G_2$ between triangulated functors (resp. bifunctors), we mean 
a morphism  $\wt{H}:\wt{G}_1\to\wt{G}_2$ between filtered lifts of $G_1$ and $G_2$
such that isomorphisms $\phi_{G_1}$ and $\phi_{G_2}$ identify morphism 
$\om\circ\wt{H}:\om\circ\wt{G}_1\to\om\circ\wt{G}_2$ with  $H\circ\om$ (resp. $H\circ(\om,\om)$).
\end{Emp} 

The following simple lemma will play an important role afterwards.

\begin{Lem} \label{L:filtr}
(i) In the notation of (b) and (d), the isomorphism $\phi_G$ identifies the functor
$gr^i\circ\wt{G}$ with $G\circ gr^i$ (resp. $G\circ gr^{-i}$)
in the case (b) and with $\oplus_{j+k=i}G\circ(gr^j,gr^k)$
(resp.  $\oplus_{j+k=i}G\circ(gr^{-j},gr^k)$)
in the case (d).

(ii) Let $H$ and $\wt{H}$ be as in (e). Then  isomorphisms of (i) identify 
$gr^i \circ \wt{H}$ with $H\circ gr^i$ 
(resp.  $H\circ gr^{-i}$) in the case (b), and with 
 $\oplus_{j+k=i}H\circ(gr^j,gr^k)$
(resp.  $\oplus_{j+k=i}H\circ(gr^{-j},gr^k)$) in the case (d).

(iii) If $H$ is an isomorphism of (bi)functors, then $\wt{H}$ is an isomorphism as well.
\end{Lem}

\begin{proof}
As an illustration, we will give the proof in the case of bifunctors
$DF(\un{A})\times DF(\un{B})\to DF(\un{C})$. Fix $\C{A}\in DF(\un{A})$ and  
$\C{B}\in DF(\un{B})$.

(i) First we claim that $\phi_G$ defines an isomorphism
$gr^{i}\wt{G}(\frak{t}^{\leq j}\!\C{A},\!\frak{t}^{\leq k}\!\C{B})\!\isom\! 
G(gr^{j}\!\C{A},\!gr^k\!\C{B})$ for each $j,k\in\B{Z}$ with $j+k=i$.

Since $\wt{G}(\frak{t}^{\leq j}\C{A},\frak{t}^{\leq k-1}\C{B})\subset DF^{\leq i-1}(\un{C})$,
we have $gr^{i}\wt{G}(\frak{t}^{\leq j}\C{A},\frak{t}^{\leq k-1}\C{B})=0$. The distinguished triangle
$\frak{t}^{[k,k]}\C{B}\to\frak{t}^{\leq k}\C{B}\to \frak{t}^{\leq k-1}\C{B}\to$ therefore induces an isomorphism
$gr^{i}\wt{G}(\frak{t}^{\leq j}\C{A},\frak{t}^{\leq k}\C{B})\overset{\sim}{\leftarrow}
gr^{i}\wt{G}(\frak{t}^{\leq j}\C{A},\frak{t}^{[k,k]}\C{B})$. Applying similar argument to  
$\frak{t}^{\leq j}\C{A}$, we obtain an isomorphism 
$gr^{i}\wt{G}(\frak{t}^{\leq j}\C{A},\frak{t}^{\leq k}\C{B})\overset{\sim}{\leftarrow}
gr^{i}\wt{G}(\frak{t}^{[j,j]}\C{A},\frak{t}^{[k,k]}\C{B})$. 
Now the required isomorphism can be written as a composition
\[
gr^{i}\wt{G}(\frak{t}^{\leq j}\C{A},\frak{t}^{\leq k}\C{B})\isom 
gr^{i}\wt{G}(\frak{t}^{[j,j]}\C{A},\frak{t}^{[k,k]}\C{B})=
\om\wt{G}(\frak{t}^{[j,j]}\C{A},\frak{t}^{[k,k]}\C{B})\overset{\phi_G}{\lra} 
G(gr^{j}\C{A},gr^k\C{B}),
\]
where the equality follows from the facts that 
$\wt{G}(\frak{t}^{[j,j]}\C{A},\frak{t}^{[k,k]}\C{B})\subset DF^{[i,i]}(\un{C})$ and 
$gr^i=\om\circ \frak{t}^{[i,i]}$.

The isomorphisms constructed above together with canonical morphisms $\C{A}\to\frak{t}^{\leq j}\C{A}$
and $\C{B}\to\frak{t}^{\leq k}\C{B}$ give rise to a morphism
\begin{equation} \label{Eq:fil}
gr^{i}\wt{G}(\C{A},\C{B})\to 
\oplus_{j+k=i} gr^{i}\wt{G}(\frak{t}^{\leq j}\C{A},\frak{t}^{\leq k}\C{B})
\isom \oplus_{j+k=i}G(gr^{j}\C{A},gr^k\C{B}),
\end{equation} 
which we claim is an isomorphism. To show it, notice that \form{fil} is a morphism of triangulated 
bifunctors, hence it will suffice to prove it under the assumption that 
$\C{A}\in DF^{[a,a]}(\un{A})$ and $\C{B}\in DF^{[b,b]}(\un{B})$ for some $a,b\in\B{Z}$.
In this case, both sides of \form{fil} vanish  unless $a+b=i$, in 
which case the map \form{fil} coincides with the isomorphism  
$\phi_G:\om\wt{G}(\C{A},\C{B})\isom G(\om\C{A},\om\C{B})$.

(ii) It will suffice to show that for each $j,k\in\B{Z}$ with $j+k=i$, the diagram 
\[
\CD
gr^{i}\wt{G}_1(\C{A},\C{B}) @> gr^{i}\wt{H}>> gr^{i}\wt{G}_2(\C{A},\C{B})\\
@AAA @AAA\\
G_1(gr^{j}\C{A},gr^k\C{B})  @> H>> G_2(gr^{j}\C{A},gr^k\C{B}),
\endCD
\]
where the vertical maps are defined in (i), is commutative. Again we can assume 
that $\C{A}\in DF^{[a,a]}(\un{A})$ and $\C{B}\in DF^{[b,b]}(\un{B})$ 
for some $a,b\in\B{Z}$. Then all objects of the diagram vanish unless  
$a=j$ and $b=k$. In this case, the vertical maps equal $\phi^{-1}_{G_1}$ and $\phi^{-1}_{G_2}$,
respectively, so the assertion follows from the definition of $\wt{H}$.

(iii) Since  $\wt{H}$ is a morphism of triangulated bifunctors, it will suffice to show that
each $\frak{t}^{[i,i]}\wt{H}$ is an isomorphism, hence it will suffice to show that each $gr^i \wt{H}$
is an isomorphism. Thus the assertion follows from (ii).
\end{proof}

\begin{Emp} \label{E:rhom}
{\bf Filtered $RHom$.} Assume that $\un{A}$ has enough injectives. 

 By \cite[V, 1.4]{Il2}, the bifunctor
 $\RHom:D(\un{A})^0\times D^+(\un{A})\to D(Ab)$ has a filtered lift 
$\wt{\RHom}:DF(\un{A})^0\times D^+ F(\un{A})\to DF(Ab)$. 
Moreover, we have a natural isomorphism of bifunctors $\Hom(\C{A},\C{B})\isom
H^0(\om\frak{t}^{\geq 0}\wt{\RHom}(\C{A},\C{B}))$ (see \cite[V, Cor. 1.4.6]{Il2}).
Furthermore, for every  $u\in \Hom(\C{A},\C{B})=H^0(\om\frak{t}^{\geq 0}\wt{\RHom}(\C{A},\C{B}))$,
the $i$-th component of the projection of $u$ to  
$H^0(gr^{0}\wt{\RHom}(\C{A},\C{B}))\cong\oplus_i\Hom(gr^{i}\C{A},gr^i\C{B})$
equals  $gr^iu$.
\end{Emp}

\subsection{Filtered  $D^b_{ctf}(X,\La)$} \label{SS:filtrace}
\begin*
\vskip 8truept
\end*
\begin{Emp}
{\bf Filtered six operations.}

(a) For a scheme $X$, denote by $D^bF_{ctf}(X,\La)$ 
the full subcategory of $D^b F(X,\La)$ consisting of objects 
$\wt{\C{F}}$ such that each $gr^i\wt{\C{F}}$ belongs to $D^b_{ctf}(X,\La)$. 

(b) Functors $f^*, f_*, f_!, f^!$ and bifunctors  $\otimes, \C{RHom}$ on 
$D^b_{ctf}(X,\La)$ have natural filtered lifts to  $D^b F_{ctf}(X,\La)$.
To show it notice first that functors $f^*, f_*$ and bifunctors $\otimes,\C{RHom}$ 
defined in \cite[V, $\S$2]{Il2} preserve  $D^bF_{ctf}(X,\La)$ (use \rl{filtr} (i)).
Next we define the ``filtered Verdier duality'' 
$\B{D}:\C{F}\mapsto\C{RHom}(\C{F}, K_X)$, where 
$K_X$ is considered as an object of $D^bF^{[0,0]}_{ctf}(X,\La)$,
and define filtered lifts of $f_!$ and $f^!$ by formulas
$f_!:=\B{D}f_*\B{D}$ and $f^!:=\B{D}f^*\B{D}$. 
\end{Emp}

\begin{Emp} \label{E:genpr}
{\bf General principle}.
{\em All canonical (iso)morphisms of functors in 
$D^b_{ctf}(X,\La)$ have filtered analogs.} 

Indeed, all these morphisms are defined recursively 
using adjointness and inverting isomorphisms. 
However adjointness also holds in the filtered case, and 
\rl{filtr} (iii) implies that every filtered lift of an 
isomorphism in $D^{b}_{ctf}(X,\La)$ is an isomorphism in
$D^{b}F_{ctf}(X,\La)$. 

\end{Emp}

For the proof of \rp{add}, we will need the following construction.

\begin{Ex} \label{E:filex}
Let $Z\subset X$ be a closed subset, $U:=X\sm Z$, and let $i:Z\hra X$ and $j:U\hra X$ be 
the inclusion maps. Then the forgetful functor
$D^b F^{[0,1]}_{ctf}(X,\La)\to D_{ctf}^b(X,\La)$ has a natural section 
$\C{F}\mapsto\wt{\C{F}}$ such that 
$gr^0\wt{\C{F}}=\C{F}_Z:=i_* i^*\C{F}$ and
$gr^1\wt{\C{F}}=\C{F}_U:=j_! j^*\C{F}$.

To construct the section, we may assume that $\La$ is finite. Denote by $Sh(X,\La)$ the category of
\'etale sheaves of $\La$-modules on $X$. 
Then the forgetful functor $C^b F^{[0,1]}(Sh(X,\La))\to C^b(Sh(X,\La))$
has a section $\C{A}\mapsto \wt{\C{A}}:=(\C{A},F^{\cdot})$ with 
$F^1\C{A}= j_! j^*\C{A}$.
Since functors $j^*$ and $j_!$ are exact, this section descends to the required
functor  $D_{ctf}^b(X,\La)\to D^b F^{[0,1]}_{ctf}(X,\La)$.
\end{Ex}

\subsection{Additivity of filtered trace maps} \label{SS:add}
\begin*
\vskip 8truept
\end*

\begin{Emp} \label{E:filtr}
Fix a correspondence $c:C\to X\times X$, an object 
$\wt{\C{F}}$ of $D^bF_{ctf}(X,\La)$, a $c$-morphism $u\in\Hom(c_{2!}c_1^*\wt{\C{F}},\wt{\C{F}})$, 
and set $\C{F}:=\om\wt{\C{F}}$. 

Recall that $\om(c_{2!}c_1^*\wt{\C{F}})$ is naturally isomorphic to 
$c_{2!}c_1^*\C{F}$, while each $gr^i(c_{2!}c_1^*\wt{\C{F}})$ is naturally isomorphic to 
$c_{2!}c_1^*gr^i\wt{\C{F}}$ (by \rl{filtr} (i)). Hence $\om u$ and $gr^iu$ are 
naturally elements of 
$\Hom(c_{2!}c_1^*\C{F},\C{F})$ and
$\Hom(c_{2!}c_1^*gr^i\wt{\C{F}},gr^i\wt{\C{F}})$, respectively. 
\end{Emp}

The following result was formulated in \cite[4.13]{Il}.

\begin{Lem} \label{L:add}
In the notation of \re{filtr}, we have an equality 
\[
\C{Tr}_{\C{F}}(\om u)=\sum_i\C{Tr}_{gr^i\wt{\C{F}}}(gr^iu)\in H^0(Fix(c),K_{Fix(c)}).
\]
\end{Lem}  
\begin{proof}
Consider the composition 
\[
\un{\C{Tr}}'_{\C{F}}:\C{RHom}(c_{2!}c_1^*\C{F},\C{F})\isom
c_{2*}\C{RHom}(c_1^*\C{F},c_2^!\C{F})\overset{c_{2*}\un{\C{Tr}}_{\C{F}}}{\lra}
c'_{*} K_{Fix(c)},
\] 
where the first map is the 
sheafification of the adjointness $\RHom(\!\C{A},\!f^!\C{B})\!\!\isom\!\!\RHom(\!f_!\C{A},\!\C{B})$.
By the general principle of \re{genpr}, the map $\un{\C{Tr}}'_{\C{F}}$ 
lifts to a filtered map 
\[
\un{\C{Tr}}'_{\wt{\C{F}}}:\C{RHom}(c_{2!}c_1^*\wt{\C{F}},\wt{\C{F}})\to c'_{*} K_{Fix(c)},
\] 
where $c'_{*} K_{Fix(c)}$ is considered as an object of  $D^bF^{[0,0]}_{ctf}(X,\La)$. Hence 
$\un{\C{Tr}}'_{\wt{\C{F}}}$ induces a morphism
\[
gr^0(\un{\C{Tr}}'_{\wt{\C{F}}}):gr^0\C{RHom}(c_{2!}c_1^*\wt{\C{F}},\wt{\C{F}})\cong
\oplus_i \C{RHom}(c_{2!}c_1^* gr^i\wt{\C{F}},gr^i\wt{\C{F}})
\to c'_{*} K_{Fix(c)}.
\]
We claim that the restriction of $gr^0(\un{\C{Tr}}'_{\wt{\C{F}}})$
to each $\C{RHom}(c_{2!}c_1^* gr^i\wt{\C{F}},gr^i\wt{\C{F}})$ equals 
$\un{\C{Tr}}'_{gr^i\wt{\C{F}}}$. 

To show it, we note that $\un{\C{Tr}}'_{\wt{\C{F}}}$ 
decomposes as 
\[
\C{RHom}(c_{2!}c_1^*\wt{\C{F}},\wt{\C{F}})\overset{\un{\C{Tr}}''_{\wt{\C{F}}}}{\lra} 
c'_{*}c'^!(\B{D}\wt{\C{F}}\otimes\wt{\C{F}})\overset{ev}{\lra}c'_{*} K_{Fix(c)}.
\]
Since $\un{\C{Tr}}''_{\wt{\C{F}}}$ is a morphism of filtered bifunctors, 
\rl{filtr} (ii) implies that 
\[
gr^0(\un{\C{Tr}}''_{\wt{\C{F}}}):\oplus_i 
\C{RHom}(c_{2!}c_1^* gr^i\wt{\C{F}},gr^i\wt{\C{F}})\to
\oplus_i c'_{*}c'^!(\B{D}(gr^i\wt{\C{F}})\otimes gr^i\wt{\C{F}})
\] 
equals $\oplus_i\un{\C{Tr}}''_{gr^i\wt{\C{F}}}$. 
Thus it remains to show that the restriction of   
\[
gr^0(ev_{\wt{\C{F}}}):gr^0(\B{D}\wt{\C{F}}\otimes\wt{\C{F}})\cong \oplus_i
(\B{D}(gr^i\wt{\C{F}})\otimes gr^i\wt{\C{F}})\to K_X
\]
to $\B{D}(gr^i\wt{\C{F}})\otimes gr^i\wt{\C{F}}$ equals $ev_{gr^i \wt{\C{F}}}$.

Since under the isomorphism $\C{RHom}(\B{D}\wt{\C{F}}\otimes\wt{\C{F}}, K_X)\isom
\C{RHom}(\wt{\C{F}},\wt{\C{F}})$, the map
$ev_{\wt{\C{F}}}$ corresponds to $\Id_{\wt{\C{F}}}$, the assertion is a reformulation of the 
obvious equality $gr^i\Id_{\C{F}}=\Id_{gr^i\wt{\C{F}}}$.

Now the assertion of the lemma is easy. Indeed, recall that $\C{Tr}_{\C{F}}(\om u)$
is the image of $u$ under the map 
\[
H^0(\om\frak{t}^{\geq 0}\un{\C{Tr}}'_{\wt{\C{F}}}):\Hom(c_{2!}c_1^*\wt{\C{F}},\wt{\C{F}})= 
H^0(X,\om\frak{t}^{\geq 0}\C{RHom}(c_{2!}c_1^*\wt{\C{F}},\wt{\C{F}}))\to H^0(X, c'_* K_{Fix(c)}).
\]
However, $H^0(\om\frak{t}^{\geq 0}\un{\C{Tr}}'_{\wt{\C{F}}})$ factors through
\[
H^0(X,gr^0(\un{\C{Tr}}'_{\wt{\C{F}}})):H^0(X,gr^0\C{RHom}(c_{2!}c_1^*\wt{\C{F}},\wt{\C{F}})) 
\to H^0(X, c'_* K_{Fix(c)}),
\]
thus the assertion follows from the last observation of \re{rhom} and the  
claim proven above.
\end{proof}

\begin{Emp}
{\bf Proof of \rp{add}.}

Consider the object $\wt{\C{F}}\in D^b F^{[0,1]}(X,\La)$ constructed in \re{filex}.
To prove the proposition, we will show that  $u\in \Hom(c_{2!}c_1^*\C{F},\C{F})$ lifts to a
unique element $\wt{u}\in \Hom(c_{2!}c_1^*\wt{\C{F}}, \wt{\C{F}})$
and that $gr^0\wt{u}=[i_Z]_!(u|_Z)$ and $gr^1\wt{u}=[j_U]_!(u|_U)$.
Then the assertion will follows from \rl{add}.

First we will show that the forgetful map 
$\om:\Hom(c_{2!}c_1^*\wt{\C{F}}, \wt{\C{F}})\to 
\Hom(c_{2!}c_1^*\C{F},\C{F})$ is an isomorphism.
Since $\Hom(c_{2!}c_1^*\wt{\C{F}}, \wt{\C{F}})= H^0(\om\frak{t}^{\geq 0}
\RHom(c_{2!}c_1^*\wt{\C{F}},\wt{\C{F}}))$, while 
$\Hom(c_{2!}c_1^*\C{F},\C{F})=H^0(\om\RHom(c_{2!}c_1^*\wt{\C{F}},\wt{\C{F}}))$,
it will suffice to check the vanishing of  
$\frak{t}^{<0}\RHom(c_{2!}c_1^*\wt{\C{F}},\wt{\C{F}})$. 
As $\wt{\C{F}}\in D^b F^{[0,1]}(X,\La)$, one sees that 
$gr^{-i}\RHom(c_{2!}c_1^*\wt{\C{F}},\wt{\C{F}})=0$ for all $i>1$, while  
$
gr^{-1}\RHom(c_{2!}c_1^*\wt{\C{F}},\wt{\C{F}})\cong\RHom(c_{2!}c_1^*\C{F}_U,\C{F}_Z)$.

Since $Z$ is $c$-invariant, the subset $c_2 (c_1^{-1}(U))$ is contained in $U$, 
hence $c_{2!}c_1^*\C{F}_U$ is supported on $c_2(c_1^{-1}(U))\subset U$, thus 
$i^* c_{2!}c_1^*\C{F}_U=0$. It follows that 
\begin{equation} \label{Eq:van}
\RHom(c_{2!}c_1^*\C{F}_U,\C{F}_Z)\cong\RHom(i^* c_{2!}c_1^*\C{F}_U,i^*\C{F})=0,
\end{equation}
implying the vanishing of $\frak{t}^{<0}\RHom(c_{2!}c_1^*\wt{\C{F}},\wt{\C{F}})$. 

For the second assertion, note that the pair $(gr^0\wt{u},gr^1\wt{u})$ can be characterized
as the unique (by \form{van}) element of
$\Hom(c_{2!}c_1^*\C{F}_Z,\C{F}_Z)\oplus\Hom(c_{2!}c_1^*\C{F}_U,\C{F}_U)$, 
whose image in 
$\Hom(c_{2!}c_1^*\C{F}_U,\C{F})\oplus\Hom(c_{2!}c_1^*\C{F},\C{F}_Z)$ equals 
the image of $u$. Since the pair $([i_Z]_!(u|_Z),[j_U]_!(u|_U))$ also satisfies this property, 
the assertion follows.
\end{Emp}

\appendix
\section{Compatibility of local terms}

The goal of the appendix is to compare the trace map, defined in \re{locterms},
with the pairing \cite[(4.2.5)]{Il} of Illusie, which was used by Pink and Fujiwara.
First we need to introduce some notation. 

\begin{Con'} \label{C:pairing}
(a) Suppose we are given a Cartesian square of schemes 
$$
\CD 
        C   @>>{b'}>                 A\\
        @V{a'}VV                        @V{a}VV\\
        B @>>{b}>                         Y,
\endCD 
$$
and two objects $\C{P},\C{Q}\in D^b_{ctf}(Y,\La)$ equipped with a morphism 
$\Phi:\C{P}\otimes\C{Q}\to K_Y$.
Denote by $c$ the composition map $b\circ a'=a\circ{b}':C\to Y$.
Then $\Phi$ gives rise to a pairing
\[
\un{\lan\cdot,\cdot\ran}:b'^*a^!\C{P}\otimes a'^*b^!\C{Q}\to c^!(\C{P}\otimes\C{Q})
\overset{\Phi}{\lra}c^! K_{Y}=K_C,
\]
where the first map is the composition
\[
b'^*a^!\C{P}\otimes a'^*b^!\C{Q}\overset{BC}{\lra}a'^!b^*\C{P}\otimes a'^*b^!\C{Q}
\overset{t_{a'^!}}{\lra}
a'^!(b^*\C{P}\otimes b^!\C{Q})\overset{t_{b'^!}}{\lra} a'^!b^!(\C{P}\otimes\C{Q}).
\]
Then $\Phi$ gives rise to a pairing 
\[
\lan\cdot,\cdot\ran:H^0(A,a^!\C{P})\otimes H^0(B,b^!\C{Q})\overset{(b'^*,a'^*)}{\lra}
H^0(C,b'^* a^!\C{P})\otimes H^0(C,a'^* b^!\C{Q})\overset{\otimes}{\lra}
\]\[
\overset{\otimes}{\lra}H^0(C,b'^*a^!\C{P}\otimes a'^*b^!\C{Q})
\overset{\un{\lan\cdot,\cdot\ran}}{\lra} H^0(C, K_C).
\]

(b) Assume that in the notation of (a) we have $Y=X_1\times X_2$,
$\C{P}=\B{D}\C{F}_1\pp \C{F}_2$, $\C{Q}=\C{F}_1\pp\B{D}\C{F}_2$ for some
 $\C{F}_1\in D_{ctf}^b(X_1,\La)$, $\C{F}_2\in D_{ctf}^b(X_2,\La)$, and $\Phi$
is the composition
\[
(\B{D}\C{F}_1\pp \C{F}_2)\otimes (\C{F}_1\pp\B{D}\C{F}_2)\cong
(\B{D}\C{F}_1\otimes \C{F}_1)\pp (\C{F}_2\otimes \B{D}\C{F}_2)
\overset{ev\pp ev}{\lra} K_{X_1}\pp K_{X_2}\overset{\diam}{\lra}
K_{X_1\times X_2},
\]
where $\diam$ is the isomorphism \form{*}.  
Put $a_i:=p_i\circ a:A\to X_i$ and $b_i:=p_i\circ b:B\to X_i$ for $i=1,2$. 
Then we have a natural isomorphism
\[
H^0(A,a^!\C{P})\isom H^0(A,\C{RHom}(a_1^*\C{F}_1,a_2^!\C{F}_2))=
\Hom(a_1^*\C{F}_1,a_2^!\C{F}_2)\isom\Hom(a_{2!}a_1^*\C{F}_1,\C{F}_2),
\]
induced by the isomorphism 
\[
a^!\C{P}\overset{\form{*}}{\lra}a^!\C{RHom}(p_1^*\C{F}_1,p_2^!\C{F}_2)
\overset{t_{f^!}}{\lra}\C{RHom}(a_1^*\C{F}_1,a_2^!\C{F}_2),
\]
and similarly an isomorphism 
$H^0(B,b^!\C{Q})\isom \Hom(b_{1!}b_2^*\C{F}_2,\C{F}_1)$, induced by a map 
\begin{equation} \label{Eq:is}
b^!\C{Q}\isom \C{RHom}(b_2^*\C{F}_2,b_1^!\C{F}_1).
\end{equation} 
Therefore the pairing $\lan\cdot,\cdot\ran$ from (a)  specializes to the pairing 
\[
\lan\cdot,\cdot\ran:\Hom(a_{2!}a_1^*\C{F}_1,\C{F}_2)\otimes\Hom(b_{1!}b_2^*\C{F}_2,\C{F}_1)
\to H^0(C, K_C).
\]

(c)  Assume that in the notation of (b) we have $B=X_1=X_2=X$, $\C{F}_1=\C{F}_2=\C{F}$,  
and $b:X\to X\times X$ is the diagonal morphism. In this case, $C=Fix(a)$, $b_1=b_2=\Id_X$, 
hence $\lan\cdot,\cdot\ran$ specializes to the pairing 
\[
\lan\cdot,\cdot\ran:\Hom(a_{2!}a_1^*\C{F},\C{F})\otimes\Hom(\C{F},\C{F})\to 
H^0(Fix(a), K_{Fix(a)}).
\]
\end{Con'}

\begin{Lem'} \label{L:Il}
In the notation of (c), for each $u\in \Hom(a_{2!}a_1^*\C{F},\C{F})$, we have
\[
\C{Tr}_a(u)=\lan u,\Id_{\C{F}}\ran.
\]
\end{Lem'} 
\begin{proof}
Note that $\C{Tr}_a:\Hom(a_{2!}a_1^*\C{F},\C{F})\to H^0(Fix(a), K_{Fix(a)})$
decomposes as
\[
H^0(\!A,\!\C{RHom}(a_1^*\C{F},\!a_2^!\C{F}))\!\overset{\Dt'^*}{\lra}\!
H^0(Fix(a),\!\Dt'^*\C{RHom}(\!a_1^*\C{F}\!,a_2^!\C{F}))\!\to\!H^0(Fix(a), K_{Fix(a)}),
\]
where the last map is induced by the composition
\[
\Dt'^*\C{RHom}(a_1^*\C{F},a_2^!\C{F})\isom \Dt'^* a^!(\B{D}\C{F}\pp\C{F})
\overset{BC}{\lra} a'^!\Dt^*(\B{D}\C{F}\pp\C{F})= a'^!(\B{D}\C{F}\otimes\C{F})
\overset{ev}{\lra} a'^! K_X.
\]
Therefore the assertion reduces to the commutativity of the following 
diagram
\begin{equation} \label{Eq:il1}
\CD
a'^!\Dt^*(\B{D}\C{F}\pp\C{F})\otimes a'^*\La_X @>t_{a'^!}>>  
a'^!(\Dt^*(\B{D}\C{F}\pp\C{F})\otimes\La_X) @>ev>> a'^! K_X\\
@VVV @VVV @|\\
 a'^!\Dt^*(\B{D}\C{F}\!\pp\!\C{F})\!\otimes\! a'^*\Dt^!(\C{F}\!\pp\!\B{D}\C{F})\!\! @>t_{a'^!}>>  
\!\!a'^!(\Dt^*(\B{D}\C{F}\!\pp\!\C{F})\!\otimes\! \Dt^!(\C{F}\!\pp\!\B{D}\C{F})) \!\!
@>\Phi\circ t_{\Dt^!}>>\!\! c^! K_{X\times X},
\endCD
\end{equation}
where both vertical arrows are induced by the composition
\[
\La_X\overset{\Id_{\C{F}}}{\lra}\C{RHom}(\C{F},\C{F})\overset{\form{is}^{-1}}{\lra} 
\Dt^!(\C{F}\pp\B{D}\C{F}).
\]
Since the left inner square of \form{il1} is commutative by  functoriality, 
it remains to check the commutativity of the right one. For this we can assume that 
$a'=\Id_A$, in which case the corresponding diagram extends to the diagram 
\begin{equation} \label{Eq:il2}
\CD
\B{D}\C{F}\!\otimes\!\C{F}\!\otimes\!\La_X @= \B{D}\C{F}\!\otimes\!\C{F} @>ev>> K_X\\
@V\Id_{\C{F}}VV @| @|\\
\B{D}\C{F}\!\otimes\!\C{F}\!\otimes \!\C{RHom}(\C{F},\C{F})\!  @>\Id_{\B{D}\C{F}}\otimes ev>> 
\B{D}\C{F}\!\otimes\!\C{F} @>ev>> K_X\\
@A\form{is}AA @A\diam'\circ(\Id\pp ev)AA @A\diam'AA\\
\Dt^*(\B{D}\C{F}\!\pp\!\C{F})\!\otimes\! \Dt^!(\C{F}\!\pp\!\B{D}\C{F})\!\! @>t_{\Dt^!}>>\!\!
\Dt^!((\B{D}\C{F}\!\otimes\!\C{F})\!\pp\!(\C{F}\!\otimes\!\B{D}\C{F}))\!\! @>ev\pp ev>>\!\! \Dt^!
(K_{X}\!\pp\! K_{X}),
\endCD
\end{equation}
where we denote by $\diam'$ the isomorphism $\Dt^!(\C{A}\pp K_X)\isom \Dt^!p_1^!\C{A}=\C{A}$,
induced by the isomorphism $\diam:\C{A}\pp K_X\isom p_1^!\C{A}$ from \form{*}.

We claim that all inner squares of \form{il2} are commutative.
Indeed, the assertion for the top right inner square is clear, the bottom right 
inner square is commutative by functoriality, while the assertion  
for the top left inner square follows from the fact that the composition
$\C{A}=\C{A}\otimes\La_X\overset{\Id_{\C{A}}}{\lra}\C{A}\otimes\C{RHom}(\C{A},\C{A})
\overset{ev}{\lra}\C{A}$ is the identity. 

Finally, the (rotated) bottom left inner square of \form{il2} extends to the diagram
\begin{equation} \label{Eq:il3}
\CD
\Dt^*(\B{D}\C{F}\pp\C{F})\otimes \Dt^!(\C{F}\pp\B{D}\C{F}) @>\form{is}>>
\B{D}\C{F}\otimes\C{F}\otimes\C{RHom}(\C{F},\C{F})\\
@Vt_{\Dt^!}VV   @V\Id_{\B{D}\C{F}}\otimes evVV\\ 
\Dt^*p_1^*\B{D}\C{F}\otimes \Dt^!(\C{F}\pp(\C{F}\otimes\B{D}\C{F}))  
@>\Id_{\B{D}\C{F}}\otimes(\diam'\circ(\Id\pp ev))>> \B{D}\C{F}\otimes\C{F}\\
 @Vt_{\Dt^!}VV @|\\
\Dt^!((\B{D}\C{F}\otimes\C{F})\pp(\C{F}\otimes\B{D}\C{F})) 
@>\diam'\circ(\Id\pp ev)>> \B{D}\C{F}\otimes\C{F}  
\endCD
\end{equation}
Since $\Dt^!(\C{F}\pp\B{D}\C{F})\overset{\form{is}}{\lra}\C{RHom}(\C{F},\C{F})$ 
is adjoint to the composition
\[
\Dt^* p_2^*\C{F}\otimes \Dt^!(\C{F}\pp\B{D}\C{F})\overset{t_{\Dt^!}}{\lra}
\Dt^!(\C{F}\pp(\C{F}\otimes\B{D}\C{F}))\overset{ev}{\lra}\Dt^!(\C{F}\pp K_X)
\overset{\diam'}{\lra}\C{F},
\]
the top inner square of \form{il3} is commutative. 
Next since $\Id_{\B{D}\C{F}\otimes\C{F}}$ decomposes as 
$t_{p_1^!}\circ t_{\Dt^!}:\Dt^!p_1^!\B{D}\C{F}\otimes\Dt^*p_1^*\C{F}\to
\Dt^!p_1^!(\B{D}\C{F}\otimes\C{F})$, while $\diam'$ decomposed as a composition 
$\Dt^!(p_1^*\C{A}\otimes p_2^* K_X)\overset{BC}{\lra}
\Dt^!(p_1^*\C{A}\otimes p_1^!\La_X)\overset{t_{p_1^!}}{\lra}\Dt^!p_1^!\C{A}$,
the commutativity the bottom inner square of \form{il3} reduces to that of 
\begin{equation} \label{Eq:il4}
\CD
\Dt^*p_1^*\B{D}\C{F}\otimes \Dt^!(p_1^*\C{F}\otimes p_1^!\La_X)  
@>\Id_{\B{D}\C{F}}\otimes t_{p_1^!}>> \Dt^*p_1^*\B{D}\C{F}\otimes\Dt^!p_1^!\C{F}\\
 @Vt_{\Dt^!}VV  @Vt_{p_1^!}\circ t_{\Dt^!}VV\\
\Dt^!(p_1^*(\B{D}\C{F}\otimes\C{F})\otimes p_1^!\La_X)
@> t_{p_1^!}>> \Dt^!p_1^!(\B{D}\C{F}\otimes\C{F}).  
\endCD
\end{equation}
Now the assertion follows from the associativity of $t_{p_1^!}$.
\end{proof}

\begin{Cor'} \label{C:locterms}
The local terms defined in \re{locterms} (b) coincide with those 
of Pink (\cite[Def. 2.1.10]{Pi}) and  Fujiwara (\cite[p. 495]{Fu}).
\end{Cor'}

\end{document}